\documentclass[reqno]{amsart}

\pdfoutput=1

\usepackage{graphicx}
\usepackage{amsfonts}
\usepackage{amssymb}
\usepackage{enumitem}
\usepackage{baskervillef}
\usepackage{pdfpages}

\newtheorem*{theorem*}{Theorem}
 
\newtheorem*{lemma*}{Lemma}

\newtheorem*{corollary*}{Corollary}

\newtheorem{innercustomgeneric}{\customgenericname}
\providecommand{\customgenericname}{}
\newcommand{\newcustomtheorem}[2]{%
  \newenvironment{#1}[1]
  {%
   \renewcommand\customgenericname{#2}%
   \renewcommand\theinnercustomgeneric{##1}%
   \innercustomgeneric
  }
  {\endinnercustomgeneric}
}
\newcustomtheorem{customthm}{Theorem}
\newcustomtheorem{customlemma}{Lemma}

\newtheorem*{mainlemma*}{Main Lemma}

\newtheorem*{adianRabin*}{The Adian-Rabin Theorem}

\theoremstyle{definition}
\newcustomtheorem{customremark}{Remark}

\newtheorem*{definition*}{Definition}

\newcommand{\fA}{\mathfrak{A}}
\newcommand{\fAqAB}{\mathfrak{A}_{q,A,B}}

\newcommand{\fP}{\mathfrak{P}}

\begin{document}
\thispagestyle{empty}
{
\centering
	
	~
	
	{\scshape\LARGE S. I. Adian's Papers on the \\ \Huge Undecidability of Algorithmic \\ Problems in Group Theory \par}
	\vspace{15pt}
	{\scshape\LARGE or \\ \vspace{15pt}\LARGE The Adian-Rabin Theorem\par}
	\vfill
    {\Large \textsc{Four papers by S. I. Adian originally published between 1955--1958}\par}
    {\Large \textsc{supplemented by two papers by A. A. Markov published in 1951}\par}
	{\Large \textsc{translated into English from the original Russian}\par}
	{\Large \textsc{by} \\ \Large \textsc{C.-F. Nyberg-Brodda}\par}
	\vspace{2.0cm}
	{\Large 2022} \\

}
\clearpage 

\thispagestyle{empty}

\begin{center}
{\large\textsc{Translator's Preface}}
\end{center}

\vspace{0.5cm}

\noindent The present work consists of an English translation of four remarkable articles, originally written in Russian, by Sergei Ivanovich Adian (1931--2020), supplemented by the translation of two articles by Andrei Andreevich Markov (1903--1979). The bibliographic details for the six articles are as follows; the page numbers stated next to the naming of each article is the page number for the translated article in the present work. \\

\

\begin{center}
\makebox[\linewidth]{\hfill \textbf{[Adi55]} \hfill \llap{(pp. 9--11)}}
\end{center}
S. I. Adian, \textit{Algorithmic unsolvability of problems of recognition of certain properties of groups}, Doklady Akademii Nauk SSSR (N. S.), \textbf{103}:4 (1955), 533--535. \\

\begin{center}
\makebox[\linewidth]{\hfill \textbf{[Adi57a]} \hfill \llap{(pp. 12--75)}}
\end{center}
S. I. Adian, \textit{Unsolvability of some algorithmic problems in the theory of groups}, \\ Proceedings of the Moscow Mathematical Society, \textbf{6} (1957), 231--298. \\

\begin{center}
\makebox[\linewidth]{\hfill \textbf{[Adi57b]} \hfill \llap{(pp. 76--79)}}
\end{center}
S. I. Adian, \textit{Finitely presented groups and algorithms}, \\ Doklady Akademii Nauk SSSR (N. S.), \textbf{117}:1 (1957), 9--12. \\

\begin{center}
\makebox[\linewidth]{\hfill \textbf{[Adi58]} \hfill \llap{(pp. 80--83)}}
\end{center}
S. I. Adian, \textit{On algorithmic problems in effectively complete classes of groups}, \\ Doklady Akademii Nauk SSSR (N. S.), \textbf{123}:1 (1958), 13--16. \\

\begin{center}
\makebox[\linewidth]{\hfill \textbf{[Mar51a]} \hfill \llap{(pp. 84--85)}}
\end{center}
A. A. Markov, \textit{Impossibility of some algorithms in the theory of associative systems}, \\ Doklady Akademii Nauk SSSR, \textbf{77}:1 (1951), 19--20. \\

\begin{center}
\makebox[\linewidth]{\hfill \textbf{[Mar51b]} \hfill \llap{(pp. 86--89)}}
\end{center}
A. A. Markov, \textit{Impossibility of algorithms for recognising some properties of associative systems}, \\ Doklady Akademii Nauk SSSR, \textbf{77}:6 (1951), 953--956. \\

\

\noindent The first four articles, when combined (see discussion below), yield a proof of the celebrated Adian-Rabin Theorem. This theorem is one of the most striking and beautiful theorems in combinatorial group theory, and represents a simultaneous development of both mathematical logic and group theory. Loosely speaking, the theorem states: \textit{from a finite presentation of a group, one can deduce almost nothing about the properties of the group it presents.} At first glance, to the modern reader this may seem trivial (if they have been raised to treat the existence of algorithmically undecidable properties of groups as obvious) or incredibly discouraging (if they are of an easily discouraged disposition). However, beneath this loose surface lies a striking depth, linking together algorithmic problems in mathematical logic and (semi)group theory. 

The aforementioned discouraged reader may think, upon reading the above, that no non-trivial property can be decided about a group, when given nothing more than its presentation. But this is not so. For example, it is decidable whether a finitely presented group coincides with its commutator subgroup or not. Part of the beauty of the Theorem is the isolating and formalising of the particular types of properties which are tractable to prove to be undecidable to recognise in general.

In modern terminology, and recognisable to anyone with some experience in combinatorial group theory, the Theorem is usually given as follows. Let $P$ be an abstract property of groups, i.e. one which is preserved under isomorphism (in the translations, such properties are called \textit{invariant}). Suppose the following two conditions hold:
\begin{enumerate}
\item There exists a finitely presented group $G_+$ which has property $P$.
\item There exists a finitely presented group $G_-$ which does not embed as a subgroup of any finitely presented group which has property $P$. 
\end{enumerate}
Then we say that $P$ is a \textit{Markov property}\footnote{The terminology ``Markov property'' was first introduced by William `Bill' Boone (1920--1983) in his MR review of M. Rabin's 1958 article (see below).}. Such properties, in the case of monoids (or \textit{associative systems}) rather than groups, were introduced by Markov in [Mar51b]. The Adian-Rabin theorem is stated in terms of these Markov properties; for this reason, the translation of [Mar51b] is also included in this work. As [Mar51b] is a strengthening of results already obtained in [Mar51a], similar to how [Adi57b] strengthens the results of [Adi57a] (see discussion below), I have decided to translate this, too. Some details on the two papers are provided at the end of this preface. 

The following is a select list of Markov properties, all of which are easily verified:
\begin{enumerate}[label=(\roman*)]
\item being a finite group,
\item being a finite group of fixed order $n$,
\item being the trivial group,
\item being an abelian group,
\item being a free group,
\item being a nilpotent group,
\item being a $k$-step nilpotent group for fixed $k$,
\item being a group with decidable word problem, 
\item being isomorphic to a fixed f.p. group $G$ with decidable word problem, \\
etc. (ad nauseam)
\end{enumerate}
The property ``being isomorphic to a fixed finitely presented group $G$'' is not always a Markov property. This can fail property (2), as e.g. by Higman's embedding theorem there exists a fixed finitely presented group $G_0$ which contains a copy of every finitely presented group. Thus ``being isomorphic to $G_0$'' is not a Markov property. 

Having given the definition of Markov properties, we can now state:

\begin{theorem*}[Adian--Rabin, 1958]
Let $P$ be a Markov property. Then there does not exist an algorithm which takes as input a finite group presentation, and outputs whether or not the group defined by this presentation has property $P$.
\end{theorem*}

The impact of the Adian--Rabin theorem on the development of combinatorial group theory can equally little be overstated as it can be comprehensively treated here. For this reason, we shall make no direct attempt to do either. We content ourselves with remarking that although one occasionally reads that no substantial extension of the theorem has been proved, demonstrating its essential completeness, this is not quite accurate. Indeed, part of Markov's results for monoids were not treated by either Adian or Rabin, and was only extended to groups at a later date. We refer the reader to the end of our discussion on Markov's theorem in this preface.

\begin{center}
\textbf{The proof of the Adian--Rabin Theorem}
\end{center}

\noindent It is very easy to state in completeness when (and where) Michael Rabin (b. 1931) proved the theorem partly named after him: this was done in a 1958 article (\textit{Ann. of Math. (2)}, \textbf{67}, (1958), 172--194). At the time, Rabin was a Ph.D. student under Alonzo Church (1903--1995) at Princeton, and the proof constituted his Ph.D. thesis, defended in 1957.\footnote{Adian and Rabin worked entirely independent of one another. In spite of this, each of the two authors acknowledges the other's contributions (Adian does this in [Adi58]).} An excellent modern reference for this proof can be found in G. Baumslag's book \textit{Topics in combinatorial group theory} (Birkhäuser, 1993), which gives a lucid two-page proof of the result (pp.~112--113). The underlying proof, though streamlined, is, however, not fundamentally different from the original proof(s).

By contrast, it is more difficult to state when (and where) S. I. Adian proved the theorem. Part of the \textit{raison d'\^etre} of the present translation was precisely this difficulty and my attempts to try and understand this state of affairs. Adian was a Ph.D. student working under the supervision of Petr Sergeevich Novikov (1901--1975); the proof of the results announced in [Adi55] constituted his Ph.D. thesis\footnote{A. S. Esenin-Volpin, one of Adian's thesis examiners, famously travelled to P. S. Novikov's dacha to try and convince him that the work in fact merited a D.Sc. degree (habilitation). Novikov agreed, but said that he did not doubt that Adian would produce even more remarkable results meriting a D.Sc. degree.}, and were later published in [Adi57a] with all details given. However, the theorem announced in [Adi55] (and proved in [Adi57a]) is not quite the Adian--Rabin theorem. 

We present an overview of the results stated in [Adi55] and proved in [Adi57a]. Let $P$ be an abstract property of (finitely presented) groups. Consider the conditions:
\begin{enumerate}
\item There exists some finitely presented group $G_+$ which has property $P$, and some finitely presented group $G_-$ which does not have property $P$.
\item $P$ is inherited by subgroups.
\item[(2$^\prime$)] any group with $P$ has decidable word problem.
\end{enumerate} 
Suppose (1) together with one of the conditions (2) or (2$^\prime$) hold for $P$. Then $P$ will be said to be a \textit{pseudo-Markov property}\footnote{This terminology does not appear anywhere else, but is used here for brevity. Adian himself calls a property satisfying condition (2) \textit{hereditary}, and a property satisfying condition (2$^\prime$) \textit{special}.}. Any pseudo-Markov property is clearly a Markov property. All of the properties (i)--(ix) listed earlier are, in fact, pseudo-Markov properties. Finding an example of a Markov property which is not pseudo-Markov is somewhat tricky. One example is: being solvable with derived length $k$ for a fixed $k \geq 3$. Here we use the result by Kharlampovich that groups with this property can have undecidable word problem (\textit{Izv. Akad. Nauk SSSR Ser. Mat.}, \textbf{45}:4, (1981), pp. 852--873).

The main theorem(s) of [Adi57a] -- first announced in [Adi55] -- can then be summarised in the following weaker form of the Adian--Rabin theorem:

\begin{theorem*}[Adian, 1955]
Let $P$ be a pseudo-Markov property. Then there does not exist an algorithm which takes as input a finite group presentation, and outputs whether or not the group defined by this presentation has the $P$.
\end{theorem*}

In particular, a very large part of the consequences of the Adian--Rabin theorem are already present in [Adi55, Adi57a], but the full theorem is not; however, any Markov property which implies having decidable word problem is a pseudo-Markov property, so many common Markov properties were already covered in [Adi55, Adi57a].

We give an outline of the proof of the main theorem of [Adi57a], which is later modified in [Adi57b] and [Adi58] to give the proof of the full Adian--Rabin theorem. The details of the proof are formidable; even minute details are treated with care. We add some notes on the modern terminology, to facilitate reading.

The main idea is as follows. Begin with the centrally-symmetric group $\fA_P$, constructed by Novikov, which is a finitely presented group with undecidable word problem. Chapter~I proves some auxiliary lemmas for this group. In modern terms, $\fA_P$ is exhibited as an HNN-extension of a group $\overline{\fA}_P^{(p)}$, with stable letter $p$. These auxiliary lemmas are not far removed from applications of Britton's lemma (cf. e.g. Lemma~6). The phrasing used in the translation in place of ``HNN-extension'' is that the letter $p$ is a \textit{pivot} letter which, together with is inverse, forms a \textit{regular system of filtering letters}. The \textit{base group} is used in the same sense as for HNN-extensions.

Next, in Chapter~II, from the group $\fA_P$ and every pair of words $A$ and $B$ over the alphabet of $\fA_P$ having a particular form, a group $\fAqAB$ is constructed, which is central to the entire paper. The group $\fAqAB$ is first given by a set of defining relations. Its structure is then dissected by constructed by producing a sequence of groups $\fA_{q_{i}, A, B}$ for $1 \leq i \leq n$,  where $\overline{\fA}^{(q_{1})}_{q_{1},A,B} = \fA_P$, $\overline{\fA}^{(q_{i})}_{q_{i},A,B} = \fA_{q_{i-1}, A, B}$, and $\overline{\fA}^{(q_{n})}_{q,A,B}$ is $\fA_{q_{n-1}, A, B}$ (see the very beginning of \S1, Chapter~II). This construction, which looks somewhat frightening, is really rather simple -- in modern terminology, $\fAqAB$ is essentially an iterated HNN-extension, where the notation $\overline{\fA}^{(q_{i})}_{q_{i},A,B}$ means ``the base group of the HNN-extension $\fA_{q_i,A,B}$ with stable letter $q_i$''. Thus, the frightening line above merely says that $\fAqAB$ can be obtained by starting with $\fA_P$, and then successively taking HNN-extensions (identifying certain subgroups) with stable letter $q_1$, then $q_2$, then $q_3$, etc. until after $n$ HNN-extensions we find $\fAqAB$ (here $n$ is the number of letters in the alphabet of $\fA_P$). This property of $\fAqAB$ is rather reminiscent of Rabin's construction. 

At this point, [Adi57a] moves on to proving properties about $\fAqAB$. Here lies the main technical challenge of the paper. The overall goal is to prove Lemma~4, \S2, Chapter~II. This says that $\fAqAB$ is a ``test group'', in the sense that: if $A=B$ in $\fA_P$, then $\fAqAB$ is the trivial group, and if $A \neq B$ in $\fA_P$, then Novikov's group $\fA_P$ embeds in $\fAqAB$. This yields all the undecidability results of Chapter~III, by some simple free products and an application of Grushko's theorem. However, to prove Lemma~4, we must first prove the ``Main Lemma''. This proof is a behemoth, spanning 23 pages in the translation, and with formidable case-by-case analysis. The Main Lemma has four parts \textbf{I}--\textbf{IV}, of which only \textbf{IV} is used in the subsequent proof of Lemma~4. However, the proof of \textbf{IV} requires \textbf{I}--\textbf{III}. All four parts are proved simultaneously by induction on a parameter $\lambda$. The inductive step for each part is itself proved by induction on a parameter $m$. Each such inductive step is subdivided into different cases, which divide into further subcases, etc. For example, in the inductive (on $m$) step as part of the inductive (on $\lambda$) step for \textbf{II}, we have a half-page proof dealing with subcase $\alpha)$ of subcase 2) of case 2$^\circ$. One needs no technical details of Chapter~II to understand Chapter~III and its undecidability results -- one needs only use Lemma~4, \S2, Chapter~II. A simplified proof of Lemma~4 (avoiding the Main Lemma) via Gröbner--Shirsov bases has been given by Bokut \& Chainikov (\textit{Discrete Math.}, \textbf{308}:21, (2008), pp. 4916--4930).

With an overview of the proof of the results of [Adi57a], it may be somewhat surprising that the subsequent articles [Adi57b] and [Adi58] contain, by comparison, almost no difficult combinatorial arguments whatsoever. These short papers are readable on their own, and summarise themselves, so we direct the interested reader to simply read the papers. However, it is worth mentioning that it is in [Adi57b] where we first see Adian make reference to Markov properties (without using the name), and in which the statement of the Adian--Rabin theorem appears. The proof, only outlined in [Adi57a], goes via constructing a group $F_{qA}$, which serves an entirely analogous purpose to the groups $\fAqAB$ above, but with no dependency on the precise structure of Novikov's centrally-symmetric group $\fA_0$ (this can be compared with Rabin's proof, which starts with an arbitrary finitely presented group with undecidable word problem). The key idea is then the same; $F_{qA}$ is a test group for triviality and, by appropriate free products and Grushko's theorem, yields a test group for any Markov property $\alpha$. One still needs an analogue of the Main Lemma from [Adi57a] -- the way to obtain this is sketched in [Adi58], by showing which lemmas in [Adi57a] regarding $\fA_0$ should be replaced by what lemmas for $F_0$ to yield the analogous Main Lemma. This yields a proof of the Adian--Rabin theorem. Finally, in [Adi57b] and [Adi58] there are also some results on ``effectively complete'' classes of groups. This is easy to read, needing no summary. 

Summarising, in [Adi55] a weaker form of the Adian--Rabin theorem is stated, and is proved in [Adi57a] with immense combinatorial details. In [Adi57b], the full Adian--Rabin theorem is stated, and a proof is sketched. In [Adi58], most details are given of this proof (but with some dependency on technical results in [Adi57a]).

We end this part by making one important remark. While it may seem like Adian's papers are significantly more technical than Rabin's, this is not necessarily the case. Adian's papers, particularly [Adi57a], are very self-contained. Essentially the only statements from [Adi57a] not proved therein concern Novikov's group with undecidable word problem, and great effort is made to clarify precisely what statements are used. Indeed, it is no exaggeration to say that, in principle, any reader can pick up [Adi57a] and from there understand the basics of combinatorial (semi)group theory, as even the theory of presentations is given in detail. If abbreviated, the paper would lose a large fraction of its pages -- but it would also lose a large fraction of its charm. 

\

\begin{center}
\textbf{Markov's Theorem(s)}
\end{center}

\noindent As mentioned before, the group-theoretic Adian--Rabin theorem is directly modelled on Markov's corresponding monoid-theoretic theorem. A \textit{Markov property} of finitely presented monoids is defined in exactly the analogous manner to how Markov properties are defined for finitely presented groups. We may then state:

\begin{theorem*}[Markov, 1951]
Let $P$ be a Markov property. Then there does not exist an algorithm which takes as input a finite monoid presentation, and outputs whether or not the monoid defined by this presentation has property $P$.
\end{theorem*}

We make brief comments on Markov's papers [Mar51a, Mar51b] and how they relate to the above theorem. The modern reader (if at all familiar with semigroup theory) may be accustomed to seeing a theorem first be proved for groups, and only later see this theorem be generalised to monoids, sometimes in a perhaps technical or convoluted manner. However, the opposite is true in this case: Markov's theorem predates the Adian--Rabin theorem. What is more, the proof of the Adian--Rabin theorem is directly inspired by, and follows exactly the same key idea as, the proof of Markov's theorem. What is more, the proof for groups is significantly more complicated than that for monoids; the proof of Markov's theorem takes only a few pages. This precedence of the semigroup-theoretic result -- just as for the word problem -- highlights the critical importance of combinatorial semigroup theory for the development of combinatorial group theory, a fact regrettably often omitted from discussions on the latter. 

Now, as for Markov's papers [Mar51a, Mar51b] themselves, we begin by noting that the modern reader should be very aware that for Markov, an \textit{associative system} is what we would today call a \textit{semigroup with identity}, or \textit{monoid} (i.e. a binary associative operation on a set with an identity element), whereas what Markov calls \textit{semigroup} is what we would today call a \textit{cancellative} monoid (i.e. one in which either one of $xy = xz$ or $yx = zx$ implies $y=z$). As kindly pointed out to me by M.\ Kambites and J.\ P.\ Wächter, Markov makes an important use of the identity element in the course of his proof, and e.g.\ it \textit{is}, by contrast, decidable whether a finitely presented semigroup is isomorphic to a free semigroup. The precise way in which the identity element is used is rather subtle; for simplicity, we will in the sequel refer only to \textit{monoids}. 

The first article, [Mar51a], proves undecidability for recognising a few properties; these are listed on p.~84. The key idea is the same as [Mar51b] below (and indeed Adian's papers), i.e. by using ``test'' monoids defined by systems of relations $\mathfrak{D}_{G,H}$ for pairs of words $G, H$. We omit the details (the paper is very short). 

A month after [Mar51a] appeared in print, the article [Mar51b] appeared. From the outset, this article sets out to give a generalisation of the theorem in [Mar51a] (in the same way [Adi57b] generalises results in [Adi57a]). The article proves the undecidability of recognising Markov properties in monoids. We give an outline with some modernisation. The outline of the proof is essentially the same as for groups (or, rather, the outline of the proof for groups is essentially the same as that for monoids!). Namely (up to some technical, but wholly inessential, details regarding the alphabets of the monoids), we take a monoid $S_0$ which does not embed in any monoid with property $\fP$; a monoid $S_1$ with undecidable word problem; and a monoid $S_4$ which has property $\fP$. We then form the monoid free product $S_0 \ast S_1 \ast S_4$. Adding to this monoid four letters $\{ a, b, c, d\}$, for any pair of words $G, H$ in $S_1$ we form a monoid $S_{G,H}$ by adding the relations $cGd = 0$ and $\xi cHd = cHd$, where $\xi$ ranges over the letters $a, b, c, d$. It is then proved: if $G=H$ in $S_1$, then $S_{G,H} \cong S_4$ (Lemma~3). If $G \neq H$ in $S_1$, then $S_0$ embeds in $S_{G,H}$ (Lemma~4). Thus $S_{G,H}$ is a ``test'' monoid for property $\fP$, i.e. $S_{G,H}$ has property $\fP$ if and only if $G=H$ in $S_1$. This completes the proof of the theorem, as $S_1$ has undecidable word problem. We remark that no proofs are provided, but the only non-trivial proof is that of Lemma~4. The beauty of Markov's proof lies in the fact that very little needs to be proved once the (very non-trivial) step of finding what statements must be proved had been discovered.

Of course, the ``test'' monoid $S_{G,H}$ plays a role entirely analogous to that of Adian's groups $\fAqAB$ in [Adi57a], or indeed the groups $F_{qA}$ in [Adi57b]. The marked increase in difficulty for groups is due to the fact that whereas monoid presentations can often be massaged with relative impunity into encoding many properties, even a light touch (i.e. the addition of a few seemingly harmless relations) on a group presentation can sometimes result in the entire group inadvertently collapsing into the trivial group.

The final part (\S17--19) of [Mar51b] is somewhat remarkable, as it contains an undecidability result for monoids which does not have a group-theoretic analogue in any of Adian's papers. Let $\fP$ be an invariant property of submonoids of monoids, such that there is some submonoid of some monoid with property $\fP$, and some submonoid of some monoid without property $\fP$. Then Markov proves (or, rather, states) that given a finitely presented monoid $S$, and a finite set of words $B_i$ over its generators, there is no algorithm which decides whether the submonoid generated by the $B_i$ has property $\fP$. The group-theoretic analogue for this result would only be proved later by Baumslag, Boone \& Neumann (\textit{Math. Scand} \textbf{7} (1959) pp. 191--201).

\clearpage

\begin{center}
\textbf{Final remarks}
\end{center}

\noindent As with any translation, there are some inconsequential technical matters that should be mentioned. First, any references to page numbers (except in bibliographic items) are always with respect to the present translation. Second, I have done my best to eliminate typographical errors. In [Adi57a], there are quite a few scattered throughout, but these were always easy to fix of a similarly harmless nature as e.g. mistyping 37$^\prime$ as 37$^1$. Any typos that remain (or were introduced in the translation process) are, in other words, my own. Third, while I have not attempted to make a facsimile of the original articles, I have attempted to preserve the layout, outline, and general design of the articles in translating them. This includes the typeface. The alignments and numbering of equations, names and numbering of lemmas, and many other similar details have been kept precisely as they were. Thus, the reader wishing to consult both the Russian original and my translation should have no difficulties in switching between the two. I have also included the original bibliographical information at the top of the first page of every article.

I only spotted one mathematical inaccuracy of a non-typographical nature. This is in [Adi57b], immediately following equation (11), at which point it is stated that ``... in a free group two elements do not commute if they are both not equal to the identity element and not equal to one another''. This is not correct --  the centraliser of any non-trivial element of a free group is infinite cyclic. However, the statement it is used to prove is still true, and can be proved with only a minor modification, which is done in [Adi58] (but makes no mention of the error in [Adi57b]). 

Some interesting miscellaneous historical tidbits which fell out as part of this translation are given below:

\begin{enumerate}
\item The undecidability of Dehn's isomorphism problem is usually attributed to Adian and Rabin. However, as Adian himself states (p.~12), this was already proved by Novikov, and is indeed not too difficult to prove directly. 
\item All articles in this translation were communicated by I. M. Vinogradov (except for [Adi57a], which was not communicated by anyone). He certainly gives off the impression of having worked at an impressive pace in reading them as soon as they were submitted. In particular, [Adi58] and [Mar51b] were both accepted after just one day, and -- remarkably -- [Mar51a] was accepted on \textit{the same day} (4 January 1951) on which it was submitted. 
\item In a manner very confusing for bibliographies, the article [Mar51a] bears the same title as his famous 1947 article (\textit{Dokl. Akad. Nauk SSSR}, \textbf{55}:7 (1947), 587--590) proving the undecidability of the word problem in a finitely presented monoid, as well as its sequel (\textit{Dokl. Akad. Nauk SSSR}, \textbf{58}:3 (1947), 353--356). The fact that all these articles were published in the same journal, with many overlapping years, does not simplify this!
\end{enumerate}

Finally, I wish to thank L. D. Beklemishev for encouragement and providing a copy of the Russian original of [Adi55], A. Dmitrieva for linguistic assistance, and M.\ Kambites and J.\ P.\ Wächter for pointing out some subtleties in defining Markov properties for semigroups rather than monoids.


\vspace{0.8cm}

\begin{flushright}
\noindent\textbf{C. F. Nyberg-Brodda}\\
\noindent{The University of Manchester}\\
\today
\end{flushright}



\section*{Supplementary material}

S. I. ADIAN

# ALGORITHMIC UNSOLVABILITY OF PROBLEMS OF RECOGNITION OF CERTAIN PROPERTIES OF GROUPS

*(Communicated by academician I. M. Vinogradov, 26 April 1955)*

We will consider only groups which can be defined by a finite number of generators and defining relations, and will not always state this assumption.

In the article [1], P. S. Novikov constructed a group with algorithmically unsolvable identity problem. Based on this result, he also proved the unsolvability of the isomorphism problem for groups, as formulated in [2].

In the present article we will prove the unsolvability of the isomorphism problem in the case that one of the groups is fixed, taking it as any group given by a finite number of generators and defining relations. More precisely, we have the following theorem:

**Theorem 1.** *Let $F_0$ be an arbitrary group. There does not exist an algorithm which decides for every group $F$ whether or not it is isomorphic to $F_0$.*

We will say that a property $\beta$ of a group $F$ is invariant, if every group $F_1$ isomorphic to $F$ also possesses this property. A property $\alpha$ of a group $F$ will be said to be hereditary, if every subgroup $F'$ of the group $F$ also possesses it. Whether a group possesses or does not possess a given invariant property does not depend on the manner in which the group is specified. Clearly, every hereditary property is invariant. A hereditary property will be said to be non-trivial, if this property is not possessed by some free group with a finite number of generators.

**Theorem 2.** *Let $\alpha_0$ be any invariant property, $\beta_0$ any non-trivial hereditary property, and suppose there exists a group $F_1$ having property $(\alpha_0 \& \beta_0)$. Then there does not exist an algorithm which decides for every group $F$ whether or not it has the property $(\alpha_0 \& \beta_0)$.*

**Theorem 3.** *Suppose there exists a group $F_0$ which possesses a given invariant property $\alpha_0$, which is not possessed by any group with unsolvable identity problem. Then there does not exist an algorithm which decides for every group $F$ whether or not it possesses the property $\alpha$.*

**Theorem 4.** *There does not exist an algorithm which recognises which groups have decidable identity problem, and which do not, i.e. the meta-identity problem is unsolvable.*

**Theorem 5.** *There does not exist an algorithm for recognising the following group properties:*

    *(1) being a simple group;*
    *(2) being a group which decomposes into a free product of $n$ subgroups, for fixed $n > 1$.*
    *(3) having a free subgroup.*

Denote by $\mathfrak{A}_0$ the centrally-symmetric group with unsolvable identity problem constructed in [1]. In [1], it is proved that there does not exist an algorithm which decides which of the words of the form $XpX^+$ (where $p$ is the pivot letter of the group $\mathfrak{A}_0$, and $X$ consists only of letters from the left alphabet) are equal, and which are not. Clearly, the same can be said regarding words of the form $pXpX^+$, where $X$ is an arbitrary word consisting of letters from the left alphabet of the group $\mathfrak{A}_0$. Two letters from the positive alphabet of the group $\mathfrak{A}_0$ will be called homotypic, if they are both from the left alphabet, or if they are both from the right alphabet, of the group $\mathfrak{A}_0$. Otherwise, the two letters will be called heterotypic. Note, that in the group $\mathfrak{A}_0$ more than half of the letters from the positive alphabet are letters which are denoted by $l_{a,i}$ and $l_{a,i}^+$ in [1]. The following lemma is easy to prove:

**Lemma.** *The entire positive alphabet of the group $\mathfrak{A}_0$ can be arranged in a row*

$$a_1, a_2, \ldots, a_{n-1}, a_n, \tag{1}$$

*which satisfies the following properties:*

*1) for $i \leq n-2$ the letters $a_i$ and $a_{i+2}$ are heterotypic;*

*2) for even $i$, $a_i$ is a letter of the form $l_{a,i}$ or $l_{a,i}^+$;*

*3) $a_{n-2}$ is the pivot letter $p$ of the group $\mathfrak{A}_0$;*

*4) $a_n$ is one of the letters $l_{a,i}$, and $a_{n-1}$ is one of the $l_{a,i}^+$.*

Every pair $(A, B)$ of words of the form $pXpX^+$ of the group $\mathfrak{A}_0$ will correspond to some group, which will be denoted by $\mathfrak{A}_{q,A,B}$. The group $\mathfrak{A}_{q,A,B}$ will be defined by giving its positive alphabet and a system of defining relations. The positive alphabet of the group $\mathfrak{A}_{q,A,B}$ will consist of all letters from the row (1) and $n$ new letters

$$q_1, q_2, \ldots, q_{n-1}, q_n. \tag{2}$$

The system of defining relations for the group $\mathfrak{A}_{q,A,B}$ will consist of all defining relations of the group $\mathfrak{A}_0$:

$$C_\nu = D_\nu \tag{3}$$

together with the following relations:

$$a_i q_i = q_i A B^{-1} \qquad (i = 1, 3, 5, \ldots, n-2), \tag{4}$$

$$q_{j-1} q_j = q_j a_{j-1} \qquad (1 < j \leq n), \tag{5}$$

$$q_n a_n = a_n q_n A B^{-1} q_n. \tag{6}$$

The following theorem is also given without proof, due to the cumbersomeness of its proof.

**Theorem 1$^\circ$.** *If $A = B$ in the group $\mathfrak{A}_0$, then $\mathfrak{A}_{q,A,B}$ is the trivial group. If $A \neq B$ in the group $\mathfrak{A}_0$, then $\mathfrak{A}_0$ is a subgroup of the group $\mathfrak{A}_{q,A,B}$.*

From this theorem it follows immediately that there does not exist any algorithm which recognises which groups are trivial, and which are not. A direct consequence of this theorem is also Theorem 4. The rest of the theorems, formulated at the beginning of the article, can also be proved by using Theorem 1$^\circ$.

We give a brief proof of Theorem 1. Suppose the group $F_0$ is given by the generators

$$b_1, b_2, \ldots, b_\lambda \tag{7}$$

and the defining relations

$$E_k = L_k \qquad (k = 1, 2, \ldots, m). \tag{8}$$

Every pair of words $(A, B)$ of the form $pXpX^+$ determines a group $\mathfrak{A}'_{q,A,B}$, defined by the generators (1), (2), and (7), together with the defining relations (3), (4), (6), and (8). The group $\mathfrak{A}'_{q,A,B}$ is a free product of the corresponding group $\mathfrak{A}_{q,A,B}$ by the group $F_0$. From Theorem $1°$ it follows, that if $A = B$ in $\mathfrak{A}_0$, then the group $\mathfrak{A}'_{q,A,B}$ will be isomorphic to the group $F_0$. We now prove that in the case $A \neq B$ in $\mathfrak{A}_0$, the group $\mathfrak{A}'_{q,A,B}$ is not isomorphic to $F_0$. Suppose $A \neq B$ in $\mathfrak{A}_0$, and that $\mathfrak{A}'_{q,A,B}$ is isomorphic to $F_0$. Then neither $\mathfrak{A}_{q,A,B}$ nor $F_0$ is the trivial group, and $\mathfrak{A}'_{q,A,B}$ is their free product. As $F_0$ is isomorphic to the group $\mathfrak{A}'_{q,A,B}$, it will itself be freely decomposable into two factors $F'_1$ and $F_1$, which are themselves isomorphic to $\mathfrak{A}_{q,A,B}$ and $F_0$, respectively. Then $\mathfrak{A}'_{q,A,B}$ will be freely decomposable into three factors $\mathfrak{A}_{q,A,B}$, $F'_1$, and $F_1$. Continuing this reasoning, we find that if, when $A \neq B$ in $\mathfrak{A}_0$, the groups $\mathfrak{A}'_{q,A,B}$ and $F_0$ are isomorphic, then the group $\mathfrak{A}'_{q,A,B}$ is freely decomposable into an arbitrarily large finite number of free factors. But this is impossible since, as proved in [2], the number of factors in a free product decomposition for any group with a finite number of generators is bounded above by the number of generators. It follows that when $A \neq B$ in $\mathfrak{A}_0$, the group $\mathfrak{A}_{q,A,B}$ is not isomorphic to $F_0$. This proves Theorem 1, as there does not exist an algorithm which can decide if $A = B$ in $\mathfrak{A}_0$.

For the proof of Theorem 2, we construct a system of group $\mathfrak{A}'_{q,A,B}$, taking instead of the group $F_0$ the group $F_1$, which possesses the property $(\alpha_0 \,\&\, \beta_0)$. Then we find, that for the property $(\alpha_0 \,\&\, \beta_0)$ to be possessed by $\mathfrak{A}'_{q,A,B}$, it is necessary and sufficient that we have $A = B$ in $\mathfrak{A}_0$.

**Sufficiency.** If $A = B$ in $\mathfrak{A}_0$, then by Theorem $1°$ the group $\mathfrak{A}'_{q,A,B}$ will be isomorphic to the group $F_1$, and, by the assumptions of Theorem 2, will hence possess the invariant property $(\alpha_0 \,\&\, \beta_0)$.

**Necessity.** It suffices to prove, that when $A \neq B$ in $\mathfrak{A}_0$, even the property $\alpha_0$ cannot be possessed by $\mathfrak{A}'_{q,A,B}$. By Theorem $1°$ when $A \neq B$ in $\mathfrak{A}_0$, the group $\mathfrak{A}_0$ will be a subgroup of the group $\mathfrak{A}_{q,A,B}$, and hence also of the group $\mathfrak{A}'_{q,A,B}$. But the group $\mathfrak{A}_0$, as is seen from [1], has a free subgroup on two generators, which in turn has a free subgroup on any finite number of generators, and hence it cannot possess the non-trivial hereditary property $\alpha_0$.

Theorem 3 and Theorem 5 are proved analogously.

**Moscow State V. I. Lenin**
**Pedagogical Institute**

S. I. ADIAN

# UNSOLVABILITY OF SOME ALGORITHMIC PROBLEMS IN THE THEORY OF GROUPS [1]

## Contents

## Introduction

The development of an accurate notion of "algorithm", given recently by a number of authors, has made it possible to prove the unsolvability of some algorithmic problems.

An algorithmic problem is always of a mass type. That is, it corresponds to a whole class of individual questions, each answerable either by "yes" or "no", and consists in finding an algorithm which can solve all the given questions. The unsolvability of an algorithmic problem means only the impossibility of solving all these questions by means of a single general method. In no way, however, does this mean that any of the individual problems is unsolvable.

Theorems about the impossibility of certain algorithmic problems are essential for the development of mathematics, if nothing else to serve as a warning to mathematicians against hopeless searches for non-existent algorithms.

Among the results in this direction, the most outstanding is that by Petr Sergeevich Novikov, who proved the undecidability of the famous identity problem for words in groups given by finitely many generators and defining relations. He also proved the unsolvability of the conjugacy problem [1] and the isomorphism problem for groups defined by finitely many generators and defining relations (P. S. Novikov's proof of the unsolvability of the isomorphism problem has not been published).

In the present article, based on the result by P. S. Novikov, we will prove the unsolvability of some other algorithmic problems in the theory of groups, some of which had been posed earlier, and some of which are new.

---

[1] The main results of the present article were presented at the meeting of the Moscow Mathematical Society held on the 18 October 1955.

This article consists of two chapters. In the first chapter, we introduce the basic concepts and establish a number of properties of the centrally-symmetric groups of P. S. Novikov, which will be necessary for what follows. The first two sections of the second chapter are dedicated to building a system of groups $\mathfrak{A}_{q,A,B}$ and proving the chapter's main lemma, regarding word transformations in these groups. In the third section, based on this main lemma, the following theorem on the unsolvability of algorithmic problems is proved.

**Theorem 1.** *Let $F_0$ be an arbitrary group, given by a finite number of generators and defining relations. There does not exist an algorithm which decides for an arbitrary group $F$, given by a finite number of generators and defining relations, whether or not it is isomorphic to $F_0$.*

A property $\alpha$ of a group $F_0$ is called *invariant* if every group $F_1$ which isomorphic to $F_0$ also has the property. A property $\beta$ of $F_0$ is called *hereditary* if every group $F_1$ which is isomorphic to a subgroup of $F_0$ also has the property. It is clear that every hereditary property is invariant. Whether a given group has or does not have a certain invariant property does not depend on the manner in which this group is specified.

Let $\alpha$ be any invariant property. The *recognition problem for property $\alpha$* is the following problem: find an algorithm, which decides for every group $F$, given by a finite number of generators and defining relations, whether $F$ possesses property $\alpha$ or not.

A hereditary property $\beta_0$ is called a *non-trivial* hereditary property, if this property is not possessed by the free group on two generators.

**Theorem 2.** *Let $\alpha_0$ be an arbitrary invariant property, and $\beta_0$ an arbitrary non-trivial hereditary property. If there exists a group $F_0$ given by a finite number of generators and defining relations and possessing the property $(\alpha_0 \, \& \, \beta_0)^1$, then the recognition problem for the property $(\alpha_0 \, \& \, \beta_0)$ is unsolvable.*

An invariant property $\alpha_0$ is called a *special invariant* property if it is not possessed by any group with unsolvable identity problem.

**Theorem 3.** *Let $\alpha_0$ be an arbitrary special invariant property. If there exists a group $F_0$, given by a finite number of generators and defining relations, which has property $\alpha_0$, then the recognition problem for $\alpha_0$ is unsolvable.*

**Theorem 4.** *The recognition problem is unsolvable for the following properties of groups:*

  *(1) being a simple group,*
  *(2) having a free subgroup of rank $k$.*

Finally, I express my sincere gratitude to my supervisor Petr Sergeevich Novikov for all the help he has provided to me during the undertaking of this work.

---

[1]The symbol & denotes the union "and".

# Chapter I
# The Centrally-symmetric Group of P. S. Novikov

## §1. Basic concepts

It is well-known that every group can be given using generators and defining relations. We will be interested in groups which can be defined with a finite number of generators and relations. The definition of such groups can be given in the following way.

First, we give a system of *letter-symbols*

$$a_1, a_2, \ldots, a_n, \tag{1}$$

which will be called the *system of generators* of the defined group $\mathfrak{A}$. We add the symbols

$$a_1^{-1}, a_2^{-1}, \ldots a_n^{-1} \tag{2}$$

to this system, which are also called letters (moreover, each of the letters of the system (2) is called the *inverse* of the letter from the system (1) to which it corresponds). The collection of symbols (1) and (2) is called the *alphabet* of the group $\mathfrak{A}$. The system (1) is called the *positive* alphabet, and (2) is called the *negative*. An arbitrary letter of the alphabet will be written as $x^\sigma$, where $\sigma$ is either $+1$ or $-1$, and $x^{+1}$ is simply $x$.

A *word* of the group $\mathfrak{A}$ is a finite sequence of letters, taken from the alphabet of this group. We also consider the empty word, which will be denoted by the symbol 1 or $A^0$.

Let $X$ and $Y$ be words of the group $\mathfrak{A}$. The *product* of the word $X$ by the word $Y$ is the word $XY$. Two words will be called *identical*, if they are spelled the same. If the words $X$ and $Y$ are identical, then we denote this by $X \equiv Y$.

If the word $X$ is of the form

$$x_1^{\sigma_1} x_2^{\sigma_2} \cdots x_n^{\sigma_n},$$

then we denote by $X^{-1}$ the word

$$x_n^{-\sigma_n} x_{n-1}^{-\sigma_{n-1}} \cdots x_2^{-\sigma_2} x_1^{-\sigma_1},$$

and call this word the *inverse* for the word $X$.

We now define a relation of equality between the words of the group $\mathfrak{A}$. For this, we give a system of defining relations:

$$\left. \begin{array}{l} a_i a_i^{-1} = 1 \\ a_i^{-1} a_i = 1 \end{array} \right\} \quad (i = 1, 2, 3, \ldots, n), \tag{3}$$

$$A_j = B_j \quad (j = 1, 2, 3, \ldots, m), \tag{4}$$

where $a_i$ and $a_i^{-1}$ are the letters of the alphabet, and $A_j$ and $B_j$ are arbitrary words of the group $\mathfrak{A}$.

The defining relations (3) are called *trivial*. They hold in any group with generators

$$a_1, a_2, a_3, \ldots, a_n.$$

Suppose that

$$R = S$$

is one of the defining relations of the group $\mathfrak{A}$. Then the transition

$$uRv \to uSv$$

and the transition

$$uSv \to uRv,$$

where $u$ and $v$ are arbitrary words of the group $\mathfrak{A}$, are both called *elementary transitions* of the group $\mathfrak{A}$.

Two words $X$ and $Y$ are said to be *equal* in the group $\mathfrak{A}$ if and only if one can be transformed into the other by some sequence of elementary transformations in $\mathfrak{A}$:

$$X \equiv X_0 \to X_1 \to X_2 \to \cdots \to X_k \equiv Y.$$

Equality, defined as such, is reflexive, symmetric, and transitive.

The set of all words of the group $\mathfrak{A}$ is thus divided into *classes* consisting of all words which are equal to each other. Each of these classes is uniquely defined by any one of the words included in it. The product of two classes $\alpha$ and $\beta$ is the class corresponding to the product of some two words taken from the classes $\alpha$ and $\beta$. Obviously, the product of the classes $\alpha$ and $\beta$ does not depend on the choice of these words.

It is easy to prove that the set of all classes of equal words forms a group relative to the product operation defined above. This group is the desired group $\mathfrak{A}$.

As the negative alphabet (2) and the defining relations (3) are uniquely determined by the system of generators (1), it follows that to completely specify the group $\mathfrak{A}$ it suffices to specify the system of generators (1) and the defining relations (4).

For this reason, it is usual to say that the group $\mathfrak{A}$ is defined by the generators

$$a_1, a_2, \ldots, a_n \tag{1}$$

and the defining relations

$$A_j = B_j \quad (j = 1, 2, \ldots, m). \tag{4}$$

If there are no defining relations (4), then $\mathfrak{A}$ is the *free group* on $n$ generators.

Words which are equal in the free group on some alphabet are equal in the free group on any alphabet which contains the letters of this word. Therefore, we can say that two words are equal in the free group regardless of which group they belong to. In all other cases, every time we speak of equality of words, we will specify the group in which this equality takes place.

The elementary transformation of the form

$$XY \to Xa^{-1}aY$$

will be called a *right-sided insertion of the letter $a$*, and a transformation of the form

$$XY \to Xaa^{-1}Y$$

will be called a *left-sided insertion of the letter $a$*.

The reversal of these transformations will be called a *left-sided deletion* resp. a *right-sided deletion* of the letter $a$.

Given an elementary transformation $XRY \to XSY$, where either $R = S$ or $S = R$ is a defining relation of the group under consideration, then the letters of the words $XRY$, which constitute the specified instance of $R$ in this word, will be said to be *affected* by the transformation $XRY \to XSY$, and the letters which constitute the specified instances of $X$ and $Y$ of the word $XRY$ will be called *unaffected* by this transformation. A segment of the word $XRY$ will be said to be unaffected by the transformation $XRY \to XSY$, if it does not contain a single letter affected by this transformation.

Consider an arbitrary group $\mathfrak{A}$ with generating alphabet $A$. Suppose $P$ is an alphabet which forms a part of the alphabet $A$ (or, as it is called, $P$ is a *subalphabet* of $A$), and let $X$ be any word made up of letters from the alphabet $A$. The word obtained from $X$ by deleting all letters belonging to the alphabet $P$ is called the *projection of the word $X$ onto the alphabet $P$*, and is denoted by $(X)_P$.

The *projection of a sequence of elementary transformations*

$$X_1 \to X_2 \to \cdots \to X_m$$

will be the sequence

$$(X_1)_P \to (X_2)_P \to \cdots \to (X_m)_P,$$

which, generally speaking, will not be transformations of the group $\mathfrak{A}$.

The subset of letters

$$p_1, p_2, \ldots, p_k, p_1^{-1}, p_2^{-1}, \ldots, p_k^{-1} \tag{5}$$

of the group $\mathfrak{A}$ will be called a *system of filtering letters*, if all defining relations of the group $\mathfrak{A}$ containing the letter $p_j^{\sigma}$ (with the exception of the trivial defining relations (3)) are of the form

$$E_i^j p_j F_i^j = H_i^j p_j G_i^j, \tag{6}$$

where the words $E_i^j, F_i^j, H_i^j$ and $G_i^j$ do not contain letters from the subset (5).

Suppose the other relations are of the form

$$A_i = B_i. \tag{7}$$

By definition, the words $A_i$ and $B_i$ do not contain the letter $p_j^{\sigma}$, and the letter $p_j^{-1}$ appears only in the trivial defining relations.

The system of filtering letters (5) forms a subalphabet of the alphabet of the group $\mathfrak{A}$. Denote this system by $P$.

Consider an arbitrary sequence of elementary transformations

$$X_1 \to X_2 \to \cdots \to X_j \to X_{j+1} \to \cdots \to X_m \tag{8}$$

of the group $\mathfrak{A}$. The letters of the system $P$ which are involved in the transformations of (8) will be said to be *identical in individuality*.

Any transformation $X_j \to X_{j+1}$ has the form

$$X_j' R X_j'' \to X_j' S X_j''. \tag{9}$$

If in the transformation (9) there are no insertions or deletions of any of the filtering letters of the system $P$, then the words $(X_j' R X_j'')_P$ and $(X_j' S X_j'')_P$ are identical. In this case the letters of the alphabet $P$ contained in the words $(X_j' R X_j'')_P$ and $(X_j' S X_j'')_P$ and having the same number, are called *identical in individuality in the transformation (9)*. If (9) is an insertion or a deletion of one of the letters from the system $P$, then we will call *identical in individuality* the letters from $P$ of the transformation (9) which have the same number in the identical segments $X_j'$ (resp. $X_j''$) of the words $X_j' R X_j''$ and $X_j' S X_j''$.

We extend, by transitivity, the definition of letters identical in individuality of the system $P$ across the entire sequence (8). Letters which are identical in individuality will also be said to have the same individuality. Two letters of the alphabet $P$, included in some word of the sequence (8) of transformations, will be called *different in individuality* if the identity of their individualities is not deducible from the above definition.

It is clear that letters which are identical in individuality are spelled the same.

We note already now that in the future, where this does not cause confusion, we will instead of "the individuality of the letter $a$" simply say "the letter $a$".

If we speak of filtering letters in the sequence of transformations (8) without distinguishing between their individuality, then it should be noted that each of them is either contained in the first word of the sequence (8), or else is introduced by an insertion and

is then present in all subsequence words, unless it disappears, having been reduced with its reverse letter. This leads to the following lemma (see [2], p. 494).

**Lemma 1.** *If the words $X$ and $Y$ are equal in the group $\mathfrak{A}$ with a system of filtering letters $P$, then their projections onto $P$ are equal in the free group.*

From the defining relation (6), using the trivial defining relations (3), we can deduce the equalities

$$H_i^{j-1} E_i^j p_j = p_j G_i^j F_i^{j-1} \tag{10}$$

and

$$E_i^{j-1} H_i^j p_j = p_j F_i^j G_i^{j-1}. \tag{11}$$

For example,

$$H_i^{j-1} E_i^j p_j \to \cdots \to H_i^{j-1} E_i^j p_j F_i^j F_i^{j-1} \to H_i^{j-1} H_i^j p_j G_i^j F_i^{j-1} \to \cdots \to p_j G_i^j F_i^{j-1}.$$

In the same way, we can obtain from either of the equalities (10) or (11) the corresponding equality (6). Hence, if we replace each defining relation (6) by the pair of corresponding relations (10) and (11), then the group $\mathfrak{A}$ will not be changed, i.e. in $\mathfrak{A}$ we find no new equalities between words, and all old equalities can still be obtained.

Suppose that in the group under consideration with system of filtering letters (5) all defining relations (6) are replaced by the corresponding (10) and (11), which in general will be written in the form

$$U^{(s)} p_j = p_j W^{(s)}. \tag{12}$$

Obviously, for every word $U^{(s)}$ there is also a word $U^{(s')}$ such that

$$U^{(s)} \equiv (U^{(s')})^{-1}.$$

The same holds true also for the words $W^{(s)}$.

Let us compose all possible products of the form

$$U^{(s_1)} U^{(s_2)} \dots U^{(s_k)}$$

and the products of the form

$$W^{(s_1)} W^{(s_2)} \dots W^{(s_k)}.$$

We will place the first set in a one-to-one correspondence with the second via:

$$U^{(s_1)} U^{(s_2)} \dots U^{(s_k)} \underline{\hspace{2cm}} W^{(s_1)} W^{(s_2)} \dots W^{(s_k)}.$$

This correspondence will be called the *correspondence $L_P$*. In what follows, the words of the form $U^{(s_1)} U^{(s_2)} \dots U^{(s_k)}$ we will denote with a $U$ with varying indices, and the corresponding word $W^{(s_1)} W^{(s_2)} \dots W^{(s_k)}$ by $W$ with the same index.

From the relation (12) we find the equality

$$Up = pW$$

for every $U$ and $W$ corresponding to each other under the correspondence $L_P$. This equality may be written in the form

$$U = pWp^{-1}.$$

It is easy to see that the set of words $U$, as well as the set of words $W$, forms a subgroup of the group $\mathfrak{A}$. From the previous relation it is clear that the correspondence $L_P$ is an isomorphism of these subgroups. Recall that the words $U$ and $W$ do not contain any occurrences of the letter $p$.

Consider a group, the alphabet of which is derived from the alphabet of the group $\mathfrak{A}$ by removing all filtering letters $p_j^\sigma$, and with defining relations (7):

$$A_i = B_i$$

(the words $A_i$ and $B_i$ do not, by definition, contain any filtering letters).

We call this group the *base of the group $\mathfrak{A}$ with respect to the system of filtering letters $P$* and denote it by $\overline{\mathfrak{A}}^{(P)}$ or simply $\overline{\mathfrak{A}}$, if there is no risk of confusion arising from this usage. The set of words $U$ and the set of words $W$ form two subgroups also of the group $\overline{\mathfrak{A}}$, but the correspondence $L_P$, generally speaking, will no longer be an isomorphism in the group $\overline{\mathfrak{A}}$.

If for every filtering letter $p_j$ the correspondence $L_P$ is an isomorphism in the group $\overline{\mathfrak{A}}$, then the system of filtering letters $P$ will be called *regular*.

We will always denote the word $W^{-1}$ by $V$ with the same indices (in the case that the word $W$ has them). If

$$W \equiv W^{(s_1)} \cdots W^{(s_k)},$$

then

$$W^{-1} \equiv V \equiv V^{(s_k)} \cdots V^{(s_1)}.$$

If $U$ and $V^{-1}$ correspond to the same word under $L_P$, then $U$ and $V$ are of the form

$$U \equiv U^{(s_1)} \cdots U^{(s_k)}, \; V \equiv V^{(s_k)} \cdots V^{(s_1)}.$$

In order to prove that the correspondence $L_P$ is an isomorphism, it suffices to prove that the words $U$ and $V$, where $U$ and $V^{-1}$ correspond to one another under $L_P$, are simultaneously equal or not equal to the identity element in the group $\overline{\mathfrak{A}}$.

In the sequel, we will as a rule consider the case when the group contains one filtering letter. Notice that the correspondence $L_P$ can easily be formed directly from the equality (6). If all these equalities for the filtering letter $p$ have the form

$$E_i p F_i = H_i p G_i \quad (i = 1, 2, \ldots, \mu), \tag{6'}$$

then the correspondence $L_P$ can be written in the form

$$U \equiv (H_{i_1}^{-1} E_{i_1})^{\sigma_1} (H_{i_2}^{-1} E_{i_2})^{\sigma_2} \cdots (H_{i_k}^{-1} E_{i_k})^{\sigma_k} \rightarrow$$
$$\rightarrow W \equiv (G_{i_1}^{-1} F_{i_1})^{\sigma_1} (G_{i_2}^{-1} F_{i_2})^{\sigma_2} \cdots (G_{i_k}^{-1} F_{i_k})^{\sigma_k},$$

where $\sigma_1, \sigma_2, \ldots, \sigma_k$ all equal $\pm 1$.

In the sequel, we will frequently make use of the following lemmas, proved in [2].

**Lemma 2.** *Suppose that*

$$X_0 p Y_0 \rightarrow X_1 p Y_1 \rightarrow \cdots \rightarrow X_m p Y_m \tag{13}$$

*is a sequence of elementary transformations of the group $\mathfrak{A}$, all words of which contain the filtering letter $p$ in its individuality (this letter is explicitly indicated in (13)). Then the word $X_0 p Y_0$ can be transformed into $X_m p Y_m$ by a sequence of transformations of the form*

$$X_0 p Y_0 \rightarrow X_0 U p V Y_0 \rightarrow \cdots \rightarrow X_m p Y_m,$$

*with the following properties:*

(1) *$U$ and $V^{-1}$ are words corresponding to each other under $L_P$.*

(2) *The projection of the obtained sequence of transformations onto the alphabet $P$ coincides up to tautological repetitions with the projection onto this alphabet of the sequence (13).*

*(3) In the sequence*

$$X_0 U p V Y_0 \rightarrow \cdots X_m p Y_m \tag{14}$$

*no transformations affect the distinguished letter p.*

We remark that in [2] this lemma was proved in a slightly more general form.

The sequence of transformations (14) will be called the *canonical sequence for (13), associated with the individuality of the letter p*.

**Lemma 3.** *If $X = Y$ in the group $\mathfrak{A}$ with a regular system of filtering letters from $P$ and the projection of $Y$ onto $P$ is irreducible, then $X$ can be transformed into $Y$ by a sequence of elementary transformations without any insertions of filtering letters.*

This statement appears in [2] as Corollary 1 of Lemma 3, proved therein.

**Corollary.** *Let $P$ be a regular system of filtering letters of the group $\mathfrak{A}$. If the word $E$ does not contain filtering letters and $E = 1$ in $\mathfrak{A}$, then $E = 1$ also in the group $\overline{\mathfrak{A}}^{(P)}$.*

Indeed, in this case we have a sequence of elementary transformations in $\mathfrak{A}$

$$E \rightarrow \cdots \rightarrow 1 \tag{15}$$

without insertions of filtering letters. As filtering letters are absent from the word $E$, they do not appear in the transformation (15). This means that among these transformations there are no instances of those defining relations which contain filtering letters on both sides, i.e. (15) is a sequence of transformations also in the group $\overline{\mathfrak{A}}^{(P)}$.

Suppose that a group $\mathfrak{A}$ is given by an alphabet $A$ an some system of defining relations. We will say that the defining relations of the group $\mathfrak{A}$ *divides the alphabet $A$* into two parts $A_1$ and $A_2$ if:

(1)  $A_1$ and $A_2$ are alphabets with empty intersection and whose union is $A$.

(2)  Every defining relation is either written over the generating alphabet $A_1$, or else over the alphabet $A_2$. This means that the system of defining relations of the group $\mathfrak{A}$ divides into two parts **I** and **II**, such that the relations from **I** are written over the alphabet $A_1$, and the relations from **II** over the alphabet $A_2$.

It is known that in the above case, $\mathfrak{A}$ is the free product of the groups $\mathfrak{A}_1$ and $\mathfrak{A}_2$, these groups being defined on the alphabet $A_1$ resp. $A_2$ with defining relations **I** resp. **II**.

The *free product* of two groups $\mathfrak{A}_1$ and $\mathfrak{A}_2$ will be denoted by $\mathfrak{A}_1 * \mathfrak{A}_2$. If the groups $\mathfrak{A}_1$ and $\mathfrak{A}_2$ are given via their generators and defining relations, then the group $\mathfrak{A}_1 * \mathfrak{A}_2$ is defined as follows. Rename, if needed, the generators of the group $\mathfrak{A}_2$ such that the alphabets of $\mathfrak{A}_1$ and $\mathfrak{A}_2$ do not intersect. Then take as defining alphabet of the group $\mathfrak{A}_1 * \mathfrak{A}_2$ the union of the alphabets of the groups $\mathfrak{A}_1$ and $\mathfrak{A}_2$, and take as system of defining relations of the group $\mathfrak{A}_1 * \mathfrak{A}_2$ the union of the systems of defining relations for the groups $\mathfrak{A}_1$ and $\mathfrak{A}_2$ (where the relations of the group $\mathfrak{A}_2$ are written over the renamed alphabet).

**Lemma 4.** *Let $\mathfrak{A} \equiv \mathfrak{A}_1 * \mathfrak{A}_2$. If $P_1$ is a regular system of filtering letters of the group $\mathfrak{A}_1$, and $P_2$ is a regular system of filtering letters of $\mathfrak{A}_2$, then their union $P_0$ is a regular system of filtering letters of the group $\mathfrak{A}$.*

We will first show that $P_0$ is a system of filtering letters of the group $\mathfrak{A}$. Suppose that some defining relation of $\mathfrak{A}$ contains some letter $p \in P_0$. This defining relation, by assumption, is either a defining relation of the group $\mathfrak{A}_1$, or else of the group $\mathfrak{A}_2$. Suppose it is a relation of the group $\mathfrak{A}_1$. The second case can be treated analogously.

Then $p \in P_1$, i.e. $p$ is a filtering letter of the group $\mathfrak{A}_1$. Therefore the defining relation under consideration has the form

$$E_i p F_i = H_i p G_i, \tag{6'}$$

where the words $E_i, F_i, H_i$ and $G_i$ do not contain the letter $p^\sigma \in P_1$, and as the relation (6') contains no letter from the alphabet of $\mathfrak{A}_2$, these words also do not contain any letter from $P_2$. We conclude that $P_0$ is a system of filtering letters of the group $\mathfrak{A}$.

By similar considerations, we conclude that for any fixed filtering letter $p \in P_0$, the correspondence $L_P$ in the group $\mathfrak{A}$ coincides with the correspondence $L_p$, defined for the letter $p$ in the group $\mathfrak{A}_1$ (or $\mathfrak{A}_2$).

We have now proved, that $P_0$ is a regular system of filtering letters, i.e. that every correspondence $L_P$ will be an isomorphism in $\overline{\mathfrak{A}}^{(P_0)}$. It is clear that $\overline{\mathfrak{A}}^{(P_0)}$ is the free product of $\overline{\mathfrak{A}}_1^{(P_1)}$ and $\overline{\mathfrak{A}}_2^{(P_2)}$, where $\overline{\mathfrak{A}}_1^{(P_1)}$ resp. $\overline{\mathfrak{A}}_2^{(P_2)}$ is the base of the group $\mathfrak{A}_1$ resp. $\mathfrak{A}_2$ with respect to the system of filtering letters $P_1$ resp. $P_2$. By the definition of the free product of groups, the equality of the word $U$ or $V$ from the correspondence $L_P$ to the identity element of $\overline{\mathfrak{A}}_1$ (resp. $\overline{\mathfrak{A}}_2$) is equivalent to its equality to the identity element of the group $\overline{\mathfrak{A}}_1 * \overline{\mathfrak{A}}_2$. Hence, if the correspondence $L_P$ is an isomorphism in $\overline{\mathfrak{A}}_1$ (resp. $\overline{\mathfrak{A}}_2$), then it is an isomorphism also in $\overline{\mathfrak{A}}_1 * \overline{\mathfrak{A}}_2$, which is what was to be shown.

**Lemma 5.** *Suppose the defining relations of the group $\mathfrak{A}$ divide the alphabet into two parts $A_1$ and $A_2$, and let $E$ be an arbitrary word of the group $\mathfrak{A}$. If $E = 1$ in $\mathfrak{A}$, then*

$$(E)_{A_1} = 1 \text{ in } \mathfrak{A}, \quad \text{and} \quad (E)_{A_2} = 1 \text{ in } \mathfrak{A}.$$

*Proof.* The equality $E = 1$ in $\mathfrak{A}$ implies that there exists a sequence of elementary transformations of the group $\mathfrak{A}$

$$E \equiv E_0 \to E_1 \to \cdots \to E_m \equiv 1. \tag{16}$$

Consider the projection of this sequence onto the alphabet $A_1$

$$(E)_{A_1} \equiv (E_0)_{A_1} \to (E_1)_{A_1} \to \cdots \to (E_m)_{A_1} \equiv 1. \tag{17}$$

We will show that (17) is also a sequence of elementary transformations in the group $\mathfrak{A}$. Every transformation in the sequence (16) is of the following form:

$$E_i' R E_i'' \to E_i' S E_i'',$$

where $R = S$ is one of the defining relations of the group $\mathfrak{A}$. The corresponding transition in the sequence (17) is hence written in the form

$$(E_i')_{A_1} (R)_{A_1} (E_i'')_{A_1} \to (E_i')_{A_1} (S)_{A_1} (E_i'')_{A_1}. \tag{18}$$

As every defining relation $R = S$ of the group $\mathfrak{A}$ is by assumption either written over the alphabet $A_1$ or over the alphabet $A_2$, the transformation (18) is either a transformation of the form $XRY \to XSY$, or else of the form $X \to X$. Thus, each transformation in the sequence (17) is in fact an elementary transformation in the group $\mathfrak{A}$. This proves that $(E)_{A_1} = 1$ in the group $\mathfrak{A}$.

Considering the projection of the sequence (16) onto the alphabet $A_2$, we can analogously prove that $(E)_{A_2} = 1$ in $\mathfrak{A}$. $\qquad\square$

## §2. The centrally-symmetric group of P. S. Novikov

P. S. Novikov constructed in [1] a centrally-symmetric group[1] with finitely many generators and defining relations and in which the identity problem is unsolvable. This group will be denoted by $\mathfrak{A}_P$.

The positive alphabet of the group $\mathfrak{A}_P$ contains a letter $p$, which is called the *pivot*. The remaining letters of the positive alphabet fall into two classes $S$ and $S^+$ containing the same number of elements. These two classes are brought into a one-to-one correspondence of the form

$$x \relbar\joinrel\relbar x^+,$$

where $x$ is an arbitrary letter from the class $S$, and $x^+$ is its corresponding letter in the class $S^+$.

Join to the letters of the class $S$ all of their respective inverse letters. Denote the resulting system by $\overline{S}$. Then $\overline{S}$ is an alphabet. Analogously, we complete $S^+$ to form the alphabet $\overline{S^+}$. In [1], the alphabet $\overline{S}$ is called the left alphabet of the group $\mathfrak{A}_P$, and $\overline{S^+}$ is called the right alphabet. For words consisting of letters from the alphabets $\overline{S}$ and $\overline{S^+}$, we define the operation denoted by $+$ as follows:

$$\left.\begin{aligned}
(x^{-1})^+ &\equiv (x^+)^{-1}, \\
((x^\sigma)^+)^+ &\equiv x^\sigma,
\end{aligned}\right\} \tag{1}$$

where $x^\sigma$ denotes an arbitrary letter from the alphabet $\overline{S}$. If $X = x_1 x_2 \cdots x_k$, where the $x_i$ are letters from $\overline{S}$ or $\overline{S^+}$, then

$$X^+ = x_k^+ x_{k-1}^+ \cdots x_1^+. \tag{2}$$

From (1) and (2) follow

1) $(X^+)^+ \equiv X$,
2) $(X^{-1})^+ \equiv (X^+)^{-1}$,

where $X$ is any word consisting of letters from the alphabets $\overline{S}$ and $\overline{S^+}$.

From the character of the defining relations of the group $\mathfrak{A}_P$ follows another property of this operation:

3) If $X$ and $Y$ are two words which do not contain the letter $p^\sigma$, and $X = Y$ in $\mathfrak{A}_P$, then $X^+$ and $Y^+$ are also equal in $\mathfrak{A}_P$.

The words $X$ and $X^+$, as well as the equalities $X = Y$ and $X^+ = Y^+$, are called mutually symmetric. Property 3) of the operation $+$ is called the *symmetry* of $\mathfrak{A}_P$.

The defining relations of $\mathfrak{A}_P$ are written using the letters of the positive alphabet. Those that contain the pivot letter have the form

$$E_i p E_i^+ = H_i p H_i^+ \quad (i = 1, 2, \ldots, \nu), \tag{3}$$

where the words $E_i$ and $H_i$ consist of letters from the class $S$. The remaining defining relations divide into two classes **I** and **II** and have the form

$$\textbf{I} \quad A_j = B_j, \qquad \textbf{II} \quad A_j^+ = B_j^+,$$

where the words $A_j$ and $B_j$ consist of letters from the class $S$. Every defining relation of type **I** corresponds to its symmetric relation of type **II**, and vice versa.

The pivot letter $p$ together with the letter $p^{-1}$ forms in the group $\mathfrak{A}_P$ a regular system of filtering letters. We will for brevity say that the letter $p$ is a regular filtering letter. All the lemmas about filtering letters given in §1 remain true for this letter.

---

[1] The term *centrally-symmetric group* was introduced by P. S. Novikov in [2].

The correspondence $L_p$ for the letter $p$ has the form

$$U \,\text{———}\, (U^+)^{-1},$$

where $U$ is a word of the form

$$(H_{i_1}^{-1} E_{i_1})^{\sigma_1} (H_{i_2}^{-1} E_{i_2})^{\sigma_2} \cdots (H_{i_k}^{-1} E_{i_k})^{\sigma_k}.$$

Note that for the word $U$ the corresponding word $V$, which appears in a canonical sequence of transformations (see p. 19), will be the word

$$((U^+)^{-1})^{-1} = U^+.$$

The defining relations of the base $\overline{\mathfrak{A}}_P$ of the group $\mathfrak{A}$ consists of the relations **I** and **II**, which divide the alphabet of $\overline{\mathfrak{A}}_P$ into $\overline{S}$ and $\overline{S^+}$.

**Lemma 1.** *Suppose $E$ is a word of the group $\mathfrak{A}_P$ not containing the letter $p^\sigma$. If $E = 1$ in $\mathfrak{A}_P$, then*

$$(E)_{\overline{S}} = 1 \text{ in } \mathfrak{A}_P \quad \text{and} \quad (E)_{\overline{S^+}} = 1 \text{ in } \mathfrak{A}_P.$$

*Proof.* As $E$ does not contain the letter $p^\sigma$, then it follows from Lemma 3 of §1 that $E = 1$ in $\overline{\mathfrak{A}}_P$. But the alphabet of the group $\overline{\mathfrak{A}}_P$ is divided by the defining relations into $\overline{S}$ and $\overline{S^+}$. Hence, by Lemma 5 of §1 we have

$$(E)_{\overline{S}} = 1 \text{ in } \overline{\mathfrak{A}}_P \quad \text{and} \quad (E)_{\overline{S^+}} = 1 \text{ in } \overline{\mathfrak{A}}_P.$$

In particular, we also have these equalities in the group $\mathfrak{A}_P$.                    □

In [1] (pp. 69 and 140) the following theorems are proved.

**Theorem 1.** *No word consisting only of letters of the positive alphabet of the group $\mathfrak{A}_P$ is equal to 1 in this group.*

**Theorem 2.** *There does not exist an algorithm which decides for an arbitrary pair of words of the form $XpX^+$ (where $X$ is any word consisting of letters from the class $S$) whether or not they are equal in the group $\mathfrak{A}_P$.*

It is clear that this theorem is also true for words of the form $pXpX^+$, i.e. we have

**Theorem 2′.** *There does not exist an algorithm which decides for an arbitrary pair of words of the form $pXpX^+$ (where $X$ is any word consisting of letters from the class $S$), whether or not they are equal in the group $\mathfrak{A}_P$.*

The number $n$ of letters in the positive alphabet of the group $\mathfrak{A}_P$ is odd. A large majority of these are letters which are denoted in [1] by $l_{a,i}$ and $l_{a,i}^+$. The letters of type $l_{a,i}$ are from the class $\overline{S}$, and the letters of the type $l_{a,i}^+$ are from the class $\overline{S^+}$. We are interested in the property of these letters which says that for every letter $l_{a,i}$, the group $\mathfrak{A}_P$ has a defining relation of the form

$$\mu_i a = a \mu_i l_{a,i}, \tag{4}$$

and for every letter $l_{a,i}^+$ a relation of the form

$$a^+ \mu_i^+ = l_{a,i}^+ \mu_i^+ a^+, \tag{4′}$$

where $a, \mu_i, a^+$ and $\mu_i^+$ are letters from the positive alphabet of $\mathfrak{A}_P$, but not of type $l_{a,i}$ or $l_{a,i}^+$.

**Remark 1.** If to the defining relations of the group $\mathfrak{A}_P$ new relations are added from which one could deduce that all letters of the group $\mathfrak{A}_P$ not of type $l_{a,i}$ or $l_{a,i}^+$ are equal to the identity of $\mathfrak{A}_P$, then in the group thus obtained all letters of the type $l_{a,i}$ and $l_{a,i}^+$ will also be equal to the identity.

Indeed, substituting the identity for $a$ and $\mu_i$ in (4) and for $a^+$ and $\mu_i^+$ in (4′), it follows from (4) that $1 = l_{a,i}$ and from (4′) that $1 = l_{a,i}^+$.

**Lemma 2.** *Let $L$ be a word containing only letters of the form $(l_{a,i})^\sigma$ and $(l_{a,i}^+)^\sigma$, where $\sigma = \pm 1$. If $L = 1$ in $\mathfrak{A}_P$, then $L = 1$ in the free group.*

*Proof.* As $p$ is a regular filtering letter in $\mathfrak{A}_P$, it follows that $L = 1$ in the group $\overline{\mathfrak{A}}_P$ (see the Corollary of Lemma 3, §1 on p. 19). The group $\overline{\mathfrak{A}}_P$ is the free product of two groups, which we will denote by $\mathfrak{A}_d$ and $\mathfrak{A}_d^+$. The group $\mathfrak{A}_d$ has as generators the left alphabet of the group $\mathfrak{A}_P$ and defining relations **I**, and $\mathfrak{A}_d^+$ has as generators the right alphabet and defining relations **II**. In [1], Chapter III, it is proved that in the group $\mathfrak{A}_d$, which is there denoted by $\mathfrak{A}_{d,\mu,l,\rho}$, the system of letters there denoted by $d_i^\sigma$ is a regular system of filtering letters. In the group $\mathfrak{A}_d^+$ the system of letters $(d_i^+)^\sigma$ is also a regular system of filtering letters. The union of all letters $d_i^\sigma$ and $(d_i^+)^\sigma$ will be denoted by $D$. By Lemma 4, §1 it follows that $D$ is a regular system of filtering letters for the group $\overline{\mathfrak{A}}_P$. As, by assumption, $L$ does not contain any of the letters $d_i^\sigma$ or $(d_i^+)^\sigma$, it follows that $L = 1$ in the group $\overline{\mathfrak{A}}_P^{(D)}$, i.e. the group which is the base of the group $\overline{\mathfrak{A}}_P$ with respect to the system $D$. For brevity, we will denote this group by $\mathfrak{B}$. The system of defining relations of this group is obtained from the system of defining relations of the group $\overline{\mathfrak{A}}_P$ by removing any defining relation which contains some letter $d_i^\sigma$ or $(d_i^+)^\sigma$.

The group $\mathfrak{B}$ is also a free product of two groups, one of which has as generators the left alphabet of this group (and which we will denote by $\mathfrak{B}_\mu$), and the other which has the right alphabet (and which we will denote by $\mathfrak{B}_\mu^+$). In [1] §3 Chapter III, it is proved that in the group $\mathfrak{B}_\mu$, which is there denoted $\mathfrak{A}_{\mu,l,\rho}$, the system consisting of all letters $\mu$ forms a regular system of filtering letters. Additionally, in $\mathfrak{B}_\mu^+$ the system of letters $\mu^+$ forms a regular system of filtering letters. The union of all $\mu$ and $\mu^+$ forms, by virtue of the same Lemma 4, §1, a regular system of filtering letters in the group $\mathfrak{B} \equiv \mathfrak{B}_\mu * \mathfrak{B}_\mu^+$. As the word $L$ does not contain any letters $\mu^\sigma$ or $(\mu^+)^\sigma$, it follows that $L = 1$ in the base $\overline{\mathfrak{B}}^{(\mu)}$ of the group $\mathfrak{B}$ with respect to the system of filtering letters $\mu^\sigma$ and $(\mu^+)^\sigma$.

The group $\overline{\mathfrak{B}}^{(\mu)}$ is a free product of two groups $\mathfrak{B}_l$ and $\mathfrak{B}_l^+$, in which the systems of letters $l_{a,i}^\sigma$ resp. $(l_{a,i}^+)^\sigma$ form regular systems of filtering letters (as is shown in [1]). By Lemma 4, §1, the system of all letters $l_{a,i}^\sigma$ and $(l_{a,i}^+)^\sigma$ forms a regular system of filtering letters in the group $\overline{\mathfrak{B}}^{(\mu)} \equiv \mathfrak{B}_l * \mathfrak{B}_l^+$. But the word $L$ contains only these letters, and $L = 1$ in $\overline{\mathfrak{B}}^{(\mu)}$. By Lemma 1, §1 it follows that $L = 1$ in the free group, which is what was to be shown. □

If we consider the subgroup of $\mathfrak{A}_P$ generated by two distinct letters of type $l_{a,i}$, then from Lemma 2 it follows that:

**Corollary.** *The group $\mathfrak{A}_P$ contains a free subgroup on two generators.*

**Definition.** Two distinct letters $x_1$ and $x_2$ of the positive alphabet of the group $\mathfrak{A}_P$ will be said to be *homotypic* if both belong to the class $S$, or if both belong to the class $S^+$. Otherwise, they will be said to be *heterotypic*.

**Lemma 3.** *The positive alphabet of the group $\mathfrak{A}_P$ can be arranged into a sequence*

$$a_1, a_2, \ldots, a_{n-2}, a_{n-1}, a_n, \tag{5}$$

*which satisfies the following properties:*

    *(1) $a_{n-2}$ is the pivot letter $p$ of the group $\mathfrak{A}_P$.*

    *(2) $a_n$ is a letter of type $l_{a,i}$, and $a_{n-1}$ is of type $l_{a,i}^+$.*

    *(3) For every $j \leq n-2$ the letters $a_j$ and $a_{j+2}$ are heterotypic.*

    *(4) For even $j$ the letter $a_j$ is of type $l_{a,i}$ or of type $l_{a,i}^+$.*

*Proof.* The sequence (5) is built as follows. We set $a_n$ to be any arbitrary letter of type $l_{a,i}$, and set $a_{n-1}$ to be its corresponding letter of type $l_{a,i}^+$. We set $a_{n-2}$ to be the pivot letter $p$. There remain $n-3$ letters, of which equally many are from $S^+$. At least half of them will be letters of type $l_{a,i}$ and $l_{a,i}^+$. The number $n-3$ is even.

  I. Suppose that $n-3$ is divisible by 4, i.e. that $n-3 = 4k_0$. Then the remaining elements of the sequence $(S)$ will be built in ascending order, for $j$ from 1 to $n-3$, according to the following rule:

    a) For $j = 4k+1$ $(k = 0, 1, \ldots, k_0 - 1)$, we will take as $a_j$ any arbitrary letter of the class $S^+$ not of type $l_{a,i}^+$, unless all such letters have been selected. If they have, then we will take $a_j$ to be any of the remaining letters of type $l_{a,i}^+$.

    b) For $j = 4k+2$ $(k = 0, 1, \ldots, k_0 - 1)$, we will take as $a_j$ any of the remaining letters of type $l_{a,i}^+$.

    c) For $j = 4k+3$ $(k = 0, 1, \ldots, k_0 - 1)$, we take $a_j$ from $S$, not of type $l_{a,i}$; if no such letters remain, then we choose any of the remaining letters $l_{a,i}$.

    d) For $j = 4k+4$ $(k = 0, 1, \ldots, k_0 - 1)$, we will take as $a_j$ any of the remaining of the letters of type $l_{a,i}$.

  As any given letter will not be selected twice, it follows that for $k = 0, 1, \ldots, k_0 - 1$ we will choose all letters from the positive alphabet of the group $\mathfrak{A}_P$. As most of the remaining letters are of type $l_{a,i}$ or $l_{a,i}^+$, there is no difficulty in doing this.

  II. If $n-3$ is not divisible by 4, then $n-3 = 4k_1 + 2$. In this case, we easily construct the sequence (5) by choosing $a_j$ for $j = 1, 2, \ldots, n-3$ by to the following rule:

    a) For $j = 4k+1$ $(k = 0, 1, \ldots, k_1)$, we choose $a_j$ as in Ic).

    b) For $j = 4k+2$ $(k = 0, 1, \ldots, k_1)$, we choose $a_j$ as in Id).

    c) For $j = 4k+3$ $(k = 0, 1, \ldots, k_1 - 1)$ – as in case Ia).

    d) For $j = 4k+4$ $(k = 0, 1, \ldots, k_1 - 1)$ – as in case Ib).

  As in case I, the letter $a_{n-3}$ is also in case II chosen by the rule Id), and clearly belongs to $S$. Hence, it is heterotypic to $a_{n-1} \in S^+$. Now $a_{n-2}$ is heterotypic with every other letter of the alphabet of $\mathfrak{A}_P$, including $a_n$ and $a_{n-4}$. Property 3) of the remaining letters of the sequence (5) and property 4) easily follow from conditions I a), b), c), and d), or from conditions II) a), b), c), and d). $\qquad\square$

We fix some concrete sequence (5) for the letters of the positive alphabet of the group $\mathfrak{A}_P$, satisfying the conditions of Lemma 3.

In what follows, by the letter $a_i$ $(i = 1, 2, \ldots, n)$ we will mean the letter of the positive alphabet of the group $\mathfrak{A}_P$, which in this fixed sequence lies in place $i$. We will, however, continue to use the old designation for the pivot letter $p$.

**Remark 1′.** Suppose some group $\mathfrak{A}$ contains among its generators all letters from the sequence (5) and among its defining relations all defining relations of the group $\mathfrak{A}_P$, written in terms of these letters. If in the group $\mathfrak{A}$ the generators $a_j$ are equal to the identity for all odd $j < n$, then the generators $a_j$ are also equal to the identity for all even $j$, as well as $a_n$.

As the set of all $a_j$, for odd $j < n$, consists of all letters from the positive alphabet of $\mathfrak{A}_P$ not of type $l_{a,i}$ or $l_{a,i}^+$, this Remark 1′ follows directly from Remark 1 (see p. 23).

**Lemma 4.** *Let $a_{i_1}$ and $a_{i_2}$ be two heterotypic letters of the group $\mathfrak{A}_P$. If $a_{i_1}^{k_1} a_{i_2}^{k_2} = 1$ in $\mathfrak{A}_P$, then $k_1 = k_2 = 0$.*

*Proof.* There are two possible cases:

1. One of these two letters is the pivot letter $a_{n-2}$.
2. One of the letters belongs to the class $S$, and the other to $S^+$. For example, $a_{i_1} \in \overline{S}$ and $a_{i_2} \in \overline{S^+}$.

**Case 1.** Suppose that $i_1 = n - 2$. Then by Lemma 1, §1, we have $a_{i_1}^{k_1} = 1$ in $\mathfrak{A}_P$. By virtue of Theorem 1, we have $k_1 = 0$. Hence $a_{i_2}^{k_2} = 1$, and again using Theorem 1 we have $k_2 = 0$.

**Case 2.** $a_{i_1}^{k_1} \cdot a_{i_2}^{k_2} = 1$ in $\mathfrak{A}_P$. By Lemma 1 of this section, we have:

$$(a_{i_1}^{k_1} a_{i_2}^{k_2})_{\overline{S}} = a_{i_1}^{k_1} = 1 \quad \text{in } \mathfrak{A}_P, \quad \text{i.e. } k_1 = 0,$$

and

$$(a_{i_1}^{k_1} a_{i_2}^{k_2})_{\overline{S^+}} = a_{i_2}^{k_2} = 1 \quad \text{in } \mathfrak{A}_P, \quad \text{i.e. } k_2 = 0.$$

Lemma 4 is thereby proved. □

**Lemma 2′.** *If in the group $\mathfrak{A}_P$ we have the equality*

$$E_0 \equiv a_{n-1}^{s_t} a_n^{-1} a_{n-1}^{s_{t-1}} a_n^{-1} \cdots a_n^{-1} a_{n-1}^{s_1} a_n^{-1} a_{n-1}^{r_1} a_n a_{n-1}^{r_2} a_n \cdots a_{n-1}^{r_\nu} a_n = 1,$$

*where $a_n$ and $a_{n-1}$ are the last letters of the sequence (5), then $r_1 = 0$ and $s_t = 0$.*

*Proof.* As $a_n$ is a letter of type $l_{a,i}$, and $a_{n-1}$ is of type $l_{a,i}^+$, by Lemma 2 of this section it follows that the word $E_0$ is equal to the identity in the free group. Hence, as is known, the word $E_0$ can be transformed into 1 by a sequence of deletions. Let

$$E_0 \to E_1 \to \cdots \to 1$$

be such a sequence of transformations.

In this sequence, all letters $a_n$ and $a_n^{-1}$ must reduce with one another. In the word $E_0$, there are only negative exponents of the letter $a_n$ at first, and after the segment $a_{n-1}^{r_1}$ there are only positive. Hence, to start reducing the letters $a_n$ with $a_n^{-1}$, it must be the case that the word $a_{n-1}^{r_1}$ is empty, i.e. $r_1 = 0$. The rightmost letter $a_n$ must reduce with the leftmost letter $a_n^{-1}$. After this reduction, the word $E_0$ has been reduced to the word $a_{n-1}^{s_t}$. Consequently, $a_{n-1}^{s_t} = 1$ in the free group, i.e. $s_t = 0$. □

**Lemma 5.** *Let $E$ be a word of the group $\mathfrak{A}_P$ which does not contain the letter $p$. If*

$$Q_0 \equiv R p^\sigma E p^{-\sigma} S \to Q_1 \to Q_2 \to \cdots \to Q_m \tag{6}$$

*is a sequence of elementary transformations in $\mathfrak{A}_P$ without insertions of the letter $p$, in which the specified letter $p^\sigma$ and $p^{-\sigma}$ reduce with one another (and where $R$ and $S$ are arbitrary words of the group $\mathfrak{A}_P$), then*

$$p^\sigma E p^{-\sigma} = (E^{-1})^+ \text{ in the group } \mathfrak{A}_P.$$

*Proof.* Suppose $\sigma = 1$. Then the part of the sequence (6) from the beginning up to the point when the letters $p^\sigma$ and $p^{-\sigma}$ reduce with one another, has the form:

$$RpEp^{-1}S \to \cdots \to Q'_k pp^{-1} Q''_k. \tag{$6'$}$$

The specified letters $p$ and $p^{-1}$ appear in all words in $(6')$. Replace $(6')$ with its canonical sequence associated with the individuality of the distinguished letter $p$ (see p. 19),

$$RU_1 p U_1^+ E p^{-1} S \to \cdots \to Q'_k pp^{-1} Q''_k. \tag{$6''$}$$

In the sequence $(6')$ the distinguished letters $p$ and $p^{-1}$ are not affected by any transformation, as this is a canonical sequence for $p$, and $p^{-1}$ is not included in any of the defining relations of $\mathfrak{A}_P$. This shows that the segment, enclosed in the original word between $p$ and $p^{-1}$, is transformed in the sequence $(6')$ into the empty word, i.e.

$$U_1^+ E = 1 \text{ in } \mathfrak{A}_P.$$

Using this equality, and the equality $p = U_1 p U_1^+$, we find

$$pEp^{-1} = U_1 p U_1^+ E p^{-1} = U_1 pp^{-1} = (U_1^+)^+ = (U_1^+ E E^{-1})^+ = (E^{-1})^+ \text{ in } \mathfrak{A}_P.$$

The case of $\sigma = -1$ is treated analogously. $\qquad\square$

**Definition.** A finite sequence of integers

$$r_1, r_2, r_3, \ldots, r_m \tag{7}$$

will be said to be *length-reducible* if $r_{j_0} = 0$ for some $j_0$ with $1 < j_0 < m$. We will call the *reduction* at $r_{j_0}$ of the sequence (7) the new sequence

$$r_1 r_2, \ldots, r_{j_0-2}, r_{j_0-1} + r_{j_0+1}, r_{j_0+2}, \ldots, r_m, \tag{8}$$

whose number of entries is two less than that of the sequence (7).

A finite sequence of integers will be called *degenerate* if it can be reduced, by some finite number of reductions, to the sequence containing only a single zero. In particular, any sequence consisting only of zeroes is degenerate.

**Lemma 6.** *Let $A \equiv pXpX^+$, $B \equiv pYpY^+$, where $X$ and $Y$ are words written with letters from the class $S$, and $A \neq B$ in the group $\mathfrak{A}_P$. If at the same time the word*

$$Q \equiv a^{k_1}(AB^{-1})^{s_1} a^{k_2}(AB^{-1})^{s_2} \cdots a^{k_\mu}(AB^{-1})^{s_\mu} a^{k_{\mu+1}}, \tag{9}$$

*where $a$ is an arbitrary letter from the alphabet of the group $\mathfrak{A}_P$, is equal to 1 in $\mathfrak{A}_P$, then the corresponding sequence of*

$$k_1, s_1, k_2, s_2, \ldots, k_\mu, s_\mu, k_{\mu+1} \tag{10}$$

*either consists entirely of zeroes, or else is length-reducible.*

*Proof.* Note that any word $Q$ of type (9) unambiguously determines a sequence of type (10), and conversely, every sequence of type (10) unambiguously determines a word of type (9). In order for (10) to be length-reducible, by definition, it is necessary and sufficient that either some $s_i$ is equal to zero, or some $k_j$ for $1 < j \leq \mu$.

We begin the proof of Lemma 6.

Suppose $Q = 1$ in $\mathfrak{A}_P$.

If $\mu = 0$, then $Q \equiv a^{k_1} = 1$ in $\mathfrak{A}_P$. Then by Lemma 4, we have $k = 0$ and the sequence (10) consists of a single zero.

For $\mu \geq 1$, consider separately the following two cases:

**I.** The letter $a$ is distinct from $p$.

**II.** The letter $a$ coincides with $p$.

**Case I.** The word $Q$ has the form

$$a^{k_1}(pXpX^+Y^{+-1}p^{-1}Y^{-1}p^{-1})^{s_1}a^{k_2}\cdots a^{k_\mu}(pXpX^+Y^{+-1}p^{-1}Y^{-1}p^{-1})^{s_\mu}a^{k_{\mu+1}}.$$

We may assume that no $s_i$ is equal to zero, for otherwise the sequence (10) would already be length-reducible. As

$$(pXpX^+Y^{+-1}p^{-1}Y^{-1}p^{-1})^{s_i} = p(XpX^+Y^{+-1}p^{-1}Y^{-1})^{s_i}p^{-1}$$

in the free group, it follows that $Q$ is equal to the word

$$Q' \equiv a^{k_1}p(XpX^+Y^{+-1}p^{-1}Y^{-1})^{s_1}p^{-1}a^{k_2}\cdots$$
$$\cdots a^{k_\mu}p(XpX^+Y^{+-1}p^{-1}Y^{-1})^{s_\mu}p^{-1}a^{k_{\mu+1}}. \tag{11}$$

Now, the word $Q'$ is also equal to 1 in $\mathfrak{A}_P$. By virtue of Lemma 3 of §1, $Q'$ can be reduced to 1 by a sequence of elementary transformations

$$Q' \equiv Q_1 \to Q_2 \cdots \to Q_m \equiv 1 \tag{12}$$

without any insertions of the letter $p^\sigma$. As $\mu \geq 1$ and $s_i \neq 0$, the word $Q'$ contains a letter $p^\sigma$. Every individual letter $p^\sigma$ included in the word $Q'$ must cancel in the sequence (12). Consider the first pair $p$ and $p^{-1}$ which cancel in this sequence. There are two cases to consider.

In the **first case**, this pair may be a pair of letters which occur in some segment of the word $Q'$ of the form $(XpX^+Y^{+-1}p^{-1}Y^{-1})^{s_i}$.

In the **second case**, this may be a pair of letters $p$ and $p^{-1}$ which surround some segment $a^{k_j}$ in (11), where $1 < j \leq \mu$.

It is easy to see that no other cases are possible. Indeed, the letters $p^\sigma$ from different factors $(XpX^+Y^{+-1}p^{-1}Y^{-1})^{s_i}$ cannot cancel first, as between them there are obviously other occurrences of letters $p^\sigma$.

**The first case.** The first-deleted letters $p$ and $p^{-1}$ in the sequence (12) match in individuality with the same letters distinguished in the segment $(XpX^+Y^{+-1}p^{-1}Y^{-1})^{s_{i_0}}$ of $Q'$. Either $s_{i_0} > 0$, or else $s_{i_0} < 0$. As these two cases are entirely analogous, we consider only the first of the two.

Suppose $s_{i_0} > 0$. Then the word $Q'$ has the form

$$Q_1'a^{k_{i_0}}p(XpX^+Y^{+-1}p^{-1}Y^{-1}XpX^+\cdots XpX^+Y^{+-1}p^{-1}Y^{-1})p^{-1}a^{k_{i_0+1}}Q_2', \tag{11'}$$

and in the sequence (12), the first deletion between some letters $p$ and $p^{-1}$ occur among those which in (11') are enclosed in brackets.

Here again we are faced with two possibilities:

1) the considered deletion is right-sided,
2) the considered deletion is left-sided.

As each of these cases is analogous with the other, we will consider only the case of right-sided deletion. Then the first deletion is of two letters $p^{-1}$ and $p$ which in $Q'$ surround a segment $Y^{-1}X$. This segment does not contain the letter $p^\sigma$. By Lemma 5, in the group $\mathfrak{A}_P$ we have the equality

$$p^{-1}Y^{-1}Xp = [(Y^{-1}X)^{-1}]^+ = (X^{-1}Y)^+ = Y^+X^{+-1}.$$

From this we find

$$XpX^+ = YpY^+$$

and, as a consequence, $A = B$ in $\mathfrak{A}_P$, which contradicts the assumption of the lemma.

**The second case.** In the sequence (11) the first pair of deleted letters $p^{-1}$ and $p$ surround a segment $a^{k_{j_0}}$ in the word $Q'$, where $1 < j_0 \leq \mu$.

Again using Lemma 5, we find that in $\mathfrak{A}_P$ we have the equality

$$p^{-1}a^{k_{j_0}}p = (a^{-k_{j_0}})^+ = a^{+k_{j_0}}. \tag{13}$$

By assumption, the letter $a$ does not coincide with $p$.

Suppose $a \in S^+$.

By virtue of Lemma 3 of §1 we have a sequence of transformations in the group $\mathfrak{A}_P$, without insertions of $p^\sigma$:

$$p^{-1}a^{k_{j_0}}p \to \cdots a^{k_{j_0}} \to Z_i'p^{-1}pZ_i'' \to Z_i'Z_i'' \to \cdots \to a^{+k_{j_0}}. \tag{14}$$

The part of the sequence (14) from the beginning to the point of the deletion of the letters $p$ and $p^{-1}$ will have the form

$$p^{-1}a^{k_{j_0}}p \to \cdots \to Z_i'p^{-1}pZ_i''. \tag{14'}$$

The letters $p$ and $p^{-1}$ represented here appear in every word of the sequence (14'). Form the canonical sequence for (14') associated to the letter $p$:

$$p^{-1}a^{k_{j_0}}UpU^+ \to \cdots \to Z_i'p^{-1}pZ_i''.$$

In this sequence there are no transformations performed on the letters $p$ and $p^{-1}$. Therefore, we have an independent transformation $a^{k_{j_0}}U \to \cdots \to 1$ in the group $\mathfrak{A}_P$. As $a \in S^+$, and as the word $U$ does not contain letters of the alphabet $\overline{S^+}$ (see the correspondence $L_p$), it follows by Lemma 1 of the present section that we have

$$(a^{k_{j_0}}U)_{\overline{S^+}} \equiv a^{k_{j_0}} = 1 \text{ in the group } \mathfrak{A}_P.$$

This implies that $k_{j_0} = 0$ for some $1 < j_0 \leq \mu$, i.e. the sequence (10) is length-reducible.

Hence, in the case that $a \in S^+$, the lemma is proved.

If instead $a \in S$, then $a^+ \in S^+$. The equality (13) can be found regardless of whether $a \in S$ or $a \in S^+$. Hence we can use it. Replacing the occurrence of $p^{-1}a^{k_{j_0}}p$ in $Q'$ by $a^{+k_{j_0}}$, we obtain the word

$$Q'' \equiv a^{k_1}p(XpX^+Y^{+^{-1}}p^{-1}Y^{-1})^{s_1}p^{-1}a^{k_2} \cdots$$
$$\cdots p(XpX^+Y^{+^{-1}}p^{-1}Y^{-1})^{s_{j_0-1}}a^{+^{-k_{j_0}}}(XpX^+Y^{+^{-1}}p^{-1}Y^{-1})^{s_{j_0}} \cdots \tag{15}$$
$$\cdots p(XpX^+Y^{+^{-1}}p^{-1}Y^{-1})^{s_\mu}p^{-1}a^{k_{\mu+1}},$$

which is also equal to 1 in the group $\mathfrak{A}_P$. Take a sequence of elementary transformations without insertions of the letter $p$ (Lemma 3, §1, Chapter I):

$$Q'' \to Q_1'' \to Q_2'' \to \cdots Q_{m_1}'' \equiv 1. \tag{16}$$

Consider the pair of letter $p$ and $p^{-1}$, which is (16) cancel with one another first. There are now only three possible cases:

1) The letters come from one of the segments of type $(XpX^+Y^{+^{-1}}p^{-1}Y^{-1})^{s_i}$ in the word $Q''$.

2) The letters come from two distinct segments of type $(XpX^+Y^{+^{-1}}p^{-1}Y^{-1})^{s_i}$ between which there are no letters $p^\sigma$. As can be seen from the expression for $Q''$ in (15), in this case we have a deletion between the final of the letters $p^{-1}$ in the segment $(XpX^+Y^{+^{-1}}p^{-1}Y^{-1})^{s_{j_0-1}}$ and the first of the letters $p$ of the segment $(XpX^+Y^{+^{-1}}p^{-1}Y^{-1})^{s_{j_0}}$.

3) The first pair of deleted letters $p$ and $p^{-1}$ surround a segment $a^{k_{j_1}}$ in the word $Q''$, where $1 < j_1 \leq \mu$.

Case 1) easily leads to a contradiction with the assumptions of the lemma, in the same way as we did when considering the first deletion in the sequence (12).

In case 2), depending on the possible combinations of signs of the numbers $s_{j_0-1}$ and $s_{j_0}$, the first letter $p$ and $p^{-1}$ to be deleted will surround one of the following words:

$$\alpha)\ Y^{-1}a^{+-k_{j_0}}X, \ \beta)\ Y^{-1}a^{+-k_{j_0}}Y, \ \gamma)\ X^{-1}a^{+-k_{j_0}}X, \ \delta)\ X^{-1}a^{+-k_{j_0}}Y.$$

All these cases are pairwise analogous, and each leads to the equality $k_{j_0} = 0$. We will only consider the case $\alpha$). As the word $Y^{-1}a^{+-k_{j_0}}X$ does not contain the letter $p$, by Lemma 5 we have in the group $\mathfrak{A}_P$ the equality

$$p^{-1}Y^{-1}a^{+-k_{j_0}}Xp = [(Y^{-1}a^{+-k_{j_0}}X)^{-1}]^+.$$

Then by Lemma 3 of §1 we can find a sequence of elementary transformations without insertions of $p^\sigma$:

$$p^{-1}Y^{-1}a^{+-k_{j_0}}Xp \to \cdots \to [(Y^{-1}a^{+-k_{j_0}}X)^{-1}]^+.$$

The part of this sequence before the deletion of the letters $p$ and $p^{-1}$ has the form

$$p^{-1}Y^{-1}a^{+-k_{j_0}}Xp \to \cdots \to Z_1 p^{-1} p Z_2, \tag{17}$$

where the letters $p$ and $p^{-1}$ appear in all these words. We form the canonical sequence for (17) associated to the letter $p$:

$$p^{-1}Y^{-1}a^{+-k_{j_0}}XUpU^+ \to \cdots \to Z_1 p^{-1} p Z_2. \tag{17'}$$

Here we have an independent transformation

$$Y^{-1}a^{+-k_{j_0}}XU \to \cdots \to 1 \ \text{ in } \mathfrak{A}_P,$$

whence by Lemma 1, taking into account the fact that the words $X, U$ and $Y^{-1}$ do not contain letters from the alphabet $\overline{S^+}$, we find

$$(Y^{-1}a^{+-k_{j_0}}XU)_{\overline{S^+}} \equiv a^{+-k_{j_0}} = 1 \text{ in } \mathfrak{A}_P.$$

This, by Lemma 4, implies that $k_{j_0} = 0$. Hence the sequence (10) is length-reducible, for $1 < j_0 \leq \mu$, and the lemma is proved.

In case 3), we find a new equality, analogous to the equality (13),

$$p^{-1}a^{k_{j_1}}p = a^{+-k_{j_1}}.$$

Replacing in the word $Q''$ the segment $p^{-1}a^{k_{j_1}}p$ by $a^{+-k_{j_1}}$, we get a word $Q'''$ of the same type as $Q''$, and which is equal to 1 in $\mathfrak{A}_P$. It follows, that in the word $Q''$ the number of segments of the form $p^{-1}a^{k_j}p$ will be one less than in the word $Q''$. The segments of type $(XpX^+Y^{+-1}p^{-1}Y^{-1})^{s_i}$ in the word $Q''$ will all be preserved in the word $Q'''$.

Repeating for the word $Q'''$ the same reasoning which we have just carried out for $Q''$, we either find the conclusion of the lemma, or we find a word $Q^{\mathrm{IV}}$ of the same type, with an even smaller number of segments of the form $p^{-1}a^{k_j}p$, etc. As all segments of type $(XpX^+Y^{+-1}p^{-1}Y^{-1})^{s_i}$ are preserved by the above process, and as there are only finitely many segments of the form $p^{-1}a^{k_j}p$ in the word $Q'$, we will after a finite number

of steps find a word $Q^{(\lambda)}$, for which we have a sequence of elementary transformations in $\mathfrak{A}_P$, without insertions of the letter $p^\sigma$:

$$Q^{(\lambda)} \to Q_1^{(\lambda)} \to Q_2^{(\lambda)} \cdots \to 1,$$

and in this sequence the first-deleted letters $p^\sigma$ occur in a segment of the word $Q^{(\lambda)}$ of type $(XpX^+Y^{+\,-1}p^{-1}Y^{-1})^{s_i}$. As this case was treated above, this proves the lemma.

**Case II.** The letter $a$ coincides with $p$ and $\mu \geq 1$. Then $Q$ has the form

$$p^{k_1}(pXpX^+Y^{+\,-1}p^{-1}Y^{-1}p^{-1})^{s_1}p^{k_2}\cdots p^{k_\mu}(pXpX^+Y^{+\,-1}p^{-1}Y^{-1}p^{-1})^{s_\mu}p^{k_{\mu+1}}.$$

We may assume that none of the $s_i$ are equal to zero, as otherwise the sequence (10) will be length-reducible and there is nothing to prove. Again, replace the segments of type $(pXpX^+Y^{+\,-1}p^{-1}Y^{-1}p^{-1})^{s_1}$ by $p(XpX^+Y^{+\,-1}p^{-1}Y^{-1})^{s_1}p^{-1}$. Doing this transforms $Q$ into the word

$$Q' \equiv p^{k_1+1}(XpX^+Y^{+\,-1}p^{-1}Y^{-1})^{s_1}p^{k_2}(XpX^+Y^{+\,-1}p^{-1}Y^{-1})^{s_2}\cdots$$
$$\cdots p^{k_\mu}(XpX^+Y^{+\,-1}p^{-1}Y^{-1})^{s_\mu}p^{k_{\mu+1}-1}. \quad (18)$$

Hence, $Q' = 1$ in $\mathfrak{A}_P$. By Lemma 3 of §1, we have a sequence of elementary transformations without insertions of the letter $p^\sigma$:

$$Q' \to Q_1 \to Q_2 \to \cdots \to Q_m \equiv 1. \quad (19)$$

In this sequence of transformations some pair of letters $p$ and $p^{-1}$ appearing in $Q'$ are the first such pair to be deleted. They cannot both be taken from the same segment of type $(XpX^+Y^{+\,-1}p^{-1}Y^{-1})^{s_i}$, because, as shown above, this leads to the equality $A = B$ in $\mathfrak{A}_P$, which contradicts the assumptions of the lemma. Consequently, one of the first-deleted letters $p^\sigma$ is taken from some segment $p^{k_{j_0}}$ of the word $Q'$, where $k_{j_0} \neq 0$. The other will be taken from a segment $(XpX^+Y^{+\,-1}p^{-1}Y^{-1})^{s_{i_0}}$ adjacent to $p^{k_{j_0}}$ (on the right or left). As the cases $k_{j_0} > 0$ and $k_{j_0} < 0$ are analogous, we will consider only the first. Suppose $k_{j_0} > 0$. Then the letter $p^{-1}$ is taken from the segment $(XpX^+Y^{+\,-1}p^{-1}Y^{-1})^{s_{i_0}}$. Recall that the deletion under consideration is the first occurring in the sequence (19). Hence, between the deletion of the letter $p$ and $p^{-1}$ in the word $Q'$ there can be no other letters $p^\sigma$. Hence, the segment $(XpX^+Y^{+\,-1}p^{-1}Y^{-1})^{s_{i_0}}$ must be that segment of the same type immediately adjacent to $p^{k_{j_0}}$ on the left, since the left end of the segment $(XpX^+Y^{+\,-1}p^{-1}Y^{-1})^{s_i}$ always contains a positive letter $p$.

Now there are two possible cases:

1) $s_{i_0} > 0$   and   2) $s_{i_0} < 0$.

Suppose $s_{i_0} > 0$. Then $(XpX^+Y^{+\,-1}p^{-1}Y^{-1})^{s_{i_0}}$ has the form

$$(XpX^+Y^{+\,-1}p^{-1}Y^{-1}XpX^+Y^{+\,-1}p^{-1}Y^{-1}\cdots XpX^+Y^{+\,-1}p^{-1}Y^{-1}),$$

and the word $Q'$ is hence written as:

$$Q_1'p^{k_{j_0}-1}(XpX^+Y^{+\,-1}p^{-1}Y^{-1}XpX^+\cdots XpX^+Y^{+\,-1}\underline{p^{-1}}Y^{-1})\times$$
$$\times \underline{p}\,p^{k_{j_0}}(XpX^+Y^{+\,-1}p^{-1}Y^{-1})^{s_{i_0}+1}Q_2', \quad (20)$$

where $Q_1'$ and $Q_2'$ are the same type of words as $Q'$.

In (20) we have underlined the individuality of the letters $p^{-1}$ and $p$ which delete with each other in the sequence (19).

By Lemma 5 we have the equality

$$p^{-1}Y^{-1}p = [(Y^{-1})^{-1}]^+ = Y^+,$$

in the group $\mathfrak{A}_P$, from which we have

$$p = YpY^+. \tag{21}$$

In the analogous way, from $s_{i_0} < 0$ we arrive at the equality

$$p = XpX^+. \tag{22}$$

Hence, it remains to consider the case when we either have the equality (21) or the equality (22) in the group $\mathfrak{A}_P$. As the considerations required in either case is essentially those of the other, we will only consider the case when we have the equality (21).

From (21) we deduce, in $\mathfrak{A}_P$, the equality: $p^{-1} = Y^{+-1}p^{-1}Y^{-1}$, from which we find $XpX^+Y^{+-1}p^{-1}Y^{-1} = XpX^+p^{-1}$. The word $Q'$ (see (18)) is equal in the group to the word

$$Q'_0 \equiv p^{k_1+1}(XpX^+p^{-1})^{s_1}p^{k_2}(XpX^+p^{-1})^{s_2}\cdots p^{k_\mu}(XpX^+p^{-1})^{s_\mu}p^{k_{\mu+1}-1}. \tag{23}$$

Here the exponents $s_i$ and $k_j$ are the same as in (18). Hence, if some $s_i$ or $k_j$ for $1 < j \leq \mu$ are equal to zero, then the sequence (10) is length-reducible.

The word $Q'_0$ is equal to 1 in $\mathfrak{A}_P$. Generally speaking, it can contain two adjacent letters $p$ and $p^{-1}$. We will now find a word $Q''_0$, equal in the free group to $Q'_0$, and which contains no two directly adjacent letters $p$ and $p^{-1}$. To write this word $Q''_0$ in the general form, we introduce the following functions:

$$\text{sign}(n) = \left\{ \begin{array}{l} 1 \text{ when } n > 0, \\ 0 \text{ when } n = 0, \\ -1 \text{ when } n < 0, \end{array} \right.$$

$$\beta_1(n) = \left\{ \begin{array}{l} 0 \text{ when } n \leq 0, \\ -1 \text{ when } n > 0, \end{array} \right. \qquad \beta_2(n) = \left\{ \begin{array}{l} 0 \text{ when } n \geq 0, \\ 1 \text{ when } n < 0. \end{array} \right.$$

With the help of these functions, every segment of $Q'_0$ having the form $(XpX^+p^{-1})^{s_i}$ can be written in the form:

$$p^{\beta_2(\delta_i)}(XpX^+p^{-1}XpX^+p^{-1}\cdots XpX^+)^{\delta_i}p^{\beta_1(\delta_i)},$$

where $\delta_i = \text{sign}(s_i)$, and in the word whose exponent is $\delta_i$ the segment $XpX^+$ is repeated $|s_i|$ times. Recall that we earlier defined $A^t$ (where $A$ is an arbitrary word) to be the empty word when $t = 0$.

After deleting the adjacent letters $p^\sigma$, the word $Q'_0$ will be transformed into:

$$Q''_0 \equiv p^{k_1+1+\beta_2(\delta_1)}(XpX^+p^{-1}XpX^+p^{-1}\cdots XpX^+)^{\delta_1}$$
$$p^{k'_2}(XpX^+p^{-1}XpX^+p^{-1}\cdots XpX^+)^{\delta_2}p^{k'_3}\cdots$$
$$\cdots p^{k'_i}(XpX^+p^{-1}XpX^+p^{-1}\cdots XpX^+)^{\delta_i}p^{k'_{i+1}}\cdots$$
$$p^{k'_\mu}(XpX^+p^{-1}XpX^+p^{-1}\cdots XpX^+)^{\delta_\mu}p^{k_{\mu+1}-1+\beta_1(s_\mu)}, \tag{24}$$

where $k'_i = k_i + \beta_1(\delta_{i-1}) + \beta_2(\delta_i)$ for $1 < i \leq \mu$, and $\delta_i = \text{sign}(i) = \pm 1$.

Note that the numbers $k'_i$ are not all required to be different from zero. However, if $\delta_{i-1} = -\delta_i$, then $k'_i = k_i$ and is hence non-zero. Indeed, if in this case $\delta_i = 1$, then $k'_i = k_i + \beta_1(-1) + \beta_2(1) = k_i + 0 + 0 = k_i$. If instead $\delta_i = -1$, then $k'_i = k_i + \beta_1(1) + \beta_2(-1) = k_i + (-1) + 1 = k_i$. Now, the word $Q''_0$ is also equal to 1 in $\mathfrak{A}_P$. We have a sequence of elementary transformations without insertions of the letter $p^\sigma$:

$$Q''_0 \to Q''_1 \to Q''_2 \to \cdots \to Q''_\mu = 1. \tag{25}$$

In this sequence, there is a first pair of letters $p$ and $p^{-1}$ to delete with one another. Each of these letters appear in the word $Q_0''$ either in a segment of type $p^{k_j'}$, or else within a segment of type

$$(XpX^+p^{-1}XpX^+p^{-1}\cdots XpX^+)^{\delta_i}. \tag{26}$$

Both of these letters cannot come from a segment of type $p^{k_j'}$, as then between them there is a letter $p^\sigma$ appearing inside some segment of type (26). We hence have three possible cases:

1) The letters $p$ and $p^{-1}$ come from a single segment of type (26).
2) The letters come from distinct segments of type (26).
3) One of the letters comes from a segment of type (26), and the other from a segment of type $p^{k_j'}$.

**Case 1).** Suppose that the letters both come from the segment

$$(XpX^+p^{-1}XpX^+p^{-1}\cdots XpX^+)^{\delta_{i_0}}.$$

Either $\delta_{i_0} = 1$, or else $\delta_{i_0} = -1$. As these cases are analogous, we will only consider the case when $\delta_{i_0} = 1$. The deletion under consideration is either right-sided, or else left-sided. These cases are also analogous.

Hence we consider the case when the deletion is right-sided. Then the deleted letters $p$ and $p^{-1}$ surround in $Q_0''$ a word $X$. By Lemma 5, we have in $\mathfrak{A}_P$ the equality

$$p^{-1}Xp = (X^{-1})^+ = X^{+^{-1}},$$

whence

$$XpX^+ = p. \tag{27}$$

If we were instead considering the case when $\delta_{i_0} = -1$, then we would in the analogous way arrive at the equality (27).

The equality (27) together with (21) yields

$$A = pXpX^+ = pp = pYpY^+ = B,$$

which contradicts the assumptions of the lemma.

**Case 2).** Suppose that the first-deleted letters $p$ and $p^{-1}$ appear, respectively, in the segments

$$(XpX^+p^{-1}XpX^+p^{-1}\cdots XpX^+)^{\delta_{i_1}} \text{ and } (XpX^+p^{-1}XpX^+p^{-1}\cdots XpX^+)^{\delta_{i_2}}.$$

We may assume without loss of generality that $i_1 < i_2$. If in this case $i_1 + 1 < i_2$, then in the word $Q_0''$ between the distinguished letters there is some whole segment of type (26), and so these letters could not be the first to cancel; hence, we have that $i_2 = i_1 + 1$. As the deleted letters $p^\sigma$ must appear on the very edges of their respective segments, we necessarily have $\delta_{i_1} = -\delta_{i_2}$. Because of this, the segment $p^{k_{i_2}'}$ enclosed between them must be empty, i.e. $k_{i_2}' = 0$. But as shown above (p. 31) from $\delta_{i_1} = +\delta_{i_2-1} = -\delta_{i_2}$ we have $k_{i_2}' = k_{i_2} = 0$. As the segment $p^{k_{i_2}'}$ has segments of type (26) to the left and to the right in the word $Q_0''$, it follows that $1 < i_2 \le \mu$. Hence, the sequence (10) must be length-reducible.

**Case 3).** One of the letters comes from a segment of type (26), and the other from a segment of type $p^{k_j'}$. As the specified deletion is the first to occur, these segments must be adjacent to one another in $Q_0''$. Suppose the segment $p^{k_j'}$ lies to the right. If it lies to

the left, then this case is treated analogously, and is omitted. We distinguish in $Q_0''$ the above segments as

$$Q_0'' \equiv Q_{0,1}''(XpX^+p^{-1}XpX^+p^{-1}\cdots XpX^+)^{\delta_{j-1}}p^{k_j'}Q_{0,2}''.$$

Either $\delta_{j-1} = 1$, or else $\delta_{j-1} = -1$. Suppose $\delta_{j-1} = 1$. Then $k_j' < 0$ and the deleted letters surround the word $X^+$. By Lemma 5, we have in $\mathfrak{A}_P$ the equality

$$XpX^+ = p. \tag{27}$$

The same conclusion can be found when $\delta_{j-1} = -1$.

However, as we saw above, this leads to a contradiction with the assumptions of the lemma. Hence, the proof of Lemma 6 is complete. $\qquad\square$

**Corollary.** *Suppose $A$ and $B$ are words of the form given in Lemma 6. In order for the word*

$$Q \equiv a^{k_1}(AB^{-1})^{s_1}a^{k_2}\cdots a^{k_\mu}(AB^{-1})^{s_\mu}a^{k_{\mu+1}} \tag{28}$$

*to be equal to the identity in $\mathfrak{A}_P$, it is necessary and sufficient that the sequence*

$$k_1, s_1, k_2, s_2, \ldots, k_\mu, s_\mu, k_{\mu+1} \tag{10}$$

*is degenerate.*

*Proof.* **Necessity.** Suppose $Q = 1$ in $\mathfrak{A}_P$. If the sequence (10) does not contain only zeroes, then by Lemma 6 it is length-reducible. This implies, that either $s_i = 0$ for some $i$, or else $k_j = 0$ for some $1 < j \leq \mu$. Suppose for definiteness that $s_{i_0} = 0$. Reducing the sequence (10), we find a new sequence

$$k_1, s_1, k_2, \ldots, s_{i_0-1}, k_{i_0} + k_{i_0+1}, s_{i_0+1}, \ldots, k_\mu, s_\mu, k_{\mu+1}. \tag{10$'$}$$

The sequence (10$'$) uniquely defines a word $Q_1$ of type (28). As in $Q$ the segment $(AB^{-1})^{s_{i_0}}$ was the empty word, we have $Q_1 = Q$ in the free group. This implies that $Q_1 = 1$ in $\mathfrak{A}_P$. The number of elements in the sequence (10$'$) is strictly less than the corresponding number in the sequence (10). If (10$'$) does not only consist of zeroes, then, applying Lemma 6 again, we find that (10$'$) is length-reducible, etc. As the number of elements of the sequences thus obtained strictly decreases, in a finite number of steps we end in a sequence consisting of only zeroes. Hence, the sequence (10) is degenerate.

**Sufficiency.** Suppose the sequence (10) is degenerate. Then, length-reducing it accordingly, we find a finite system of sequences, each of which is obtained from another by a length-reduction; the last length-reduction produces a sequence consisting only of zeroes. Denote the sequence (10) by $\lambda_1$, and the sequence obtained from $\lambda_1$ by a length-reduction by $\lambda_2$, etc. Our system of sequences is then written as:

$$\lambda_1, \lambda_2, \ldots, \lambda_t, \tag{29}$$

where $\lambda_t$ is a sequence of zeroes.

Every sequence $\lambda_i$ in the system (29) uniquely defines a word $Q_i$ of type (28), where $Q_1$ is simply our word $Q$, and $Q_t$ is the empty word. The sequence $\lambda_{i+1}$ is obtained from $\lambda_i$ by a single length-reduction. This implies that $Q_i = Q_{i+1}$ in the free group. By transitivity, $Q = 1$ in the free group, and consequently also in $\mathfrak{A}_P$. $\qquad\square$

**Lemma 7.** *Let $A = pXpX^+$ and $B = pYpY^+$, where $X$ and $Y$ are words over the alphabet $S$, and $A \neq B$ in $\mathfrak{A}_P$. Let further $C$ and $E$ be arbitrary words from the group $\mathfrak{A}_P$ not containing the letter $p^\sigma$. If the word*

$$Q \equiv BA^{-1}EAB^{-1}C$$

*is equal to 1 in $\mathfrak{A}_P$, then in $\mathfrak{A}_P$ we also have the equalities:*

$$1) \ (E)_{\overline{S}} = 1, \quad 2) \ (E)_{\overline{S^+}} = 1, \quad 3) \ C = 1,$$

*where $(E)_{\overline{S}}$ and $(E)_{\overline{S^+}}$ are the projections of the word $E$ onto the alphabets $\overline{S}$ and $\overline{S^+}$.*

*Proof.* Suppose the word

$$Q \equiv pYpY^+X^{+^{-1}}p^{-1}X^{-1}p^{-1}EpXpX^+Y^{+^{-1}}p^{-1}Y^{-1}p^{-1}C$$

is equal to 1 in $\mathfrak{A}_P$. By Lemma 3 of §1, we can find a sequence of elementary transformations

$$Q \to Q_1 \to Q_2 \to \cdots \to 1 \tag{30}$$

without insertions of the letter $p^\sigma$.

Consider the first deleted letter $p^\sigma$ in this sequence. This deletion is either right-sided, or left-sided.

Suppose it is left-sided. Then the deleted individuality of $p^\sigma$ in the word $Q$ surrounds either a segment $Y^+X^{+^{-1}}$, or else a factor $X^+Y^{+^{-1}}$. Let us say that it surrounds the segment $Y^+X^{+^{-1}}$. The other case is entirely analogous. By Lemma 5, we have in $\mathfrak{A}_P$ the equality

$$pY^+X^{+^{-1}}p^{-1} = [(Y^+X^{+^{-1}})^{-1}]^+ = (X^+Y^{+^{-1}})^+ = Y^{-1}X,$$

whence

$$YpY^+ = XpX^+$$

and

$$B = pYpY^+ = pXpX^+ = A,$$

which contradicts the assumptions of the lemma.

Suppose now that the deletion is right-sided. Then the deleted letters surround in $Q$ a segment $E$. By Lemma 5, in the group $\mathfrak{A}_P$ we have

$$p^{-1}Ep = (E^{-1})^+ = E^{+^{-1}}. \tag{31}$$

In the group $\mathfrak{A}_P$ this gives a sequence of elementary transformations without insertions of $p^\sigma$:

$$p^{-1}Ep \to \cdots \to E^{+^{-1}}.$$

The part of this sequence from the beginning to the moment of deletion of the letters $p$ and $p^{-1}$ has the form

$$p^{-1}Ep \to \cdots \to Z_1p^{-1}pZ_2. \tag{32}$$

By Lemma 2 of §1 we form the canonical sequence for (32)

$$p^{-1}EU_1pU_1^+ \to \cdots \to Z_1p^{-1}pZ_2, \tag{32'}$$

where no transformation affects the letters $p$ and $p^{-1}$.

From the sequence (32') we have

$$EU_1 = 1 \text{ in } \mathfrak{A}_P,$$

whence follows the equality

$$U_1 = E^{-1}. \tag{33}$$

From the equalities (31) and (33) we find

$$p^{-1}Ep = (E^{-1})^+ = U_1^+.$$

Replace in the word $Q$ the fixed segment $p^{-1}Ep$ by the word $U_1^+$, and denote the resulting word by $Q'$. In the group $\mathfrak{A}_P$, we find the equality

$$Q = pYpY^+X^{+^{-1}}p^{-1}X^{-1}U_1^+XpX^+Y^{+^{-1}}p^{-1}Y^{-1}C \equiv Q'.$$

The word $Q'$ is equal to 1 in $\mathfrak{A}_P$.

By considering a sequence of elementary transformations $Q' \to \cdots \to 1$, without insertions of the letter $p^\sigma$, and the first deleted letters $p^\sigma$ in this sequence, we may, analogously to how we did it for the word $Q$, find the equality

$$p^{-1}X^{-1}U_1^+Xp = [(X^{-1}U_1^+X)^{-1}]^+.$$

Consider further the sequence of transformations without insertions of $p^\sigma$,

$$p^{-1}X^{-1}U_1^+Xp \to \cdots \to [(X^{-1}U_1^+X)^{-1}]^+,$$

and using the corresponding canonical sequence of transformations, we find that for some $U_2$ we have the equality

$$X^{-1}U_1^+XU_2 = 1 \tag{34}$$

in the group $\mathfrak{A}_P$. As the words $U_1$ and $U_2$ contain only letters from the alphabet $\overline{S}$, from (34) and based on Lemma 1 we have

$$(X^{-1}U_1^+XU_2)_{\overline{S^+}} \equiv U_1^+ = 1 \quad \text{in } \mathfrak{A}_P.$$

By virtue of the symmetry of the group $\mathfrak{A}_P$ (see p. 21), the words 1 and $(U_1^+)^+ = U_1$ are equal. From (33) we have

$$E^{-1} = 1, \quad \text{i.e. } E = 1.$$

By Lemma 1, §2, Chapter I we find

$$(E)_{\overline{S}} = 1 \quad \text{and} \quad (E)_{\overline{S^+}} = 1 \quad \text{in } \mathfrak{A}_P.$$

Replacing in $Q$ the word $E$ by the empty word, we have

$$Q = BA^{-1}AB^{-1}C = C = 1$$

in the group $\mathfrak{A}_P$. Lemma 7 is proved.     □

# Chapter II
# Unsolvability of Some Algorithmic Problems in the Theory of Groups

## §1. The system of groups $\mathfrak{A}_{q,A,B}$ and $\mathfrak{A}_{q_i,A,B}$

For every pair of words $A$ and $B$ of the form $pXpX^+$, where $p$ is the pivot letter of the group $\mathfrak{A}_P$, and $X$ is any word written with letters from the class $S$, we will define a group denoted by $\mathfrak{A}_{q,A,B}$. For fixed words $A$ and $B$, the group $\mathfrak{A}_{q,A,B}$ will be given by a finite alphabet and a finite system of defining relations.

Of course, to specify the generating alphabet of any group, it is enough to specify its positive part.

The positive alphabet of the group $\mathfrak{A}_{q,A,B}$ will consist of all letters from the positive alphabet of the group $\mathfrak{A}_P$:

$$a_1, a_2, \ldots, a_n \tag{1}$$

together with $n$ new letters

$$q_1, q_2, \ldots, q_n, \tag{2}$$

and where the sequence (1) satisfies the conditions of Lemma 3, §2, Chapter I.

Suppose that

$$A_j = B_j \qquad (j = 1, 2, \ldots, m) \tag{3}$$

is a system of defining relations for the group $\mathfrak{A}_P$, written over the alphabet (1). The system of defining relations for the group $\mathfrak{A}_{q,A,B}$ is obtained from (3) by adding the new relations:

$$a_i q_i = q_i A B^{-1} \qquad (i = 1, 3, 5, \ldots, n-2), \tag{4}$$

$$q_{i-1} q_i = q_i a_{i-1} \qquad (1 < i \le n), \tag{5}$$

$$a_n q_n = q_n a_n q_n^{-1} B A^{-1}. \tag{6}$$

The letters $q_n^\sigma$ appear only in two relations of the group $\mathfrak{A}_{q,A,B}$: in the relation (6), and in the relation

$$q_{n-1} q_n = q_n a_{n-1}. \tag{5'}$$

Delete from the alphabet of the group $\mathfrak{A}_{q,A,B}$ the letters $q_n$ and $q_n^{-1}$, and delete from the system of defining relations those relations containing the letters $q_n^\sigma$. The group obtained in this manner will be denoted $\mathfrak{A}_{q_{n-1},A,B}$. The letter $q_{n-1}$ is a filtering letter in the group $\mathfrak{A}_{q_{n-1},A,B}$, as it appears only in a single defining relation

$$q_{n-2} q_{n-1} = q_{n-1} a_{n-2}. \tag{7}$$

The base $\overline{\mathfrak{A}}_{q_{n-1},A,B}^{(q_{n-1})}$ of the group $\mathfrak{A}_{q_{n-1},A,B}$ with respect to the filtering letter $q_{n-2}$ will be denoted by $\mathfrak{A}_{q_{n-2},A,B}$. In this group the letter $q_{n-2}$ appears only in the relations (4) and (5) for $i = n-2$. It is hence a filtering letter of the group $\mathfrak{A}_{q_{n-2},A,B}$. The base of the group $\mathfrak{A}_{q_{n-2},A,B}$ with respect to the letter $q_{n-2}$ is denoted by $\mathfrak{A}_{q_{n-3},A,B}$, etc. Hence, for a fixed pair of words $(A, B)$ we have defined the groups $\mathfrak{A}_{q_i,A,B}$ $(i = 1, 2, \ldots, n-1)$ in the following way:

$$1.\ \overline{\mathfrak{A}}_{q_1,A,B}^{(q_1)} \approx \mathfrak{A}_P, \qquad 2.\ \overline{\mathfrak{A}}_{q_i,A,B}^{(q_i)} \approx \mathfrak{A}_{q_{i-1},A,B} \quad (1 < i \le n-1). \tag{8}$$

The sign $\approx$ here signifies that we have two different designations of the same group. The alphabet of the group $\mathfrak{A}_P$ is contained in the alphabet of every group $\mathfrak{A}_{q_i,A,B}$. In

the group $\mathfrak{A}_{q_i,A,B}$ the letter $q_i$ is a filtering letter. The correspondence $L_{q_1}$ in the group $\mathfrak{A}_{q_1,A,B}$ has the form

$$a_1^k \underline{\hspace{2cm}} (AB^{-1})^k. \tag{9}$$

The correspondence $L_{q_i}$ in the group $\mathfrak{A}_{q_i,A,B}$ $(1 < i \leq n-1)$ for even $i$ will be

$$q_{i-1}^k \underline{\hspace{2cm}} a_{i-1}^k, \tag{10}$$

and for odd $i$ will be

$$q_{i-1}^{r_1} a_i^{s_1} q_{i-1}^{r_2} a_i^{s_2} \cdots q_{i-1}^{r_\mu} a_i^{s_\mu} q_{i-1}^{r_{\mu+1}} \underline{\hspace{2cm}}$$
$$\underline{\hspace{2cm}} a_{i-1}^{r_1} (AB^{-1})^{s_1} a_{i-1}^{r_2} (AB^{-1})^{s_2} \cdots a_{i-1}^{r_\mu} (AB^{-1})^{s_\mu} a_{i-1}^{r_{\mu+1}}. \tag{11}$$

**Lemma 1.** *Let $A = pXpX^+$, $B = pYpY^+$, where $X$ and $Y$ are arbitrary words over the alphabet $S$, and $A \neq B$ in the group $\mathfrak{A}_P$. Then for every group $\mathfrak{A}_{q_i,A,B}$ $(1 \leq i \leq n-1)$ the following holds:*

    I. *The letter $q_i$ is a regular filtering letter of the group $\mathfrak{A}_{q_i,A,B}$.*

    II. *If $E$ is a word consisting only of letters from alphabet of the group $\mathfrak{A}_P$, then the equality $E = 1$ in $\mathfrak{A}_{q_i,A,B}$ is equivalent to the equality $E = 1$ in $\mathfrak{A}_P$.*

*Proof.* The proof is by induction on the index $i$.

Suppose $i = 1$. The letter $q_1$ will be a regular filtering letter in $\mathfrak{A}_{q_1,A,B}$ if the correspondence (9) is an isomorphism in the group $\mathfrak{A}_P$. To show this, it suffices to show that the equality to 1 of an arbitrary word on the left-hand side of (9) is equivalent to the equality to 1 in $\mathfrak{A}_P$ of the corresponding word on the right-hand side. Suppose $a_1^k = 1$ in $\mathfrak{A}_P$. Then by Lemma 4, §2, Chapter I we have $k = 0$, whence $(AB^{-1})^k = 1$ in $\mathfrak{A}_P$. Suppose $(AB^{-1})^k = 1$ in $\mathfrak{A}_P$. As by assumption $A \neq B$ in $\mathfrak{A}_P$, by Lemma 6, §2, Chapter I we have $k = 0$, i.e. $a_1^k = 1$ in $\mathfrak{A}_P$. Hence **I** is proved for $i = 1$.

We now prove part **II** for $i = 1$. Suppose $E$ is a word over the letters of the alphabet of the group $\mathfrak{A}_P$. If $E = 1$, then obviously $E = 1$ also in $\mathfrak{A}_{q_1,A,B}$, as all defining relations of the group $\mathfrak{A}_P$ are also relations of the group $\mathfrak{A}_{q_1,A,B}$.

Suppose $E = 1$ in $\mathfrak{A}_{q_1,A,B}$. As proved, $q_1$ is a regular filtering letter in $\mathfrak{A}_{q_1,A,B}$, and so by the corollary of Lemma 3, §1, Chapter I we have $E = 1$ in $\overline{\mathfrak{A}}_{q_1,A,B}^{(q_1)}$, i.e. in $\mathfrak{A}_P$.

Suppose now that the lemma has already been proved for $i < i_0$.

We will now prove it for $i = i_0$. Part **I** is proved separately in two cases:

    1) $i_0$ is even,
    2) $i_0$ is odd.

**Case 1).** Suppose $i_0$ is even. Then the correspondence $L_{q_{i_0}}$ has the form (10). We must prove that this is an isomorphism in the group $\mathfrak{A}_{i_0-1,A,B}$. Suppose that $q_{i_0-1}^k = 1$ in $\mathfrak{A}_{q_{i_0-1},A,B}$. As $q_{i_0-1}$ is a filtering letter in $\mathfrak{A}_{q_{i_0-1},A,B}$, by Lemma 1, §1, Chapter I we have $k = 0$. Hence the word $a_{i_0-1}^k$ is also equal to 1 in the group $\mathfrak{A}_{q_{i_0-1},A,B}$.

Suppose now that $a_{i_0-1}^k = 1$ in $\mathfrak{A}_{q_{i_0-1},A,B}$. By the inductive hypothesis this implies that $a_{i_0-1}^k = 1$ in $\mathfrak{A}_P$, whence by Lemma 4, §2, Chapter I we have $k = 0$. Hence obviously also the word $q_{i_0-1}^k$ is equal to 1 in $\mathfrak{A}_{q_{i_0-1},A,B}$.

**Case II).** Suppose $i_0$ is odd. The correspondence $L_{q_{i_0}}$ has the form

$$U \equiv q_{i_0-1}^{r_1} a_{i_0}^{s_1} q_{i_0-1}^{r_2} a_{i_0}^{s_2} \cdots q_{i_0-1}^{r_\mu} a_{i_0}^{s_\mu} q_{i_0-1}^{r_{\mu+1}} \underline{\hspace{2cm}}$$
$$\underline{\hspace{2cm}} a_{i_0-1}^{r_1} (AB^{-1})^{s_1} a_{i_0-1}^{r_2} (AB^{-1})^{s_2} \cdots a_{i_0-1}^{r_\mu} (AB^{-1})^{s_\mu} a_{i_0-1}^{r_{\mu+1}} \equiv W. \tag{12}$$

We denote words appearing on the left-hand side of this correspondence by $U$ (with indices), and the corresponding word on the right-hand side by $W$ (indexed identically). Every word $W$ is a word of type (28) from §2, Chapter I. It uniquely defines a sequence of integers

$$r_1, s_1, r_2, s_2, \ldots, r_\mu, s_\mu, r_{\mu+1}. \tag{13}$$

The corresponding word $U$ is uniquely defined from the sequence (13).

We need to prove that $L_{q_{i_0}}$ is an isomorphism in the group $\mathfrak{A}_{q_{i_0-1}, A, B}$, i.e. that $W = 1$ in $\mathfrak{A}_{q_{i_0-1}, A, B}$ if and only if the corresponding word $U$ is equal to 1 in $\mathfrak{A}_{q_{i_0-1}, A, B}$. As the word $W$ contains only letters from the alphabet of the group $\mathfrak{A}_P$, by the inductive hypothesis we have that $W = 1$ in $\mathfrak{A}_{q_{i_0-1}, A, B}$ if and only if $W = 1$ in $\mathfrak{A}_P$. By assumption $A \neq B$ in $\mathfrak{A}_P$. Hence, by the corollary of Lemma 6, §2, Chapter I, it follows that the equality $W = 1$ holds in the group $\mathfrak{A}_P$ if and only if the sequence (13) corresponding to $W$ is degenerate. We will hence have proved that the correspondence $L_{q_{i_0}}$ is an isomorphism, once we have proved the following claim:

*For the word*

$$U \equiv q_{i_0-1}^{r_1} a_{i_0}^{s_1} q_{i_0-1}^{r_2} a_{i_0}^{s_2} \cdots q_{i_0-1}^{r_\mu} a_{i_0}^{s_\mu} q_{i_0-1}^{r_{\mu+1}} \tag{14}$$

*to be equal to the identity in 1 in the group $\mathfrak{A}_{q_{i_0-1}, A, B}$, it is necessary and sufficient for the corresponding sequence (13) to be degenerate.*

**Sufficiency.** Suppose that the sequence (13) is degenerate. This means that we can build a system of sequences of integers

$$\lambda_1, \lambda_2, \ldots, \lambda_t,$$

where $\lambda_1$ is the sequence (13), $\lambda_t$ consists only of zeroes, and $\lambda_{i+1}$ is obtained for every $i$ from $\lambda_i$ by one step of length-reduction. These sequences give rise to a corresponding system of words of the form (14)

$$U_1, U_2, \ldots, U_t,$$

where $U_1 \equiv U$, $U_t$ is the empty word, and every $U_{i+1}$ is equal to $U_i$ in the free group. Consequently, $U = 1$ in the free group, and hence also in $\mathfrak{A}_{q_{i_0-1}, A, B}$.

**Necessity.** Suppose $U = 1$ in $\mathfrak{A}_{q_{i_0-1}, A, B}$. We will prove first that the sequence (13) either consists only of zeroes, or else is length-reducible. If $\mu = 0$, then $u \equiv q_{i_0-1}^{r_1} = 1$ in $\mathfrak{A}_{q_{i_0-1}, A, B}$. As $q_{i_0-1}$ is a filtering letter in $\mathfrak{A}_{q_{i_0-1}, A, B}$, we have $r_1 = 0$, and the sequence (13) consists of a single zero.

Suppose $\mu \geq 1$. To prove the length-reducibility of the sequence (13), we must show that either $s_i = 0$ for some $i$, or else $r_j = 0$ for some $j$ satisfying $1 < j \leq \mu$.

If the word $U$ does not contain the letter $q_{i_0-1}^\sigma$, then all $r_j$ are equal to zero. Then for $\mu > 1$ the sequence (13) is length-reducible, and for $\mu = 1$ we have

$$U = a_{i_0}^{s_1} = 1 \quad \text{in} \quad \mathfrak{A}_{q_{i_0-1}, A, B}.$$

By the inductive hypothesis (part **II**) this implies that $a_{i_0}^{s_1} = 1$ in $\mathfrak{A}_P$. Consequently by Lemma 4, §2, Chapter I we have $s_1 = 0$. The sequence (13) hence only consists of zeroes.

The only remaining case is when $\mu \geq 1$ and the word $U$ contains some occurrence of the letter $q_{i_0-1}$. As $q_{i_0-1}$ is, by the inductive hypothesis, a regular filtering letter of the group $\mathfrak{A}_{q_{i_0-1}, A, B}$, by Lemma 3, §2, Chapter I we have a sequence of elementary transformations, without insertions of the letter $q_{i_0-1}^\sigma$, in $\mathfrak{A}_{q_{i_0-1}, A, B}$:

$$U \rightarrow \cdots \rightarrow 1. \tag{15}$$

In this sequence, every letter $q_{i_0-1}^\sigma$ contained in the word $U$ must eventually be deleted. Consider the first pair of letters $q_{i_0-1}^\sigma$ and $q_{i_0-1}^{-\sigma}$ to delete with one another in the sequence (15). We distinguish them in the word $U$ as:

$$U' q_{i_0-1}^\sigma a_{i_0}^{s_k} q_{i_0-1}^{-\sigma} U''.$$

The part of the sequence (15) from the beginning to the moment of deletion of the letters $q_{i_0-1}^\sigma$ and $q_{i_0-1}^{-\sigma}$ will be a sequence of the form

$$U' q_{i_0}^\sigma a_{i_0}^{s_k} q_{i_0-1}^{-\sigma} U'' \to \cdots \to Z' q_{i_0-1}^\sigma q_{i_0-1}^{-\sigma} Z''. \tag{16}$$

Suppose $\sigma = 1$. In the sequence

$$U' q_{i_0-1} a_{i_0}^{s_k} q_{i_0-1}^{-1} U'' \to \cdots \to Z' q_{i_0-1} q_{i_0-1}^{-1} Z'' \tag{16$'$}$$

the distinguished letters $q_{i_0-1}$ and $q_{i_0-1}^{-1}$ appear also in every intermediate word. Replace (16$'$) by its canonical sequence associated to the distinguished letter $q_{i_0-1}$ (see Lemma 2, §1, Chapter I),

$$U' q_{i_0-2}^{t_1} \underline{q_{i_0-1}} a_{i_0-2}^{-t_1} a_{i_0}^{s_k} \underline{q_{i_0-1}^{-1}} U'' \to \cdots \to Z' \underline{q_{i_0-1}} \; \underline{q_{i_0-1}^{-1}} Z''.$$

Since in this sequence no transformations are performed on the underlined letters $\underline{q_{i_0-1}}$ and $\underline{q_{i_0-1}^{-1}}$, the segment $a_{i_0-2}^{-t_1} a_{i_0}^{s_k}$ is independently transformed into 1 in the group $\mathfrak{A}_{q_{i_0},A,B}$. By the inductive hypothesis (part **II**), this implies that $a_{i_0-2}^{-t_1} a_{i_0}^{s_k} = 1$ in $\mathfrak{A}_P$. The letters $a_{i_0-2}$ and $a_{i_0}$ are heterotypic (Lemma 3, §2, Chapter I). By Lemma 4, §2, Chapter I we hence have $s_k = 0$, i.e. the sequence (13) is length-reducible.

Suppose $\sigma = -1$. Then the sequence (16) has the form

$$U' q_{i_0-1}^{-1} a_{i_0}^{s_k} q_{i_0-1} U'' \to \cdots \to Z' q_{i_0-1}^{-1} q_{i_0-1} Z'' \tag{16$''$}$$

Form from (16$''$) the canonical sequence associated with the letter $q_{i_0-1}$,

$$U' \underline{q_{i_0-1}^{-1}} a_{i_0}^{s_k} q_{i_0-2}^{t_1} \underline{q_{i_0-1}} a_{i_0-2}^{-t_1} U'' \to \cdots \to Z' \underline{q_{i_0-1}^{-1}} \; \underline{q_{i_0-1}} Z''. \tag{17}$$

In the sequence (17), the underlined letters $\underline{q_{i_0-1}}$ and $\underline{q_{i_0-1}^{-1}}$ are not affected by any of the transformations, and there are no insertions of the letter $q_{i_0-1}^\sigma$. This implies that the segment $a_{i_0}^{s_k} q_{i_0-2}^{t_1}$ is transformed into the identity only by transformations in the group $\mathfrak{A}_{q_{i_0-2},A,B}$. But in $\mathfrak{A}_{q_{i_0-2},A,B}$ the letter $q_{i_0-2}$ is a filtering letter. Hence $r_1 = 0$ and $a_{i_0}^{s_k} = 1$ in $\mathfrak{A}_{q_{i_0-2},A,B}$. By the inductive hypothesis, it follows that $a_{i_0}^{s_k} = 1$ in $\mathfrak{A}_P$. By Lemma 4, §2, Chapter I, this implies that $s_k = 0$, i.e. the sequence (13) is length-reducible.

Thus, we have proved that if $U = 1$ in $\mathfrak{A}_{q_{i_0-1},A,B}$, then the sequence (13) corresponding to $U$ either consists only of zeroes, or else is length-reducible. From this, it easily follows that if $U = 1$ in $\mathfrak{A}_{q_{i_0-1},A,B}$, then the sequence (13) is degenerate. Indeed, suppose $U = 1$ in $\mathfrak{A}_{q_{i_0-1},A,B}$. Then, as proved, if the sequence (13) does not consist only of zeroes, then it is length-reducible. Reducing it once, we obtain a sequence which defines a word $U_1$ of type (14), and which is equal to $U$ in the free group. This implies that also $U_1 = 1$ in $\mathfrak{A}_{q_{i_0-1},A,B}$. Hence, the length-reduced sequence either consists only of zeroes, or else is length-reducible. If it does not consist only of zeroes, then we can reduce it, obtaining a new sequence, etc.

As every length-reduction indeed reduces the length of the sequence, after a finite number of steps we, by the method above, obtain a sequence which consists only of zeroes. Hence, the sequence (13) is degenerate. The necessity is proved.

It remains to prove part **II** of the lemma for $i = i_0$, under the assumption that $i < i_0$. Suppose $E$ is an arbitrary word over the alphabet of the group $\mathfrak{A}_P$, and that $E = 1$ in $\mathfrak{A}_{q_{i_0},A,B}$. Then from $E = 1$ in $\mathfrak{A}_{q_{i_0},A,B}$ follows $E = 1$ in $\mathfrak{A}_{q_{i_0-1},A,B} \approx \overline{\mathfrak{A}}_{q_{i_0},A,B}^{(q_{i_0})}$ (see the Corollary of Lemma 3, §1, Chapter I). Hence, by the inductive hypothesis, we have that $E = 1$ in $\mathfrak{A}_P$. The converse assertion is obvious, as every defining relation of the group $\mathfrak{A}_P$ is also a relation in the group $\mathfrak{A}_{q_{i_0},A,B}$.                                                                        $\square$

**Lemma 2.** *Suppose $A = pXpX^+, B = pYpY^+$, where $X$ and $Y$ are arbitrary words over the alphabet $S$, and $A \neq B$ in the group $\mathfrak{A}_P$. Consider the word*

$$R \equiv BA^{-1}EAB^{-1}U, \tag{18}$$

*where $E$ is an arbitrary word over the alphabet of the group $\mathfrak{A}_P$ not containing the letter $p^\sigma$, and $U$ is a word of the form*

$$q_{n-1}^{r_1} a_n^{s_1} q_{n-1}^{r_2} a_n^{s_2} \cdots q_{n-1}^{r_\mu} a_n^{s_\mu} q_{n-1}^{r_{\mu+1}}. \tag{19}$$

*If the word $R$ is equal to 1 in the group $\mathfrak{A}_{q_{n-1},A,B}$, then the sequence*

$$r_1, s_1, r_2, s_2, \ldots, r_\mu, s_\mu, r_{\mu+1} \tag{13}$$

*corresponding to $U$ is degenerate.*

The proof of this lemma is obtained by almost word-by-word repetition of the part of the proof of the previous lemma in which is proved that: for the word (14) to be equal to 1 in the group $\mathfrak{A}_{q_{i_0-1},A,B}$, it is necessary that the corresponding sequence (13) is degenerate. Here, the only difference is that the index $i_0$ is replaced by $n$.

In the case that the word $U$ contains the letters $q_{n-1}^\sigma$, the presence of the segment $BA^{-1}EAB^{-1}$ does not necessitate any difference to the proof. A new part of the proof appears only in the case when $\mu = 1$ and the word $U$ does not contain the letter $q_{n-1}^\sigma$. In this case we have

$$R \equiv BA^{-1}EAB^{-1}a_n^{s_1} = 1$$

in the group $\mathfrak{A}_{q_{n-1},A,B}$. By the previous lemma, this implies $R = 1$ in $\mathfrak{A}_P$, whence by Lemma 7, §2, Chapter I we find:

$$a_n^{s_1} = 1 \quad \text{in} \quad \mathfrak{A}_P.$$

We obtained the same equality in the analogous case in the proof of the previous lemma. In all other points, the previous proof can be carried over completely, so we do not do this again.

**Corollary 1.** *Suppose the word $R$, which appears in Lemma 2, is equal to 1 in the group $\mathfrak{A}_{q_{n-1},A,B}$. Then in $\mathfrak{A}_P$ we have the equalities*

$$(E)_{\overline{S}} = 1 \quad and \quad (E)_{\overline{S^+}} = 1.$$

*Proof.* From $R = 1$ in $\mathfrak{A}_{q_{n-1},A,B}$ it follows by Lemma 2 that the sequence (13) is degenerate. From the fact that the sequence (13) degenerates in exactly the same way as in Lemma 1, it follows that $U = 1$ in the free group. Consequently

$$R \equiv BA^{-1}EAB^{-1}U = BA^{-1}EAB^{-1}$$

in the free group.

From this, by assumption, we have the equality

$$BA^{-1}EAB^{-1} = 1$$

in the group $\mathfrak{A}_{q_{n-1},A,B}$, and by Lemma 1 we find the same equality holds in the group $\mathfrak{A}_P$. As by assumption we have $A \neq B$, hence by Lemma 7, §2, Chapter I we have the equalities

$$(E)_{\overline{S}} = 1 \quad \text{and} \quad (E)_{\overline{S^+}} = 1.$$

The corollary is proved. □

**Corollary 2.** *Suppose $U$ is the word (19) and is equal to 1 in the group $\mathfrak{A}_{q_{n-1},A,B}$. Then its corresponding sequence (13) is degenerate.*

### §2. Transformations in the groups $\mathfrak{A}_{q,A,B}$

In the previous section, we defined the system of groups $\mathfrak{A}_{q,A,B}$. The alphabet of each of these groups contained the letter $q_n$, which appeared in two defining relations:

$$q_{n-1}q_n = q_n a_{n-1}, \tag{1}$$

$$a_n q_n = q_n a_n q_n^{-1} B A^{-1}. \tag{2}$$

Obviously, the letter $q_n$ is not a filtering letter of the group $\mathfrak{A}_{q,A,B}$, but bears some resemblance to one. We will say that it is a *quasi-filtering letter*.

In the current section we will introduce concepts for the quasi-filtering letter $q_n$, such as identity in individuality, canonical sequences of transformations, and other concepts, which P. S. Novikov introduced for filtering letters in [1] and [2]. We will then prove lemmas analogous to the lemmas about filtering letters given in §1, Chapter I.

From the relations (1) and (2) follow the relations:

$$\left. \begin{aligned} q_n &= q_{n-1}^{-1} q_n a_{n-1}, \\ q_n &= q_{n-1} q a_{n-1}^{-1}, \\ q_n &= a_n^{-1} q_n a_n q_n^{-1} B A^{-1}, \\ q_n &= a_n q_n A B^{-1} q_n a_n^{-1}, \end{aligned} \right\} \tag{3}$$

and, conversely, from the equalities (3) we can deduce the equalities (1) and (2). Hence, if we replace the defining relations (1) and (2) in the group $\mathfrak{A}_{q,A,B}$ by the relations (3), then we obtain the same group $\mathfrak{A}_{q,A,B}$.

Note that in sequences of transformations using to the equalities (3), it suffices to consider only the transformations:

$$\left. \begin{aligned} Xq_n Y &\to X q_{n-1}^{-1} q_n a_{n-1} Y, \\ Xq_n Y &\to X q_{n-1} q a_{n-1}^{-1} Y, \\ Xq_n Y &\to X a_n^{-1} q_n a_n q_n^{-1} B A^{-1} Y, \\ Xq_n Y &\to X a_n q_n A B^{-1} q_n a_n^{-1} Y, \end{aligned} \right\} \tag{3'}$$

as any reverse transformations can easily be eliminated. Indeed, the transformation

$$X q_{n-1}^{-1} q_n a_{n-1} Y \to X q_n Y$$

can be replaced by the sequence of transformations

$$X q_{n-1}^{-1} q_n a_{n-1} Y \to X q_{n-1}^{-1} q_{n-1} q_n a_{n-1}^{-1} a_{n-1} Y \to \cdots \to X q_n Y$$

where we first carry out the transformation

$$X q_n Y \to X q_{n-1} q_n a_{n-1}^{-1} Y$$

and after this we carry out two deletions.

Analogously, we can also remove the reverse transformations in the remaining cases.

*In the sequel, we will assume that such reverse transformations, which we can eliminate by the above schema, are never performed.*

Consider an arbitrary sequence of elementary transformations, satisfying the above requirement, in the group $\mathfrak{A}_{q,A,B}$:

$$E_1 \to E_2 \to \cdots \to E_m. \tag{4}$$

The *index* of the sequence (4) of transformations will be the total number $\lambda$ of transformations in the sequence (4) which are either insertions of the letters $q_n^{-1}$ and $q_n$, or transformations from the schema (3′).

Analogously to how it was done in §1, Chapter I for filtering letters, we will now define the notion of being identical in individuality for quasi-filtering letters $q_n$ or $q_n^{-1}$ appearing in distinct words from the sequence (4) of elementary transformations.

Every transformation $E_j \to E_{j+1}$ in the sequence (4) has the form

$$XRY \to XSY. \tag{5}$$

Any two occurrences of a letter $q_n$ or $q_n^{-1}$ which appear in exactly the same place inside both the specified segments $X$ (or $Y$) inside both of the words $XRY$ and $XSY$, will be said to have the same individuality. Moreover, if (5) is a transformation of type (3′), then we will say that the letter $q_n$, which in this case corresponds to $R$, will have the same individuality as the letter $q_n$ inside the word $S$ on the right. We extend this definition by transitivity, and say that two letters of type $q_n^\sigma$ in some words from the sequence (4) are different in individuality, if the identity of their individualities does not follow from this definition. It is clear that the only letter which can be identical in individuality to $q_n$ $(q_n^{-1})$ is $q_n$ $(q_n^{-1})$.

Unlike for filtering letters, new individualities of quasi-filtering letters $q_n^\sigma$ can appear in the sequence (4) not only as a result of insertions, but also from some of the transformations from (3′).

Suppose that in the sequence (4) we have a transformation according to the schema

$$X\underline{q_n}Y \to Xa_n^{-1}\underline{q_n}a_nq_n^{-1}BA^{-1}Y.$$

The underlined two letters $\underline{q_n}$ have the same individuality. The individuality $q_n^{-1}$, which appears on the right side of this transformation, will be said to have been *directly engendered by* the individuality $\underline{q_n}$, and the letter $\underline{q_n}$ will be said to be the direct engendering letter for the individuality of $q_n^{-1}$.

If we have a transformation

$$X\underline{q_n}Y = Xa_n\underline{q_n}AB^{-1}q_na_n^{-1}Y \tag{5′}$$

then we will say that non-underlined letter $q_n$ was *directly engendered by* that individuality of the letter $q_n$ which in (5′) is underlined. The individuality underlined in (5′) will be said to be the direct engendering letter for the one which was not underlined. The letter $q_n^{-1}$ can itself be engendered, but cannot be an engendering letter. Note that an engendered letter $q_n^\sigma$ always lies to the right of its engendering letter $q_n$. The individuality $q_n$ engendered in the sequence of transformations (4) can itself engender new individualities of letters $q_n^\sigma$, which when $\sigma = 1$ can engender yet more, etc. Any sequence of individualities of letters $q_n^\sigma$

$$q_n^{(1)}, q_n^{(2)}, q_n^{(3)}, \ldots, q_n^{(k)},$$

appearing in the sequence of transformations (4), such that every letter $q_n^{(i+1)}$ is directly engendered by the letter $q_n^{(i)}$, will be called a *hereditary chain of letters $q_n$ engendering the individuality $q_n^{(k)}$ in* (4).

Every individuality $q_n^{(i)}$ $(i > 1)$ will be called *hereditarily engendered* (or simply engendered) from the individuality $q_n^{(1)}$.

**Remark 1.** If the sequence of elementary transformations

$$E_1 \to E_2 \to \cdots \to E_m \tag{4}$$

in the group $\mathfrak{A}_{q,A,B}$ does not contain insertions of the letters $q_n^\sigma$, then every individuality $q_n^\sigma$ occurring in the sequence (4) and which is not contained in $E_1$, is necessarily engendered in (4) by one of the letters $q_n$ appearing in the word $E_1$.

Analogously to how we previously defined the correspondence $L_p$ for any filtering letter $p$, we will construct a correspondence $L_{q_n}$ for any quasi-filtering letter $q_n$. In this case, however, the words on the right-hand side of the correspondence can contain the letters $q_n^\sigma$. The correspondence $L_{q_n}$ will have the form

$$U \equiv q_{n-1}^{r_1} a_n^{s_1} q_{n-1}^{r_2} a_n^{s_2} \cdots q_{n-1}^{r_\mu} a_n^{s_\mu} q_{n-1}^{r_{\mu+1}} \underline{\hspace{3cm}}$$
$$\underline{\hspace{3cm}} a_{n-1}^{r_1} (a_n q_n^{-1} B A^{-1})^{s_1} a_{n-1}^{r_2} (a_n q_n^{-1} B A^{-1})^{s_2} \cdots a_{n-1}^{r_\mu} (a_n q_n^{-1} B A^{-1})^{s_\mu} a_{n-1}^{r_{\mu+1}}.$$

The words which can appear on the left-hand side of this correspondence will henceforth be denoted by $U$ with varying indices; its corresponding word will be denoted by $W$ with the same index, and $W^{-1}$ will be denoted by $V$, also with the same index.

When corresponding to the word

$$U \equiv q_{n-1}^{r_1} a_n^{s_1} q_{n-1}^{r_2} a_n^{s_2} \cdots q_{n-1}^{r_\mu} a_n^{s_\mu} q_{n-1}^{r_{\mu+1}},$$

the word $V$ will have the form

$$a_{n-1}^{r_1} (a_n q_n^{-1} B A^{-1})^{s_1} a_{n-1}^{r_2} (a_n q_n^{-1} B A^{-1})^{s_2} \cdots a_{n-1}^{r_\mu} (a_n q_n^{-1} B A^{-1})^{s_\mu} a_{n-1}^{r_{\mu+1}}.$$

Obviously, if $U = U_1 U_2$, we have $V = V_2 V_1$.

It is easily seen, that if the words $U$ and $V^{-1}$ correspond to each other under $L_{q_n}$, then by (3') we have the equality

$$q_n = U q_n V.$$

**Lemma 1.** *Suppose the sequence of elementary transformations in the group $\mathfrak{A}_{q,A,B}$*

$$E_1' \underline{q_n} E_1'' \to E_2 \to \cdots \to E_m \tag{6}$$

*does not contain insertions of the letters $q_n^\sigma$. Then the word $E_1' q_n E_1''$ can be transformed into $E_m$ by a sequence of transformations*

$$E_1' \underline{q_n} E_1'' \to \cdots \to E_1' U \underline{q_n} V E_1'' \to \cdots \to E_m, \tag{7}$$

*satisfying the following properties:*

1. *$U$ and $V^{-1}$ correspond to each other under $L_{q_n}$.*
2. *The part*

   $$E_1' \underline{q_n} E_1'' \to \cdots \to E_1' U \underline{q_n} V E_1''$$

   *of the sequence (7) consists only of transformations of type (3') which affect the underlined individuality $\underline{q_n}$. In the remaining part of the sequence (7), the underlined letter $\underline{q_n}$ is not affected by any transformations of type (3').*
3. *In the sequence (7) there are also no insertions of the letters $q_n^\sigma$, and its index coincides with the index of the sequence (6).*

*Proof.* The proof will be carried out by induction on the number $\eta$ of applications of transformations of type $(3')$ to the distinguished individuality $\underline{q}_n$ in the sequence (6). If the number $\eta$ is zero, then the lemma is obvious. Suppose that for all $\eta < \eta_0$ the lemma has been proved, and that the number of transformations considered in (6) is $\eta_0 > 0$. Suppose $E_{i_1} \to E_{i_1+1}$ is the first of these transformations. For definiteness, we assume that $E_{i_1} \to E_{i_1+1}$ is a transformation by the schema $(5')$. This assumption may be made without any loss of generality. Then the sequence (6) has the form

$$E_1'\underline{q}_n E_1'' \to \cdots \to E_{i_1}'\underline{q}_n E_{i_1}'' \to E_{i_1}'a_n\underline{q}_n AB^{-1}q_n a_n^{-1}E_{i_1}'' \to \cdots \to E_m,$$

and moreover in this sequence, up to the $i_1$-th transformation the underlined letter $\underline{q}_n$ is not affected, and all transformations are carried out on $E_1'$ and $E_1''$. Rearrange the transformations in (6) by placing the initial $i_1$ transformations first. We obtain the sequence

$$E_1'a_n\underline{q}_n AB^{-1}q_n a_n^{-1}E_1'' \to \cdots \to E_{i_1}'a_n\underline{q}_n AB^{-1}q_n a_n^{-1}E_{i_1}'' \to \cdots \to E_m, \qquad (6')$$

in which $\eta_0 - 1$ transformations of type $(3')$ are already performed on the letter $\underline{q}_n$. If the index of the sequence (6) is $\lambda$, then the index of the sequence $(6')$ will be $\lambda - 1$. The sequence $(6')$ does not contain insertions of the letters $q_n^\sigma$. By the inductive hypothesis, for the sequence $(6')$ there will be a sequence

$$E_1'a_n\underline{q}_n AB^{-1}q_n a_n^{-1}E_1'' \to \cdots \to E_1'a_n U_1\underline{q}_n V_1 AB^{-1}q_n a_n^{-1}E_1'' \to \cdots \to E_m, \quad (7')$$

which has index $\lambda - 1$ and which satisfies the requirements of the lemma. In particular, the words $U_1$ and $V_1^{-1}$ correspond under $L_{q_n}$.

The sequence of transformations

$$E_1'\underline{q}_n E_1'' \to E_1'a_n\underline{q}_n AB^{-1}q_n a_n^{-1}E_1'' \to \cdots$$
$$\cdots \to E_1'a_n U_1\underline{q}_n V_1 AB^{-1}q_n a_n^{-1}E_1'' \to \cdots \to E_m, \qquad (7'')$$

which, starting from the second transformation, coincides completely with the sequence $(7')$, will be the desired sequence. First, the words $a_n U_1$ and $(V_1 AB^{-1}q_n a_n^{-1})^{-1}$ correspond under $L_{q_n}$. The index of the sequence $(7'')$ will be $\lambda$, and there are no insertions of the letters $q_n^\sigma$. Property 2 of the lemma also follows, as it holds for the sequence $(7')$. Lemma 1 is proved. $\qquad\square$

A sequence of transformations satisfying all properties of Lemma 1 will be called a *canonical sequence for (6), with respect to the letter $\underline{q}_n$.*

Lemma 1 implies that for every sequence of elementary transformations (6) without insertions of the letters $q_n^\sigma$, there exists a canonical sequence with respect to any letter $q_n$ which is contained in the initial word of the sequence (6). In order not to mix up the distinct individualities of the letters $q_n^\sigma$ in the sequence of transformations (4), we will sometimes assign them indices, as superscripts in parentheses.

Suppose we are given a sequence of elementary transformations (4) and a hereditary chain

$$q_n^{(1)}, q_n^{(2)}, q_n^{(3)}, \ldots, q_n^{(k)}, \qquad (8)$$

of letters $q_n^\sigma$ from this sequence, where $q_n^{(1)}$ is an individuality which is contained in the first word of the sequence (4). We will call a *canonical sequence for (4), with respect to the hereditary chain (8),* any sequence of elementary transformations

$$E_1 \to \cdots \to E_m,$$

in which all transformations of type $(3')$ in the sequence (4) are performed first on the individuality $q_n^{(1)}$, then all transformations of type $(3')$ on the individuality $q_n^{(2)}$, etc, and finally all transformations of type $(3')$ on the individuality $q_n^{(k)}$, if in the sequence (4) there are such, and only after this are the remaining transformations in the sequence (4) performed.

**Lemma 2.** *For any sequence of elementary transformations* (4) *without insertions of the letters* $q_n^\sigma$, *and any given individuality* $q_n^{\sigma_1}$ *which is engendered in the sequence* (4), *there exists a canonical sequence with respect to the hereditary chain of letters* $q_n^\sigma$ *which engenders the individuality* $q_n^{\sigma_1}$ *in* (4).

*Proof.* Suppose $\sigma_1 = 1$. As the given individuality $q_n$, by assumption, is engendered in the sequence (4), it follows by Remark 1 that it is engendered by some letter $q_n$ which appears in the word $E_1$. Let (8) be the hereditary chain which engenders the given individuality $q_n$ in (4), in which case $q_n^{(1)}$ appears in the word $E_1$. Then the word $E_1$ has the form $E_1' q_n^{(1)} E_1''$. By Lemma 1, for the sequence (4) there exists a canonical sequence of transformations with respect to the letter $q_n^{(1)}$. We write one down:

$$T_1) \qquad E_1' q_n^{(1)} E_1'' \to \cdots \to E_1' U_1 q_n^{(1)} V_1 E_1'' \to \cdots \to E_m. \qquad (9)$$

All the individualities of the letter $q_n$, directly engendered in (4) by the letter $q_n^{(1)}$, appear in the word $V_1$. In particular, the word $V_1$ contains the letter $q_n^{(2)}$. Then $V_1$ has the form

$$V_1' A B^{-1} q_n^{(2)} a_n^{-1} V_1'',$$

where $V_1'$ and $V_1''$ are some words of type $V$.

In the second part of the sequence (9),

$$E_1' U_1 q_n^{(1)} V_1 E_1'' \equiv E_1' U_1 q_n^{(1)} V_1' A B^{-1} q_n^{(2)} a_n^{-1} V_1'' E_1'' \to \cdots \to E_m, \qquad (9')$$

no transformations, except for deletions, are performed on the individuality $q_n^{(1)}$. By Lemma 1, there exists for sequence $(9')$ a canonical sequence of transformations with respect to the letter $q_n^{(2)}$, which we write as

$$T_2) \quad E_1' U_1 q_n^{(1)} V_1' A B^{-1} q_n^{(2)} a_n^{-1} V_1'' E_1'' \to \cdots$$
$$\cdots \to E_1' U_1 q_n^{(1)} V_1' A B^{-1} U_2 q_n^{(2)} V_2 a_n^{-1} V_1'' E_1'' \to \cdots \to E_m. \qquad (10)$$

Every individuality of the letter $q_n$, directly engendered in (4) by the letter $q_n^{(2)}$, appears in the word $V_2$, so $V_2$ has the form

$$V_2' A B^{-1} q_n^{(3)} a_n^{-1} V_2''.$$

In the second part of the sequence (10),

$$E_1' U_1 q_n^{(1)} V_1' A B^{-1} U_2 q_n^{(2)} V_2' A B^{-1} q_n^{(3)} a_n^{-1} V_2'' a_n^{-1} V_1'' E_1'' \to \cdots \to E_m \qquad (10')$$

no transformations, except for deletions, are performed on the individualities $q_n^{(1)}$ and $q_n^{(2)}$. By considering the canonical sequence with respect to $q_n^{(3)}$, we continue to $q_n^{(4)}$, etc. until we reach $q_n^{(k)}$. The canonical sequence with respect to $q_n^{(k)}$ has the form

$$T_k) \quad E_1' U_1 q_n^{(1)} V_1' A B^{-1} U_2 q_n^{(2)} V_2' \cdots q_n^{(k-1)} V_{k-1}' A B^{-1} q_n^{(k)} a_n^{-1} V_{k-1}'' a_n^{-1} \cdots$$
$$\cdots a_n^{-1} V_2'' a_n^{-1} V_1'' E_1'' \to \cdots \to E_1' U_1 q_n^{(1)} V_1' A B^{-1} U_2 q_n^{(2)} V_2' \cdots$$
$$\cdots q_n^{(k-1)} V_{k-1}' A B^{-1} U_k q_n^{(k)} V_k a_n^{-1} V_{k-1}'' a_n^{-1} V_{k-2}'' a_n^{-1} \cdots a_n^{-1} V_2'' a_n^{-1} V_1'' E_1'' \to \cdots$$
$$\cdots \to E_m.$$

In the second part of this sequence, no transformations, except for deletions, are performed on the individualities $q_n^{(1)}, q_n^{(2)}, \ldots, q_n^{(k)}$. Using the first part of the sequences $T_1, T_2$, etc. until $T_{k-1}$, and the full sequence of transformations $T_k$, we obtain

$$E_1 \equiv E_1' q_n^{(1)} E_1'' \to \cdots \to E_1' U_1 q_n^{(1)} V_1' A B^{-1} q_n^{(2)} a_n^{-1} V_1'' E_1'' \to \cdots$$
$$\cdots \to E_1' U_1 q_n^{(1)} V_1' A B^{-1} U_2 q_n^{(2)} V_2' q_n^{(3)} a_n^{-1} V_2'' a_n^{-1} V_1'' E_1'' \to \cdots$$
$$\to E_1' U_1 q_n^{(1)} V_1' A B^{-1} U_2 q_n^{(2)} V_2' A B^{-1} U_3 q_n^{(3)} \cdots q_n^{(k-1)} V_{k-1}' A B^{-1} q_n^{(k)} a_n^{-1} V_{k-1}'' \cdots$$
$$\cdots a_n^{-1} V_2'' a_n^{-1} V_1'' E_1'' \to \cdots \to E_1' U_1 q_n^{(1)} V_1' A B^{-1} U_2 q_n^{(2)} V_2' A B^{-1} U_3 q_n^{(3)} \cdots$$
$$\cdots q_n^{(k-1)} V_{k-1}' A B^{-1} U_k q_n^{(k)} V_k a_n^{-1} V_{k-1}'' \cdots a_n^{-1} V_2'' a_n^{-1} V_1'' E_1'' \to \cdots \to E_m. \tag{11}$$

This is the desired canonical sequence of transformations.

If $\sigma_1 = -1$, then $q_n^{(k)}$ is some individuality of the letter $q_n^{-1}$, on which no transformations, except for deletions, are performed in the sequence (4).

In this case, in the analogous way we find a canonical sequence of transformations, which now has the form

$$E_1 \equiv E_1' q_n^{(1)} E_1'' \to \cdots \to E_1' U_1 q_n^{(1)} V_1' A B^{-1} q_n^{(2)} a_n^{-1} V_1'' E_1'' \to \cdots \to$$
$$\to E_1' U_1 q_n^{(1)} V_1' A B^{-1} U_2 q_n^{(2)} V_2' A B^{-1} q_n^{(3)} a_n^{-1} V_2'' a_n^{-1} V_1'' E_1'' \to \cdots \to$$
$$\to E_1' U_1 q_n^{(1)} V_1' A B^{-1} U_2 q_n^{(2)} V_2' A B^{-1} U_3 q_n^{(3)} \cdots U_{k-1} q_n^{(k-1)} V_{k-1}' a_n q_n^{-1} \to$$
$$\to B A^{-1} V_{k-1}'' a_n^{-1} \cdots a_n^{-1} V_2'' a_n^{-1} V_2'' a_n^{-1} V_1'' E_1'' \to \cdots \to E_m. \tag{11'}$$

This completes the proof of Lemma 2. $\qquad\square$

**Remark 2.** The original and canonical sequences of transformations always have the same index, as they differ from one another only by a rearrangement of the transformations.

In what follows, when writing the canonical sequence of transformations, with respect to one or another hereditary chain, we will omit the superscripted indices for the letters $q_n^\sigma$, if this does not cause confusion.

**Lemma 3.** *Suppose the sequence of elementary transformations in the group $\mathfrak{A}_{q,A,B}$*

$$XY \to \cdots \to 1, \tag{12}$$

*where $X$ and $Y$ are arbitrary words of the group $\mathfrak{A}_{q,A,B}$, has index $\lambda$ and does not contain insertions of the letters $q_n^\sigma$. Then we can specify a sequence of elementary transformations in the group $\mathfrak{A}_{q,A,B}$:*

$$YX \to \cdots \to 1, \tag{12'}$$

*which also has index $\lambda$, and does not contain insertions of the letters $q_n^\sigma$.*

*Proof.* Suppose $\lambda = 0$. Then in the sequence (12), only deletions are performed on the letters $q_n^\sigma$. The proof on this case proceeds by induction on the number $k$ of letters $q_n^\sigma$ appearing in the word $XY$.

Suppose $k = 0$. Then the sequence (12') is obtained in the following way:

$$YX \to \cdots \to X^{-1} XYX \to \cdots \to X^{-1} X \to \cdots \to 1.$$

Here we first insert the word $X^{-1} X$, and then transform the segment $XY$ into the empty word by using the sequence (12), and finally transform the word $X^{-1} X$ into the identity via deletions. As, by assumption, the word $XY$ does not contain any letters $q_n^\sigma$, we did not perform any insertions of letters $q_n^\sigma$ when doing the above.

Suppose that for all $k < k_1$ and $\lambda = 0$ the lemma is proved. We will prove the lemma for $k = k_1$ and $\lambda = 0$.

Suppose

$$X = X_1 q_n^{\sigma_1} X_2 q_n^{\sigma_2} X_3 \cdots X_\mu q_n^{\sigma_\mu} X_{\mu+1},$$

$$Y = Y_1 q_n^{\delta_1} Y_2 q_n^{\delta_2} Y_3 \cdots Y_\nu q_n^{\delta_\nu} Y_{\nu+1},$$

where $\sigma_i = \pm 1, \delta_i = \pm 1$, and the words $X_i$ and $Y_i$ do not contain the letters $q_n^\sigma$. All letters $q_n^\sigma$ appearing in the word $XY$ must be deleted in the sequence (12). Consider those which are deleted first.

There are three possible cases:

1) both letters $q_n^\sigma$ which are deleted first in (12) come from the word $X$;
2) both come from the word $Y$;
3) one comes from the word $X$, and the other from $Y$.

Consider **Case 1)**. Suppose the first deleted letters are $q_n^{\sigma_i}$ and $q_n^{\sigma_{i+1}}$. Then evidently $\sigma_{i+1} = -\sigma_i$ and the sequence (12) has the form

$$X_1' q_n^{\sigma_i} X_{i+1} q_n^{-\sigma_i} X_2' Y \to \cdots \to Z_1 q_n^{\sigma_i} q_n^{-\sigma_i} Z_2 \to Z_1 Z_2 \to \cdots \to 1.$$

As, by assumption, $\lambda = 0$, there are no transformations affecting the letters $q_n^{\sigma_i}$ and $q_n^{-\sigma_i}$ except deletions. Therefore, we have the following independent transformations:

1. $X_1' \to \cdots \to Z_1,$          3. $X_2' Y \to \cdots \to Z_2,$

2. $X_{i+1} \to \cdots \to 1,$          4. $Z_1 Z_2 \to \cdots \to 1,$

each of which has index 0.

By combining the sequences $1, 3$, and $4$, we find the following sequence with index 0:

$$X_1' X_2' Y \to \cdots \to Z_1 X_2' Y \to \cdots \to Z_1 Z_2 \to \cdots \to 1.$$

As the word $X_1' X_2 Y$ contains $k_1 - 2$ letters $q_n^\sigma$, by the inductive hypothesis it follows that we have a sequence of transformations

$$Y X_1' X_2' \to \cdots \to 1 \tag{12''}$$

with index 0.

By using the sequences (2) and (12''), we obtain the following sequence with index 0:

$$Y X \equiv Y X_1' q_n^{\sigma_i} X_{i+1} q_n^{-\sigma_i} X_2' \to \cdots \to Y X_1' q_n^{\sigma_i} q_n^{-\sigma_i} X_2' \to Y X_1' X_2' \to \cdots \to 1.$$

**Case 2)** is proved analogously.

**Case 3).** As we consider the first deletion to occur in the sequence (12), in this case the deletion is between the letters $q_n^{\sigma_\mu}$ in $X$ and $q_n^{\delta_1}$ in $Y$. It follows that $\sigma_\mu = -\delta_1$. The sequence (12) has the form

$$X' q_n^{-\delta_1} X_{\mu+1} Y_1 q_n^{\delta_1} Y' \to \cdots \to Z_1 q_n^{-\delta_1} q_n^{\delta_1} Z_2 \to Z_1 Z_2 \to \cdots \to 1, \tag{13}$$

where no transformations are performed on the letters $q_n^{-\delta_1}$ and $q_n^{\delta_1}$ except for the deletion. In the sequence (13) we hence have the following independent transformations:

1. $X_1' \to \cdots \to Z_1,$          3. $X_{\mu+1} Y_1 \to \cdots \to 1,$

2. $Y' \to \cdots \to Z_2,$          4. $Z_1 Z_2 \to \cdots \to 1,$

each of which has index 0.

As the sequence (13) does not contain insertions of the letters $q_n^\sigma$ or transformations from the schema $(3')$, the number of letter $q_n^\sigma$ in the word $Z_1 Z_2$ is not greater than

$k_1 - 2$. By the inductive hypothesis, the existence of the sequence 4 implies that there exists a sequence $4'$: $Z_2 Z_1 \to \cdots \to 1$ with index 0. By combining the sequences $1, 2, 4'$, and 3, we obtain a sequence of transformations with index $\lambda = 0$:

$$YX \equiv Y_1 q_n^{\delta_1} Y' X' q_n^{-\delta_1} X_{\mu+1} \to \cdots \to Y_1 q_n^{\delta_1} Y' Z_1 q_n^{-\delta_1} X_{\mu+1} \to \cdots$$

$$\cdots \to Y_1 q_n^{\delta_1} Z_2 Z_1 q_n^{-\delta_1} X_{\mu+1} \to \cdots \to Y_1 q_n^{\delta_1} q_n^{-\delta_1} X_{\mu+1} \to Y_1 X_{\mu+1} \to \cdots \to 1.$$

In this way, we have proved the lemma for $\lambda = 0$. Assume now that the lemma has been proved for all $0 < \lambda < \lambda_1$. We will prove it for $\lambda = \lambda_1$. As the sequence (12) does not contain insertions of letters $q_n^\sigma$, and as the index $\lambda > 0$, there is some transformation from the schema $(3')$ applied to some letter $q_n$ appearing in the word $XY$. We denote one of these letters by $\underline{q}_n$. This appears either in $X$, or in $Y$. The two cases are analogous. We will hence only consider the case when this letter $\underline{q}_n$ appears in the word $X$. Then the word $XY$ has the form

$$X_1' \underline{q}_n X_2' Y.$$

We construct a canonical sequence for (12) with respect to the selected letter $\underline{q}_n$:

$$X_1' \underline{q}_n X_2' Y \to \cdots \to X_1' U \underline{q}_n V X_2' Y \to \cdots \to 1. \tag{14}$$

By the way which we chose the letter $\underline{q}_n$, the part

$$X_1' \underline{q}_n X_2' Y \to \cdots \to X_1' U \underline{q}_n V X_2' Y \tag{14'}$$

of the sequence (14) has index $\lambda_2 > 0$. Hence, the second part

$$X_1' U \underline{q}_n V X_2' Y \to \cdots \to 1$$

of the sequence (14) has index $\lambda_3 = \lambda_1 - \lambda_2 < \lambda_1$. By the inductive hypothesis, we can find a sequence of elementary transformations with index $\lambda_3$:

$$Y X_1' U \underline{q}_n V X_2' \to \cdots \to 1. \tag{14''}$$

By combining the sequences $(14')$ and $(14'')$, we find the desired sequence

$$YX \equiv Y X_1' \underline{q}_n X_2' \to \cdots \to Y X_1' U \underline{q}_n V X_2' \to \cdots \to 1$$

with index $\lambda = \lambda_3 + \lambda_2 = \lambda_1$. Lemma 3 is proved. $\qquad\square$

**Main Lemma.** *Suppose the words $A$ and $B$ are as in Lemma 1, §1, Chapter II, and $A \neq B$ in the group $\mathfrak{A}_P$. Then for any $\lambda$ the following properties about transformations in the group $\mathfrak{A}_{q,A,B}$ hold:*

  I. *If in the sequence of elementary transformations with index $\lambda$*

$$E_1 \to E_2 \to \cdots \to E_m \tag{15}$$

  *the word $E_m$ does not contain the letter $q_n^{-1}$, then we can specify a sequence of elementary transformations with index $\lambda_1 \leq \lambda$*

$$E_1 \to E_2' \to \cdots \to E_m,$$

  *which does not contain insertions of the letter $q_n^\sigma$.*

 II. *Suppose $U$ and $V^{-1}$ are two words corresponding under $L_{q_n}$, and that there is a sequence of elementary transformations with index $\lambda$*

$$V \to \cdots \to 1, \tag{16}$$

  *which does not contain insertions of the letter $q_n$. Then $U = 1$ in the free group.*

III. *There does not exist any sequence of elementary transformations with index $\lambda$ which transforms into the empty word some word of the form*

$$Q_1 \equiv V_1 a_n^{-1} V_2 a_n^{-1} \cdots a_n^{-1} V_{k-1} a_n^{-1} V_k,$$

*where $k > 1$, and $V_i$ is an arbitrary word from the right-hand side of the correspondence $L_{q_n}$.*

IV. *There does not exist any sequence of elementary transformations with index $\lambda$ which transforms into the empty word some word of the form*

$$Q_2 \equiv U_1 q_n G_1 AB^{-1} U_2 q_n G_2 AB^{-1} \cdots U_k q_n G_k AB^{-1} U_{k+1},$$

*where $k \geq 1$; $G_i$ is a word of the form $V_{i,1} a_n V_{i,2} a_n \cdots V_{i,t_i-1} a_n V_{i,t_i}$; $U_j$ is an arbitrary word from the left-hand side of the correspondence $L_{q_n}$; and $V_{i,j}$ is an arbitrary word from the right side of the correspondence $L_{q_n}$.*

*Proof.* The proof of the lemma will be by induction on the index $\lambda$ which appears in all four parts of the lemma.

Suppose $\lambda = 0$. Then there is nothing to show for **Part I**, as the sequence (15) does not contain insertions of the letter $q_n^\sigma$.

We proceed to **Part II**. We have a sequence of elementary transformations (16) with index $\lambda = 0$. Then in (16) there are no insertions of letters $q_n^\sigma$, and we do not apply the schema (3′), i.e. there are no transformations, except deletions, performed on the letters $q_n^\sigma$ appearing in the word $V$. The word $V$ has the form

$$a_{n-1}^{-r_{\mu+1}} (a_n q_n^{-1} BA^{-1})^{-s_\mu} a_{n-1}^{-r_\mu} \cdots a_{n-1}^{-r_2} (a_n q_n^{-1} BA^{-1})^{-s_1} a_{n-1}^{-r_1},$$

where all $s_i$ can be assumed to be non-zero. The proof of **II** for $\lambda = 0$ will be by induction on the number $m$ of letters $q_n^\sigma$ contained in the word $V$. If $m = 0$, then $V = a_{n-1}^{-r_1}$. In this case, (16) is a sequence of transformations in the group $\mathfrak{A}_{q_{n-1},A,B}$. By Lemma 1, §1, Chapter II it follows that $a_{n-1}^{-r_1} = 1$ in $\mathfrak{A}_P$. This implies that $r_1 = 0$, whence we deduce $U = q_{n-1}^{r_1} = 1$ in the free group.

Suppose now that for all $m$, less than some $m_0 > 0$, property **II** has been proved when $\lambda = 0$. We will now prove it when $\lambda = 0$ and $m = m_0$. The word $V$ contains letters $q_n^\sigma$, which in the sequence (16) must be deleted. Consider the first of these to be deleted. Between them in the word $V$ there can be no other letters $q_n^\sigma$. Hence the word $V$ has the form

$$V \equiv V'(a_n q_n^{-1} BA^{-1})^\delta a_{n-1}^{-r_{j_0}} (a_n q_n^{-1} BA^{-1})^{-\delta} V'',$$

where $\delta = \pm 1$, and $V', V''$ are words of type $V$. The cases $\delta = 1$ and $\delta = -1$ are analogous. Suppose $\delta = 1$. Then the sequence (16) has the form

$$V \equiv V' a_n \underline{q}_n^{-1} BA^{-1} a_{n-1}^{-r_{j_0}} AB^{-1} \underline{q}_n a_n^{-1} V'' \to \cdots$$
$$\cdots \to Z_1 \underline{q}_n^{-1} \underline{q}_n Z_2 \to Z_1 Z_2 \to \cdots \to 1.$$

As no transformations are here performed on the letters $\underline{q}_n$ and $\underline{q}_n^{-1}$, except for deletions, we have the following independent transformations with index 0:

1. $V' a_n \to \cdots \to Z_1$,     3. $BA^{-1} a_{n-1}^{-r_{j_0}} AB^{-1} \to \cdots \to 1$,

2. $a_n^{-1} V'' \to \cdots \to Z_2$,     4. $Z_1 Z_2 \to \cdots \to 1$,

The word $BA^{-1} a_{n-1}^{-r_{j_0}} AB^{-1}$ does not contain any letters $q_n^\sigma$. Hence the sequence 3 is a sequence of transformations in the group $\mathfrak{A}_{q_{n-1},A,B}$. By Lemma 1, §1, Chapter II, it

follows that $BA^{-1}a_{n-1}^{-r_{j_0}}AB^{-1} = 1$ in $\mathfrak{A}_P$, whence by Lemma 6, §2, Chapter I we have $r_{j_0} = 0$. Hence the word $V$ has the form

$$V'a_n q_n^{-1} BA^{-1}AB^{-1} q_n a_n^{-1} V''.$$

The word $U$, corresponding under $L_{q_n}$ to the word $V^{-1}$, will have the form $U''a_n a_n^{-1} U'$, where $U''$ and $U'$ correspond, respectively, to $(V'')^{-1}$ and $(V')^{-1}$.

By combining the sequences 1, 2, and 4 , we obtain a sequence of elementary transformations with index 0 as:

$$V'V'' \to V'a_n a_n^{-1} V'' \to \cdots \to Z_1 a_n^{-1} V'' \to \cdots \to Z_1 Z_2 \to \cdots \to 1.$$

As the word $V'V''$ contains $m_0 - 2$ occurrences of letters $q_n^\sigma$, it follows by the inductive hypothesis that we have $U''U' = 1$ in the free group. Hence, also the word $U = U''a_n a_n^{-1} U'$ is equal to 1 in the free group.

**Part III** is proved by contradiction. Suppose we have some sequence of elementary transformations with index 0:

$$Q_1 \equiv V_1 a_n^{-1} V_2 a_n^{-1} \cdots V_{k-1} a_n^{-1} V_k \to \cdots \to 1, \tag{17}$$

where $k > 1$. If the word $Q_1$ does not contain letters $q_n^\sigma$, then it has the form

$$a_{n-1}^{r_1} a_n^{-1} a_{n-1}^{r_2} a_n^{-1} \cdots a_{n-1}^{r_{k-1}} a_n^{-1} a_{n-1}^{r_k}.$$

Then (17) is a sequence of transformations in the group $\mathfrak{A}_{q_{n-1},A,B}$. By Lemma 1, §1, Chapter II, it follows that

$$a_{n-1}^{r_1} a_n^{-1} a_{n-1}^{r_2} a_n^{-1} \cdots a_{n-1}^{r_{k-1}} a_n^{-1} a_{n-1}^{r_k} = 1 \quad \text{in } \mathfrak{A}_P$$

and by Lemma 1, §2, Chapter I we find

$$(a_{n-1}^{r_1} a_n^{-1} a_{n-1}^{r_2} a_n^{-1} \cdots a_{n-1}^{r_{k-1}} a_n^{-1} a_{n-1}^{r_k})_{\overline{S}} = a_n^{-k+1} = 1 \quad \text{in } \mathfrak{A}_P,$$

which is impossible, as by assumption $k > 1$.

Suppose now that **III** is proved when $\lambda = 0$ for all words $Q_1$ in which the number of letters $q_n^\sigma$ is less than some number $m_0 > 0$. We will prove that it remains true when $\lambda = 0$ also for words $Q_1$ with $m_0$ letters $q_n^\sigma$. Suppose the sequence (17) transforms such a word into the empty word. Every letter $q_n^\sigma$ appearing in the word $Q_1$ must eventually be deleted in the sequence (17). As $\lambda = 0$, no new letters $q_n^\sigma$ can appear in the sequence (17). The first deletion of letters $q_n^\sigma$ in the sequence (17) is either right-sided, or else left-sided. As these cases are entirely analogous, we will consider only the case of a right-sided deletion. We then have two cases.

1. The first deleted letters $q_n^\sigma$ in (17) appear in the same segment $V_i$ of the word $Q_1$.

2. The two letters are contained in different segments $V_{i_1}$ and $V_{i_2}$ of the word $Q_1$.

**Case 1).** The word $Q_1$ will have the form

$$Q_1' a_n q_n^{-1} BA^{-1} a_{n-1}^{r_j} AB^{-1} q_n a_n^{-1} Q_1'',$$

where the explicitly written letters $q_n$ and $q_n^{-1}$ are the first to be deleted in the sequence (17). Then, analogously to how we did in **II** for the word $V$, we find a sequence of elementary transformations with index 0,

$$Q_1' Q_1'' \to Q_1' a_n a_n^{-1} Q_1'' \to \cdots \to 1. \tag{17'}$$

But $Q_1'Q_1''$ is a word of type $Q_1$, containing $m_0 - 2$ occurrences of the letters $q_n^\sigma$. By the inductive hypothesis, the sequence (17') cannot exist. This is a contradiction.

**Case 2).** The letters $\underline{q}_n$ and $\underline{q}_n^{-1}$, the first to be deleted, appear in different segments $V_{i_1}$ and $V_{i_1 \pm t}$ of the word $Q_1$. Obviously, between them there can be no letters $q_n^\sigma$. Hence $Q_1$ has the form

$$Q_1' a_n \underline{q}_n^{-1} BA^{-1} a_{n-1}^{r_{i_1}} a_n^{-1} a_{n-1}^{r_{i_1}+1} \cdots a_n^{-1} a_{n-1}^{r_{i_1}+t} AB^{-1} \underline{q}_n a_n^{-1} Q_1''.$$

By the assumption that the underlined letters $\underline{q}_n$ and $\underline{q}_n^{-1}$ are the first to be deleted in (17), and that in (17) no transformations are performed on them, we have the independent sequence of transformations

$$BA^{-1} a_{n-1}^{r_{i_1}} a_n^{-1} a_{n-1}^{r_{i_1}+1} \cdots a_n^{-1} a_{n-1}^{r_{i_1}+t} AB^{-1} \to \cdots \to 1$$

with index 0. By Lemma 1, §1, Chapter II, this implies that

$$BA^{-1} a_{n-1}^{r_{i_1}} a_n^{-1} a_{n-1}^{r_{i_1}+1} \cdots a_n^{-1} a_{n-1}^{r_{i_1}+t} AB^{-1} = 1$$

in the group $\mathfrak{A}_P$. Hence by Lemma 7, §2, Chapter I we have

$$(a_{n-1}^{r_{i_1}} a_n^{-1} a_{n-1}^{r_{i_1}+1} \cdots a_n^{-1} a_{n-1}^{r_{i_1}+t})_{\overline{S}} = a_n^{-t} = 1 \quad \text{in } \mathfrak{A}_P.$$

Consequently, we have $t = 0$. This contradicts the assumption that $q_n^{-1}$ and $q_n$ appear in different segments $V_i$.

**Part IV.** Again, we will assume that we have a sequence of elementary transformations with index 0:

$$Q_2 \equiv U_1 q_n G_1 AB^{-1} U_2 q_n G_2 AB^{-1} \cdots U_k q_n G_k AB^{-1} U_{k+1} \to \cdots \to 1, \qquad (18)$$

where $k \geq 1$, and then prove by induction on the number $m$ of letters $q_n^\sigma$ in the word $Q_2$ that such an assumption leads to a contradiction.

The word of type $Q_2$ always contains at least one letter $q_n$, as $k \geq 1$. Suppose $Q_2$ in the sequence (18) has the minimal number of letters $q_n^\sigma$ possible, i.e. a single letter $q_n$. Then

$$Q_2 \equiv U_1 q_n G_1 AB^{-1} U_2,$$

where $G_1$ does not contain letters $q_n^\sigma$.

As $\lambda = 0$, in the sequence (18) no new letters $q_n^\sigma$ can appear, and consequently the explicitly written letter $q_n$ in $Q_2$ cannot be deleted. This is a contradiction.

Suppose that we have proved the impossibility of the sequence (18) for all words $Q_2$ whose number $m$ of letters $q_n^\sigma$ is less than some $m_0 > 0$. We will now prove the impossibility of the sequence (18) when $m = m_0$. Consider the letters $q_n$ and $q_n^{-1}$, contained in $Q_2$, which are deleted first in the sequence (18). We have two possibilities:

1. both these letters appear in one and the same segment $G_i$;
2. the letter $q_n$ is one of the letters $q_n$ distinguished in (18), and $q_n^{-1}$ appears in $G_i$.

The letters $q_n^\sigma$ which are the first to be deleted cannot be taken from different segments $G_i$, as then between them will be one of the letters $q_n^\sigma$ distinguished in (18).

In **Case 1**, we repeat the reasoning which was given when considering **III**, and we get a contradiction to the inductive hypothesis. The fact that in the place of the letters $a_n^{-1}$ we here have letters $a_n$ is immaterial to the proof.

**Case 2**. We divide this into two subcases.

$\alpha$) the deletion is right-sided,

$\beta$) the deletion is left-sided.

**Subcase $\alpha$).** The letter $q_n^{-1}$ appears in the segment $G_i$ of the word $Q_2$, and to the left of the letter $q_n$ with which it deletes. As the considered deletion is the first deletion to occur in (18), $G_i$ has the form

$$G_i \equiv G_i' a_n q_n^{-1} B A^{-1} a_{n-1}^{r_1} a_n a_{n-1}^{r_2} \cdots a_n a_{n-1}^{r_\nu}.$$

Then the sequence (18) is written as:

$$U_1 q_n G_1 A B^{-1} \cdots U_i q_n G_i' a_n \underline{q}_n B A^{-1} a_{n-1}^{r_1} a_n a_{n-1}^{r_2} \cdots$$
$$\cdots a_n a_{n-1}^{r_\nu} A B^{-1} U_{i+1} \underline{q}_n G_{i+1} A B^{-1} \cdots U_k q_n G_k A B^{-1} U_{k+1} \to \cdots$$
$$\cdots \to Z_1 \underline{q}_n^{-1} \underline{q}_n Z_2 \to Z_1 Z_2 \to \cdots \to 1, \qquad (18')$$

where the underlined letters $\underline{q}_n$ and $\underline{q}_n^{-1}$ delete with one another. As in (18') no transformations, other than deletions, are performed on the letters $q_n^\sigma$, we have the following independent sequences of transformations with index 0:

1. $U_1 q_n G_1 A B^{-1} \cdots U_i q_n G_i' a_n \to \cdots \to Z_1$,
2. $G_{i+1} A B^{-1} \cdots U_k q_n G_k A B^{-1} U_{k+1} \to \cdots \to Z_2$,
3. $Z_1 Z_2 \to \cdots \to 1$.

By combining the sequences 1, 2, and 3, we obtain a sequence of elementary transformations with index $\lambda = 0$:

$$Q_{2,0} \equiv U_1 q_n G_1 A B^{-1} \cdots U_i q_n G_i' a_n G_{i+1} A B^{-1} \cdots U_k q_n G_k A B^{-1} U_{k+1} \to \cdots \to$$
$$\to Z_1 G_{i+1} A B^{-1} U_{i+2} q_n \cdots U_k q_n G_k A B^{-1} U_{k+1} \to \cdots \to Z_1 Z_2 \to \cdots \to 1.$$

As $Q_{2,0}$ is a word of type $Q_2$ and contains $m_0 - 2$ of the letters $q_n^\sigma$, we hence have a contradiction to the inductive hypothesis.

**Subcase $\beta$).** The deletion is left-sided. Then the segment $G_i$, whence $\underline{q}_n^{-1}$ was taken, lies to the right of the letter $\underline{q}_n$ and has the form

$$G_i \equiv a_{n-1}^{r_1} a_n a_{n-1}^{r_2} \cdots a_n a_{n-1}^{r_{\nu+1}} a_n \underline{q}_n BA^{-1} G_i''.$$

The sequence (18) hence has the form

$$U_1 q_n G_1 A B^{-1} \cdots U_i \underline{q}_n a_{n-1}^{r_1} a_n a_{n-1}^{r_2} \cdots a_n a_{n-1}^{r_{\nu+1}} a_n \underline{q}_n^{-1} BA^{-1} G_i'' AB^{-1} U_{i+1} q_n \cdots$$
$$\cdots U_k q_n G_k A B^{-1} U_{k+1} \to \cdots \to Z_1 \underline{q}_n \underline{q}_n^{-1} Z_2 \to Z_1 Z_2 \to \cdots \to 1.$$

As no transformations, other than deletion, are performed on the $\underline{q}_n^\sigma$, it follows that

$$a_{n-1}^{r_1} a_n a_{n-1}^{r_2} \cdots a_n a_{n-1}^{r_{\nu+1}} a_n \to \cdots \to 1$$

with index 0. This is a sequence of transformations in $\mathfrak{A}_{q_{n-1}, A, B}$. By Lemma 1, §1, Chapter II it follows that in the group $\mathfrak{A}_P$ we have the equality

$$a_{n-1}^{r_1} a_n a_{n-1}^{r_2} \cdots a_n a_{n-1}^{r_{\nu+1}} a_n = 1,$$

whence by Lemma 1, §2, Chapter II we find

$$(a_{n-1}^{r_1} a_n a_{n-1}^{r_2} \cdots a_n a_{n-1}^{r_{\nu+1}} a_n)_{\overline{S}} = a^{\nu+1} = 1 \quad \text{in } \mathfrak{A}_P.$$

But this is not possible, as $\nu + 1 \geq 1$.

Hence, we have shown the lemma when $\lambda = 0$.

Suppose the lemma is proved for all $\lambda < \lambda_0$, where $\lambda_0 \geq 1$. We prove it for $\lambda = \lambda_0$.

**Part I.** The sequence (15) has index $\lambda_0$. Clearly, **I** will be proved, if we can find a new sequence of transformations changing the word $E_1$ into $E_m$, and having index $\lambda_1 < \lambda_0$. If there are no insertions of letters $q_n^\sigma$ in (15), then there is nothing to prove.

Suppose, then, that there are such insertions in (15). Let us say that the last of these is the transformation $E_i \to E_{i+1}$.

We will first consider the case when the insertion is **right-sided**. Then

$$E_i \equiv E_i' E_i''; \quad E_{i+1} \equiv E_i' q_n^{-1} q_n E_i''.$$

The sequence (15) will have the form

$$E_1 \to \cdots \to E_i' E_i'' \to E_i' q_n^{-1} q_n E_i'' \to \cdots \to E_m. \tag{15'}$$

As the word $E_m$ by assumption does not contain the letter $q_n^{-1}$, the letter $q_n^{-1}$ must eventually be deleted in (15'). There are three cases to consider:

1°. The considered letter $q_n^{-1}$ deletes with some $q_n$ appearing in (15') to its left.

2°. The considered letter $q_n^{-1}$ deletes with a letter $q_n$ whose individuality is the same as that of the letter $q_n$ with which the $q_n^{-1}$ was inserted.

3°. The letter $q_n^{-1}$ from the final insertion deletes with a $q_n$ which appears in (15') to its right, but which is distinct in individuality from that which it was inserted with.

**Case 1°.** The sequence (15') has the form

$$E_1 \to \cdots \to E_i' E_i'' \to E_i' \underline{q_n}^{-1} q_n E_i'' \to \cdots$$
$$\cdots \to E_j' q_n \underline{q_n}^{-1} E_j'' \to E_j' E_j'' \to \cdots \to E_m, \tag{19}$$

where no transformations, except deletions, are performed on the letter $\underline{q_n}^{-1}$. Then from (19) we have the following independent transformations:

1. $E_1 \to \cdots \to E_i' E_i''$     3. $q_n E_i'' \to \cdots \to E_j''$,

2. $E_i' \to \cdots \to E_j' q_n$,     4. $E_j' E_j'' \to \cdots \to E_m$.

By combining these, we obtain from (19) a new sequence (19'), in which the same transformations are performed in a different order, but without the final insertion of the letter $q_n^\sigma$:

$$E_1 \to \cdots \to E_i' E_i'' \to \cdots \to E_j' q_n E_i'' \to \cdots \to E_j' E_j'' \to \cdots \to E_m. \tag{19'}$$

The sequence (19') now clearly has index $\lambda_0 - 1$.

**Case 2°.** Write the part of (15') which begins at the moment of the final insertion:

$$E_i' \underline{q_n}^{-1} \underline{q_n} E_i'' \to \cdots \to E_j' \underline{q_n}^{-1} \underline{q_n} E_j'' \to E_j' E_j'' \to \cdots \to E_m. \tag{20}$$

Replace the sequence (20) by a canonical sequence with respect to the letter $\underline{q_n}$.

$$E_i' \underline{q_n}^{-1} \underline{q_n} E_i'' \to \cdots \to E_i' \underline{q_n}^{-1} U \underline{q_n} V E_i'' \to \cdots$$
$$\cdots \to E_j' \underline{q_n}^{-1} \underline{q_n} E_j'' \to E_j' E_j'' \to \cdots \to E_m. \tag{20'}$$

From (20') we have the following independent sequences of transformations without insertions of the letter $q_n^\sigma$:

1. $\underline{q_n} \to \cdots \to U q_n V$     3. $V E_i'' \to \cdots \to E_j''$,

2. $E_i' \to \cdots \to E_j'$,     4. $U \to \cdots \to 1$. $\tag{20''}$

As the word $U$ does not contain any letters $q_n^\sigma$, sequence 4 is a sequence of elementary transformations in the group $\mathfrak{A}_{q_{n-1},A,B}$. By virtue of Corollary 2 of Lemma 2, §1, Chapter II we have that the sequence of exponents corresponding to the word $U$ is degenerate. Then the word $W$, corresponding to $U$ under $L_{q_n}$, and consequently also

$V = W^{-1}$, is equal to 1 in the free group. This implies that we can find a sequence consisting only of insertions:

$$1 \to \cdots \to V. \tag{21}$$

The index of the sequence (21) is no greater than the number of letters $q_n^\sigma$ in the word $V$, which is in turn no greater than the index of the transition $q_n \to \cdots \to U q_n V$.

By combining (15′), (21), and (20″), we find a new sequence of transformations

$$E_1 \to \cdots \to E_i' E_i'' \to \cdots \to E_i' V E_i'' \to \cdots \to E_j' V E_i'' \to \cdots$$
$$\cdots \to E_j' E_j'' \to \cdots \to E_m, \tag{22}$$

the index of which is strictly less than $\lambda_0$. Indeed, all transformations of this sequence are also in the sequence (15′), with the exception of the transition

$$E_i' E_i'' \to \cdots \to E_i' V E_i'',$$

the index of which is strictly less than the index of the transformation

$$E_i' E_i'' \to E_i' q_n^{-1} q_n E_i'' \to \cdots \to E_i' q_n^{-1} U q_n V E_i'',$$

which is present in (15′) but absent in (22). However, the index of any sequence of transformations does not depend on the total number of transformations, but is only dependent on the number of insertions of letters $q_n^\sigma$ and transformations from the schema (3′). Consequently, the index of the sequence (22) is indeed strictly less than the index of the sequence (15′), i.e. less than $\lambda_0$.

**Case 3°.** The letter $q_n^{-1}$ deletes with some letter $q_n$ which appears to the right of it, but which is distinct from that letter $q_n$ together with which it was inserted. The part of the sequence (15′), beginning at the moment of the specified insertion and continuing until the end, has the form

$$E_i' q_n^{-1} q_n E_i'' \to \cdots \to E_j' q_n^{-1} \underline{q}_n E_j'' \to E_j' E_j'' \to \cdots \to E_m. \tag{23}$$

Denote the index of the sequence (23) by $\lambda'$, and the index of the sequence

$$E \to \cdots \to E_i' E_i'' \tag{24}$$

by $\lambda''$. Then we have $\lambda' + \lambda'' + 1 = \lambda_0$. The letter $\underline{q}_n$, which deletes in (23) with the specified letter $q_n^{-1}$, is either contained in the word $E_i''$, or else is engendered in the sequence (23) from one of the letter $q_n$ appearing in the segment $q_n E_i''$. There are two possible cases:

1) the letter $\underline{q}_n$ is engendered by the letter $q_n$ from the last insertion;
2) the letter $\underline{q}_n$ is engendered in (23) by some letter $q_n$ from the word $E_i''$, or is itself contained in $E_i''$.

**Case 1).** Replace the sequence (23) by a canonical sequence with respect to the hereditary chain which engenders the letter $\underline{q}_n$ (see Lemma 2 of the present section). We have

$$E_i' q_n^{-1} q_n E_i'' \to \cdots \to E_i' \underline{q}_n^{-1} U_1 q_n V_1' A B^{-1} U_2 q_n V_2' \cdots$$
$$\cdots V_{k-1}' A B^{-1} U_k \underline{q}_n V_k a_n^{-1} V_{k-1} a_n^{-1} \cdots a_n^{-1} V_2'' a_n^{-1} V_1'' E_i'' \to \cdots$$
$$\cdots \to E_j' \underline{q}_n^{-1} \underline{q}_n E_j'' \to E_j' E_j'' \to \cdots \to E_m. \tag{23′}$$

In (23′) the underlined letters, which delete with one another, are unaffected by transformations from the moment they are underlined until the moment they are deleted.

Hence the word lying between them is independently transformed into 1. This yields a sequence of transformations

$$U_1 q_n V_1' A B^{-1} U_2 q_n V_2' \cdots V_{k-1}' A B^{-1} U_k \to \cdots \to 1$$

with index strictly less than $\lambda_0$, as the sequence (23$'$) has index $\lambda < \lambda_0$ (cf. Remark 2 on p. 46). This yields a contradiction to the inductive hypothesis (for **IV**). In other words, Case 1) is impossible.

**Case 2).** The word $E_i''$ has the form $E_{i,1}'' q_n E_{i,2}''$, where the distinguished letter $q_n$, which is either engendered in (23) from the letter $\underline{q}_n$ which deletes with the letter $q_n^{-1}$ from the last insertion, or else it itself deletes with the letter $q_n^{-1}$ from the last insertion. Replace the sequence (23) by a canonical sequence with respect to the hereditary chain which engenders the letter $\underline{q}_n$. This yields the sequence

$$E_i' q_n^{-1} E_{i,1}'' q_n E_{i,2}'' \to \cdots \to E_i' \underline{q}_n^{-1} q_n E_{i,1}'' U_1 q_n V_1' A B^{-1} U_2 q_n V_2' \cdots$$
$$\cdots V_{k-1}' A B^{-1} U_k \underline{q}_n V_k a_n^{-1} V_{k-1}'' \cdots V_2'' a_n^{-1} V_1'' E_{i,2}'' \to \cdots$$
$$\cdots \to E_j' \underline{q}_n^{-1} \underline{q}_n E_j'' \to E_j' E_j'' \to \cdots \to E_m \qquad (23'')$$

which has index $\lambda'$.

From (23$''$) we find the following independent transitions:

1. $q_n \to \cdots \to U_1 q_n V_1' A B^{-1} U_2 q_n V_2' ... V_{k-1}' A B^{-1} U_k q_n V_k a_n^{-1} V_{k-1}'' ... V_2'' a_n^{-1} V_1''$,
2. $q_n E_{i,1}'' U_1 q_n V_1' A B^{-1} U_2 q_n V_2' \cdots V_{k-1}' A B^{-1} U_k \to \cdots \to 1$,
3. $V_k a_n^{-1} V_{k-1}'' \cdots V_2'' a_n^{-1} V_1'' E_{i,2}'' \to \cdots \to E_j''$,
4. $E_i' \to \cdots \to E_j'$,
5. $E_j' E_j'' \to \cdots \to E_m$.

The sum of the indices of these sequences equals $\lambda'$. By virtue of Lemma 3 of the present section applied to sequence 2, we have the sequence

$$2'. \quad E_{i,1}'' U_1 q_n V_1' A B^{-1} U_2 q_n V_2' \cdots V_{k-1}' A B^{-1} U_k q_n \to \cdots \to 1$$

which has the same index as sequence 2.

Combining the transformations (24), 1, 2$'$, 3, 4, and 5, we find the following sequence of elementary transformations:

$$E_1 \to \cdots \to E_i' E_{i,1}'' q_n E_{i,2}'' \to \cdots \to E_i' E_{i,1}'' U_1 q_n V_1' A B^{-1} U_2 q_n V_2' \cdots$$
$$\cdots U_{k-1} A B^{-1} U_k q_n V_k a_n^{-1} V_{k-1}'' \cdots V_2'' a_n^{-1} V_1'' E_{i,2}'' \to \cdots$$
$$\cdots \to E_i' V_k a_n^{-1} V_{k-1}'' \cdots V_2'' a_n^{-1} V_1'' E_{i,2}'' \to \cdots \to E_i' E_j'' \to \cdots$$
$$\cdots \to E_j' E_j'' \to \cdots \to E_m. \qquad (24')$$

The sequence (24$'$) transforms $E_1$ into $E_m$ and has index $\lambda'' + \lambda' = \lambda_0 - 1 < \lambda_0$.

**The case of a left-sided insertion.** The sequence (15), in this case, appears as:

$$E_1 \to \cdots \to E_i' E_i'' \to E_i' q_n q_n^{-1} E_i'' \to \cdots \to E_m. \qquad (15'')$$

The letter $q_n^{-1}$ from the final insertion in the sequence (15$''$) must eventually be deleted. There are three possible cases:

1°. The letter $q_n^{-1}$ deletes with some $q_n$ which in (15$''$) appears to the right of it.
2°. The letter $q_n^{-1}$ deletes with the letter $q_n$ with which it was inserted.
3°. The letter $q_n^{-1}$ is deleted on the left, but with a letter $q_n$ distinct in individuality from that of the letter with which it was inserted.

**Case 1°.** This is analogous to Case 1° for right-sided insertions, and is omitted.

**Case 2°.** The sequence (15″) has the form

$$E_1 \to ... \to E_i'E_i'' \to E_i'\underline{q}_n\underline{q}_n^{-1}E_i'' \to ... \to E_j'\underline{q}_n\underline{q}_n^{-1}E_j'' \to E_j'E_j'' \to ... \to E_m, \quad (25)$$

where the distinguished letters $\underline{q}_n$ have one and the same individuality. The part of this sequence, starting from the moment of the last insertion and continuing until the end, has the form

$$E_i'\underline{q}_n\underline{q}_n^{-1}E_i'' \to \cdots \to E_j'\underline{q}_n\underline{q}_n^{-1}E_j'' \to E_j'E_j'' \to \cdots \to E_m. \quad (26)$$

Consider a canonical sequence for (26), with respect to the underlined individuality $\underline{q}_n$,

$$E_i'\underline{q}_n\underline{q}_n^{-1}E_i'' \to \cdots \to E_i'U\underline{q}_n V\underline{q}_n^{-1}E_i'' \to \cdots \to E_j'\underline{q}_n\underline{q}_n^{-1}E_j'' \to E_j'E_j'' \to \cdots \to E_m.$$

Here we have the following independent transitions without insertions of letters $q_n^\sigma$:

1. $E_i'U \to \cdots \to E_j'$,        3. $E_i'' \to \cdots \to E_j''$,
2. $V \to \cdots \to 1$,        4. $E_j'E_j'' \to \cdots \to E_m$.

The index of transition **2** is strictly less than $\lambda_0$. By the inductive hypothesis (for **II**) it follows that $U = 1$ in the free group. Hence there exists a sequence consisting only of insertions:

$$2'. \ 1 \to \cdots \to U,$$

which has index 0, as the word $U$ does not contain any letters $q_n^\sigma$. By combining the transitions **2′, 1, 3**, and **4**, we obtain the sequence of transformations:

$$E_1 \to \cdots \to E_i'E_i'' \to \cdots \to E_i'UE_i'' \to \cdots \to E_j'E_i'' \to \cdots \to E_j'E_j'' \to \cdots \to E_m.$$

Clearly, the index of this sequence is strictly less than $\lambda_0$, as it does not contain the final insertion of letters $q_n^\sigma$ from the sequence (25).

**Case 3°.** Just like when we considered the case of a right-sided insertion, here we have two subcases:

1) the $q_n$ which deletes with $q_n^{-1}$ is engendered from the $q_n$ of the final insertion;
2) the $q_n$ which deletes with $q_n^{-1}$ is either engendered from some letter $q_n$ which is contained in the word $E_i'$, or is itself contained in $E_i'$.

**Subcase 1).** The part

$$E_i'q_nq_n^{-1}E_i'' \to \cdots \to E_m$$

of the sequence (15″) has index $\lambda < \lambda_0$ and does not contain insertions of letters $q_n^\sigma$. Replace this part with its canonical subsequence with respect to the hereditary chain engendered from the letter $q_n$, which deletes with the distinguished letter $q_n^{-1}$. We obtain a sequence with index $\lambda < \lambda_0$ (Remark 2, 46):

$$E_i'q_nq_n'E_i'' \to \cdots \to E_i'U_1q_nV_1'AB^{-1}U_2q_nV_2' \cdots$$
$$\cdots AB^{-1}U_k\underline{q}_nV_ka_n^{-1}V_{k-1}'' \cdots a_n^{-1}V_2''a_n^{-1}V_1''\underline{q}_n^{-1}E_i'' \to \cdots$$
$$\cdots \to E_j'\underline{q}_n\underline{q}_n^{-1}E_j'' \to E_j'E_j'' \to \cdots \to E_m.$$

From this sequence, where the letters $q_n$ and $q_n^{-1}$ which delete with one another are underlined, we obtain the sequence

$$V_ka_n^{-1}V_{k-1}'' \cdots V_2''a_n^{-1}V_1'' \to \cdots \to 1$$

with index $\lambda < \lambda_0$. This is a contradiction to the inductive hypothesis (of **III**), as $k > 1$. Consequently, subcase 1) is impossible.

**Subcase 2)** is handled similarly as the analogous case for right-sided insertions.

**Part II.** By assumption, for $\lambda < \lambda_0$ the lemma is proved. Suppose the sequence

$$V \to \cdots \to 1 \tag{16}$$

has index $\lambda_0$ and does not contain insertions of letters $q_n^\sigma$. We then have to prove that the word $U$, corresponding to the word $V^{-1}$ under $L_{q_n}$, is equal to 1 in the free group.

The proof will be carried out by induction on the number $m$ of letters $q_n^\sigma$ appearing in the word $V$. If $m = 0$, i.e. the word $V$ does not contain any letters $q_n^\sigma$, then $\lambda_0 = 0$. This case has already been considered. Suppose then that for all words $V$, for which $m < m_0$, for $\lambda = \lambda_0$ part **II** has already been proved, and that the word $V$ in the sequence (16) contains $m_0$ letters $q_n^\sigma$. We may assume that $m_0 > 0$. The letters $q_n^\sigma$, appearing in the word $V$ and engendered in the sequence (16), must all be deleted in (16). Consider the first deleted pair of letters $q_n$ and $q_n^{-1}$ in (16). We have three cases:

1°. The first-deleted individualities $q_n^\sigma$ both appear in the word $V$.

2°. One of the two appears in $V$, and the other is engendered in (16) via an application of the schema (3′).

3°. Both are engendered in the sequence (16).

**Case 1°.** The word $V$ has the form

$$V \equiv a_{n-1}^{-r_{\mu+1}} (a_n q_n^{-1} B A^{-1})^{-s_\mu} a_{n-1}^{-r_\mu} \cdots a_{n-1}^{-r_2} (a_n q_n^{-1} B A^{-1})^{-s_1} a_{n-1}^{-r_1},$$

where we may assume that all $s_i$ are non-zero. The first-deleted letters $q_n$ and $q_n^{-1}$ in the sequence (16) must appear in neighbouring segments of type $(a_n q_n^{-1} B A^{-1})^{-s_j}$.

We have two possibilities:

1) the considered deletion is right-sided,

2) the considered deletion is left-sided.

**Case 1).** The word $V$ has the form

$$V_1 a_n q_n^{-1} B A^{-1} a_{n-1}^{-r_{j_0}} A B^{-1} \underline{q}_n a_n^{-1} V_2,$$

where the underlined letters $q_n^{-1}$ are the first to be deleted in (16). The word $U$, corresponding to $V^{-1}$ under $L_{q_n}$, will be $U_2 a_n q_{n-1}^{r_{j_0}} a_n^{-1} U_1$, where $U_1$ and $U_2$ correspond to the words $V_1^{-1}$ and $V_2^{-1}$ under $L_{q_n}$. The sequence (16) has the form

$$V_1 a_n \underline{q}_n^{-1} B A^{-1} a_{n-1}^{-r_{j_0}} A B^{-1} \underline{q}_n a_n^{-1} V_2 \to \cdots \to Z_1 \underline{q}_n^{-1} \underline{q}_n Z_2 \to Z_1 Z_2 \to \cdots \to 1. \tag{16′}$$

Replace (16′) by a canonical sequence with respect to the underlined letter $\underline{q}_n$. We obtain the sequence

$$V_1 a_n \underline{q}_n^{-1} B A^{-1} a_{n-1}^{-r_{j_0}} A B^{-1} \underline{q}_n a_n^{-1} V_2 \to \cdots \to V_1 a_n \underline{q}_n^{-1} B A^{-1} a_{n-1}^{-r_{j_0}} \cdots$$
$$A B^{-1} U' \underline{q}_n V' a_n^{-1} V_2 \to \cdots \to Z_1 \underline{q}_n^{-1} \underline{q}_n Z_2 \to \cdots \to Z_1 Z_2 \to \cdots \to 1 \tag{16″}$$

with index $\lambda_0$.

In this sequence we have the following independent sequences of transformations without insertions of letters $q_n^\sigma$:

| | | |
|---|---|---|
| 1. | $q_n \to \cdots \to U' q_n V'$ | with index $\lambda_1$, |
| 2. | $B A^{-1} a_{n-1}^{-r_{j_0}} A B^{-1} U' \to \cdots \to 1$ | »　　»　$\lambda_2$, |
| 3. | $V_1 a_n \to \cdots \to Z_1$ | »　　»　$\lambda_3$, |
| 4. | $V' a_n^{-1} V_2 \to \cdots \to Z_2$ | »　　»　$\lambda_4$, |
| 5. | $Z_1 Z_2 \to \cdots \to 1$ | »　　»　$\lambda_5$, |

where $\lambda_1 + \lambda_2 + \lambda_3 + \lambda_4 + \lambda_5 = \lambda_0$.

The word $U'$ has the form

$$q_{n-1}^{k_1} a_n^{t_1} q_{n-1}^{k_2} a_n^{t_2} \cdots q_{n-1}^{k_\nu} a_n^{t_\nu} q_{n-1}^{k_{\nu+1}},$$

and $V'^{-1}$ corresponds with the word $W$ under $L_{q_n}$. As in transition 2 there are no insertions of letters $q_n^\sigma$, we have $\lambda_2 = 0$ and the transition 2 is a sequence of transformations in the group $\mathfrak{A}_{q_{n-1}, A, B}$. By Lemma 2, §1, Chapter II, the sequence of numbers

$$k_1, t_1, k_2, t_2, \ldots, k_\nu, t_\nu, k_{\nu+1}, \tag{27}$$

corresponding to the word $U'$, is degenerate, and by virtue of the corollary of this lemma

$$(a_{n-1}^{-r_{j_0}})_{\overline{S^+}} \equiv a_{n-1}^{-r_{j_0}} = 1 \text{ in } \mathfrak{A}_P.$$

Hence it follows that $r_{j_0} = 0$.

From the degeneracy of the sequence of numbers (27), it follows that $V' = 1$ in the free group, and consequently we may find a sequence of transformations containing only insertions

$$1 \to \cdots \to V'. \tag{28}$$

The index $\lambda'$ of the sequence (28) is no greater than the number of letters $q_n^\sigma$ appearing in the word $V'$, which in turn is no greater than the index $\lambda_1$ of transition 1. Combining sequence (28) with the transitions 3, 4, and 5, we find a new sequence of transformations

$$V_1 V_2 \to V_1 a_n a_n^{-1} V_2 \to \cdots \to V_1 a_n V' a_n^{-1} V_2 \to \cdots$$
$$\cdots \to Z_1 V' a_n V_2 \to \cdots \to Z_1 Z_2 \to \cdots \to 1,$$

the index of which is equal to $\lambda' + \lambda_3 + \lambda_4 + \lambda_5 \le \lambda_0$. As the word $V_1 V_2$ contains $m_0 - 2$ letters $q_n^\sigma$, it follows by the inductive hypothesis that $U_2 U_1 = 1$ in the free group. As $r_{j_0} = 1$, the word $U$ is equal to $U_2 a_n a_n^{-1} U_1$ and is hence also equal to $1$ in the free group.

**Case 2).** The considered deletion is left-sided. In this case, the sequence (16) has the form

$$V \equiv V_1 A B^{-1} \underline{q}_n a_n^{-1} a_{n-1}^{-r_{i_0}} a_n \underline{q}_n^{-1} B A^{-1} V_2 \to \cdots$$
$$\cdots \to Z_1 \underline{q}_n \underline{q}_n^{-1} Z_2 \to \cdots Z_1 Z_2 \to \cdots \to 1, \tag{29}$$

where we have underlined the individualities which delete with one another.

Under $L_{q_n}$, the word $V^{-1}$ corresponds to the word $U = U_2 a_n^{-1} q_{n-1}^{-r_{i_0}} a_n U_1$, where $U_1, U_2$ correspond under $L_{q_n}$ with the words $V_1^{-1}$ and $V_2^{-1}$. Replace (29) by a canonical sequence with respect to the underlined letter $q_n$, yielding

$$V_1 A B^{-1} \underline{q}_n a_n^{-1} a_{n-1}^{-r_{i_0}} a_n \underline{q}_n^{-1} B A^{-1} V_2 \to \cdots$$
$$\cdots \to V_1 A B^{-1} U' \underline{q}_n V' a_n^{-1} \cdots a_{n-1}^{-r_{i_0}} a_n \underline{q}_n^{-1} B A^{-1} V_2 \to \cdots$$
$$\cdots \to Z_1 \underline{q}_n \underline{q}_n^{-1} Z_2 \to Z_1 Z_2 \to \cdots \to 1,$$

where we have the following independent transitions without insertions of letters $q_n^\sigma$:

$$\begin{array}{llll}
1. & V' a_n^{-1} a_{n-1}^{-r_{i_0}} a_n \to \cdots \to 1 & \text{with index } \lambda_1, \\[4pt]
2. & V_1 A B^{-1} U' \to \cdots \to Z_1 & \text{»} \quad \text{»} \quad \lambda_2, \\[4pt]
3. & B A^{-1} V_2 \to \cdots \to Z_2 & \text{»} \quad \text{»} \quad \lambda_3, \\[4pt]
4. & Z_1 Z_2 \to \cdots \to 1 & \text{»} \quad \text{»} \quad \lambda_4,
\end{array}$$

and with $\lambda_1 + \lambda_2 + \lambda_3 + \lambda_4 \leq \lambda_0$.

As the considered deletion occurs first, in transition 1 we do not have insertions of letters $q_n^\sigma$, and consequently the word $V'$ does not contain letters $q_n^\sigma$. Thus $V' = a_{n-1}^{-k_1}$ and $U' = q_{n-1}^{k_1}$. The sequence

$$1'. \quad a_{n-1}^{-k_1} a_n^{-1} a_{n-1}^{-r_{i_0}} a_n \to \cdots \to 1$$

does not contain insertions of letters $q_n^\sigma$. Hence, $a_{n-1}^{-k_1} a_n^{-1} a_{n-1}^{-r_{i_0}} a_n = 1$ in the group $\mathfrak{A}_{q_{n-1}, A, B}$, whence by Lemma 1, §1, Chapter II we have $a_{n-1}^{-k_1} a_n^{-1} a_{n-1}^{-r_{i_0}} a_n = 1$ in $\mathfrak{A}_P$. By Lemma 2, §2, Chapter I, we have $k_1 = 0$ and $r_{i_0} = 0$. Hence $U'$ equals 1 in the free group and from transition 2 we find

$$2'. \quad V_1 A B^{-1} \to \cdots \to V_1 A B^{-1} U' \to \cdots \to Z_1$$

with index $\lambda_2$.

Combining the sequences 2', 3, and 4, we obtain a sequence of transformations

$$V_1 V_2 \to \cdots \to V_1 A B^{-1} B A^{-1} V_2 \to \cdots \to Z_1 B A^{-1} V_2 \to \cdots Z_1 Z_2 \to \cdots \to 1$$

with index $\lambda = \lambda_2 + \lambda_3 + \lambda_4 \leq \lambda_0$.

The word $V_1 V_2$ contains $m_0 - 2$ letters $q_n^\sigma$. By the inductive hypothesis, we have $U_2 U_1 = 1$ in the free group. As $r_{i_0} = 0$, the word $U = U_2 a_n^{-1} q_{n-1}^{r_{i_0}} a_n U_2$ is also equal to 1 in the free group.

**Case** $2°$. One of the two letters $q_n^\sigma$, deleted first in (16), appears in $V$, and the other is engendered in (16).

We distinguish two possible subcases:

    1) the engendered letter is $\underline{q}_n$,

    2) the engendered letter is $\underline{q}_n^{-1}$.

**Subcase 1).** The letter $\underline{q}_n$, being the first to be deleted in (16), is engendered from some letter $q_n$ which appears in the word $V$. This letter $q_n$ in the word $V$ must lie next to the letter $\underline{q}_n^{-1}$ which in (16) deletes with $\underline{q}_n$. Clearly, it can lie only to the left of $\underline{q}_n^{-1}$, as otherwise it would lie between the letters $\underline{q}_n^{-1}$ and $\underline{q}_n$, and would have to be deleted before them. The word $V$ then has the form

$$V_1 A B^{-1} q_n a_n^{-1} a_{n-1}^{r_{i_1}} a_n \underline{q}_n^{-1} B A^{-1} V_2.$$

For the sequence (16), consider the canonical sequence with respect to the hereditary chain engendering the letter $\underline{q}_n$:

$$V_1 A B^{-1} q_n a_n^{-1} a_{n-1}^{r_{i_1}} a_n \underline{q}_n^{-1} B A^{-1} V_2 \to \cdots \to V_1 A B^{-1} U_1 q_n V_1' A B^{-1} U_2 q_n V_2' \cdots$$
$$\cdots U_k \underline{q}_n V_k a_n^{-1} V_{k-1}'' a_n^{-1} \cdots a_n^{-1} V_1'' a_n^{-1} a_{n-1}^{r_{i_1}} a_n \underline{q}_n^{-1} B A^{-1} V_2 \to \cdots$$
$$\cdots \to Z_1 \underline{q}_n \underline{q}_n^{-1} Z_2 \to Z_1 Z_2 \to \cdots \to 1. \quad (30)$$

As the individuality $\underline{q}_n$ does not appear in the word $V$, we have $k > 1$. In the sequence (30), no transformations, except deletions, are performed on the letters $\underline{q}_n$ and $\underline{q}_n^{-1}$ after the point that they are underlined. Hence, we have the following independent transition without insertions of letters $q_n^\sigma$:

$$V_k a_n^{-1} V_{k-1}'' a_n^{-1} \cdots a_n^{-1} V_1'' a_n^{-1} a_{n-1}^{r_{i_1}} a_n \to \cdots \to 1. \quad (31)$$

Consequently, the words $V_k$ and $V_j''$ $(j = 1, 2, \ldots, k-1)$ do not contain letters $q_n^\sigma$, i.e.

$$V_k = a_{n-1}^{s_k}, \quad V_j'' = a_{n-1}^{s_j}.$$

The transition (31) is a sequence of transformations in $\mathfrak{A}_{q_{n-1},A,B}$. By Lemma 1, §1, Chapter II, we have

$$a_{n-1}^{s_k}a_n^{-1}a_{n-1}^{s_{k-1}}a_n^{-1}\cdots a_n^{-1}a_{n-1}^{s_1}a_n^{-1}a_{n-1}^{r_{i_1}}a_n = 1$$

in the group $\mathfrak{A}_P$. By virtue of Lemma 1, §2, Chapter I, we have

$$(a_{n-1}^{s_k}a_n^{-1}a_{n-1}^{s_{k-1}}a_n^{-1}\cdots a_n^{-1}a_{n-1}^{s_1}a_n^{-1}a_{n-1}^{r_{i_1}}a_n)_{\overline{S}} = a_n^{-k+1} = 1 \quad \text{in } \mathfrak{A}_P,$$

which is impossible, as $k > 1$.

**Subcase 2).** We have two possibilities.

$\alpha$) The considered letter $\underline{q}_n^{-1}$ deletes with a letter $\underline{q}_n$ which is distinct in individuality from that which engendered it.

$\beta$) The considered letter $\underline{q}_n^{-1}$ deletes with that letter $\underline{q}_n$ which engendered it.

Suppose we are in $\alpha$). Then in the same way as we did in subcase 1), we reach the conclusion that in the group $\mathfrak{A}_P$ there is some equality

$$a_{n-1}^{s_k}BA^{-1}a_{n-1}^{s_{k-1}}a_n^{-1}a_{n-1}^{s_{k-2}}a_n^{-1}\cdots a_{n-1}^{s_1}a_n^{-1}a_{n-1}^{r_{i_1}}AB^{-1} = 1,$$

where $k > 1$.

By Lemma 7, §2, Chapter I we find

$$(a_{n-1}^{s_{k-1}}a_n^{-1}a_{n-1}^{s_{k-2}}a_n^{-1}\cdots a_{n-1}^{s_1}a_n^{-1}a_{n-1}^{r_{i_1}})_{\overline{S}} = a_n^{-k+1} = 1 \quad \text{in } \mathfrak{A}_P,$$

which is impossible, as $k > 1$.

In the case $\beta$), we distinguish in $V$ that letter $q_n$, which deletes with the letter $\underline{q}_n^{-1}$ engendered by it:

$$V \equiv V_1 AB^{-1}\underline{q}_n a_n^{-1}V_2.$$

Replace the sequence (16) by a canonical sequence with respect to the letter $\underline{q}_n$:

$$V_1 AB^{-1}\underline{q}_n a_n^{-1}V_2 \to \cdots \to V_1 AB^{-1}U'\underline{q}_n V_1' a_n \underline{q}_n^{-1} BA^{-1}V_1'' a_n^{-1}V_2 \to \cdots$$
$$\cdots \to Z_1 \underline{q}_n \underline{q}_n^{-1}Z_2 \to Z_1 Z_2 \to \cdots \to 1.$$

From this sequence we have a transition without deletions of the letters $q_n^\sigma$

$$V_1' a_n \to \cdots \to 1. \tag{32}$$

Consequently $V_1' = a_{n-1}^{s'}$, and (32) is a sequence of transformations in the group $\mathfrak{A}_{q_{n-1},A,B}$. By Lemma 1, §1, Chapter II, we have $a_{n-1}^{s'}a_n = 1$ in $\mathfrak{A}_P$, which is impossible by virtue of Lemma 4, §2, Chapter I.

Thus, all subcases of Case 2° have been proved to be impossible.

**Case 3°.** The letters $\underline{q}_n$ and $\underline{q}_n^{-1}$, which delete first in the sequence (16), are both engendered in (16). Clearly, they must be engendered from one and the same letter $q_n$, which is contained in the word $V$, because otherwise they cannot be the first to delete in the sequence (16). We write out in the word $V$ the letter $q_n^{(1)}$ which engenders them:

$$V \equiv V'AB^{-1}q_n^{(1)}a_n^{-1}V''.$$

Consider a canonical sequence for (16) associated with the letter $q_n^{(1)}$:

$$V'AB^{-1}q_n^{(1)}a_n^{-1}V'' \to \cdots \to V'AB^{-1}U_1 q_n^{(1)}V_1 a_n^{-1}V'' \to \cdots \to 1. \tag{33}$$

All letters $q_n^\sigma$, directly engendered from the letter $q_n^{(1)}$ in the sequence (33), appear in the word $V_1$. Every letter $q_n^\sigma$, engendered (but not directly engendered) from the letter $q_n^{(1)}$ in the sequence (33) is engendered in (33) from some letter $q_n$, which appears in the word $V_1$. Hence we can see that the letters $\underline{q}_n$ and $\underline{q}_n^{-1}$, the first to be deleted in (33),

either appear in $V_1$, or else are engendered from some letter $q_n$ which appears in $V_1$. We have three possible cases, analogous to the cases $1°, 2°,$ and $3°$:

$1^{\mathbf{1}}$. The letters $\underline{q}_n$ and $\underline{q}_n^{-1}$ both appear in the word $V_1$.

$2^{\mathbf{1}}$. One of the letters appears in $V_1$, and the other is engendered from some letter $q_n$ in the word $V_1$.

$3^{\mathbf{1}}$. Both letters are engendered from the same letter $q_n$, appearing in the word $V_1$.

Case $2^{\mathbf{1}}$, as we have already proved (cf. case $2°$), is impossible.

If we are in case $3^{\mathbf{1}}$, then we distinguish in the word $V_1$ the letter $q_n^{(2)}$ which engenders the letters $\underline{q}_n$ and $\underline{q}_n^{-1}$, i.e. the first two letters to be deleted in (33), and replace the sequence (33)

$$V'AB^{-1}U_1 q_n^{(1)} V_1' AB^{-1} q_n^{(2)} a_n^{-1} V_1'' a_n^{-1} V'' \to \cdots \to 1 \qquad (33')$$

by a canonical sequence for it with respect to the letter $q_n^{(2)}$. We obtain

$$V'AB^{-1}U_1 q_n^{(1)} V_1' AB^{-1} q_n^{(2)} a_n^{-1} V_1'' a_n^{-1} V' \to \cdots$$
$$\cdots \to V'AB^{-1}U_1 q_n^{(1)} V_1' AB^{-1} U_2 q_n^{(2)} V_2 a_n^{-1} V_1'' a_n^{-1} V'' \to \cdots \to 1.$$

Again, we find that the letters $\underline{q}_n$ and $\underline{q}_n^{-1}$, the first two to be deleted, either both appear in $V_2$, or else both engendered from the same letter $q_n$, appearing in $V_2$, etc.

Clearly, we find in a finite number of steps a word $V_k$, engendered by a letter $q_n^{(k)}$, in which we find both letters $\underline{q}_n$ and $\underline{q}_n^{-1}$, i.e. the first two letters to delete in (33). In doing so, we find a hereditary chain of letters

$$q_n^{(1)}, q_n^{(2)}, q_n^{(3)}, \ldots, q_n^{(k)}, \qquad (34)$$

engendered from each other in the sequence (33). Replace the sequence (33) by a canonical sequence for it with respect to the hereditary chain (34),

$$V'AB^{-1}q_n^{(1)} a_n^{-1} V'' \to \cdots \to V'AB^{-1}U_1 q_n^{(1)} V_1' AB^{-1} U_2 q_n^{(2)} V_2' \cdots$$
$$\cdots AB^{-1} U_k q_n^{(k)} \underline{V}_k a_n^{-1} \cdots V_1'' a_n^{-1} V'' \to \cdots \to 1. \qquad (35)$$

The sequence (35), just as (33), has index $\lambda_0$.

Based on the above, we have that the two letters $\underline{q}_n$ and $\underline{q}_n^{-1}$, which delete first in the sequence (35), both appear in the segment $\underline{V}_k$. As in case $1°$, here we have two possible subcases:

1) the considered deletion is right-sided,

2) the considered deletion is left-sided.

**Subcase 1).** The word $\underline{V}_k$ has the form

$$V_k' a_n \underline{q}_n^{-1} BA^{-1} a_{n-1}^{r'} AB^{-1} \underline{q}_n a_n^{-1} V_k''.$$

Replace the second part of the second (35) by a canonical sequence with respect to the letter $\underline{q}_n$,

$$V'AB^{-1}U_1 q_n^{(1)} V_1' AB^{-1} U_2 q_n^{(2)} V_2' \cdots$$
$$\cdots AB^{-1} U_k q_n^{(k)} V_k' a_n \underline{q}_n^{-1} BA^{-1} a_{n-1}^{r'} AB^{-1} \underline{q}_n a_n^{-1} V_k'' a_n^{-1} V_{k-1}'' \cdots$$
$$\cdots a_n^{-1} V_1'' a_n^{-1} V'' \to \cdots \to V'AB^{-1}U_1 q_n^{(1)} V_1' AB^{-1} U_2 q_n^{(2)} V_2' \cdots$$
$$\cdots AB^{-1} U_k q_n^{(k)} V_k' a_n \underline{q}_n^{-1} BA^{-1} a_{n-1}^{r'} AB^{-1} U_{k+1} \underline{q}_n V_{k+1} a^{-1} V_k'' a_n^{-1} \cdots$$
$$\cdots V_1'' a_n^{-1} V'' \to \cdots \to Z_1 \underline{q}_n^{-1} \underline{q}_n Z_2 \to Z_1 Z_2 \to \cdots \to 1. \qquad (35')$$

In the sequence $(35')$ we have the following independent transition without insertions of letters $q_n^\sigma$:

$$BA^{-1}a_{n-1}^{r'}AB^{-1}U_{k+1} \to \cdots \to 1,$$

whence, analogously to case $1°$, we find that $r' = 0$ and that there is a sequence of elementary transformations, consisting only of insertions:

$$1 \to \cdots \to a_n V_{k+1} a_n^{-1}. \tag{36}$$

The index of the sequence $(36)$ is no greater than the index of the transition

$$1. \quad \underline{q}_n \to \cdots \to U_{k+1}\underline{q}_n V_{k+1},$$

which occurs at the beginning of the sequence $(35')$. From $(35')$, moreover, we have the transitions:

$$2. \quad V'AB^{-1}U_1 q_n^{(1)} V_1' AB^{-1} U_2 q_n^{(2)} V_2' \cdots AB^{-1} U_k q_n^{(k)} V_k' a_n \to \cdots \to Z_1,$$

$$3. \quad V_{k+1}a_n^{-1}V_k'' a_n^{-1} \cdots V_1'' a_n^{-1}V'' \to \cdots \to Z_2,$$

$$4. \quad Z_1 Z_2 \to \cdots \to 1.$$

From $r' = 0$ it follows that

$$V_k = V_k' a_n q_n^{-1} BA^{-1} AB^{-1} q_n a_n^{-1} V_k''.$$

Thus

$$U_k = U_k'' a_n a_n^{-1} U_k'.$$

By combining the sequences of transformations $(35)$, $(36)$, and transition $2$, $3$, and $4$, we find a new sequence

$$V \equiv V'AB^{-1}q_n^{(1)} a_n^{-1} V'' \to \cdots \to V'AB^{-1}U_1 q_n^{(1)} V_1' AB^{-1} U_2 q_n^{(2)} V_2' \cdots$$

$$\cdots V_{k-1}'AB^{-1}U_k'' U_k' q_n^{(k)} V_k' V_k'' a_n^{-1} V_{k-1}'' \cdots V_1'' a_n^{-1} V'' \to \cdots$$

$$\cdots \to V'AB^{-1}U_1 q_n^{(1)} V_1' AB^{-1} U_2 q_n^{(2)} V_2' \cdots$$

$$\cdots V_{k-1}'AB^{-1} \underbrace{U_k'' a_n a_n^{-1} U_k'}_{U_k} q_n^{(k)} V_k' a_n V_{k+1} a_n^{-1} V_k'' a_n^{-1} \cdots$$

$$\cdots V_1'' a_n^{-1} V'' \to \cdots \to Z_1 V_{k+1} a_n^{-1} V_k'' a_n^{-1} \cdots V_1'' a_n^{-1} V'' \to \cdots \to Z_1 Z_2 \to \cdots \to 1. \tag{37}$$

In the sequence $(37)$ the letter $q_n^{(k)}$ appears in the word $V_k' V_k''$, which contains two fewer letters $q_n^\sigma$ fewer than the word $V_k$. All transformations of the sequence $(37)$ appear also in $(35)$, except for the insertion of $a_n a_n^{-1}$ and the transition $(36)$, the index of which is no greater than the transition $1$, which is missing from $(37)$. Consequently, the index of the sequence $(37)$ is no greater than $\lambda_0 - 2$, and by the inductive hypothesis the word $U$ is equal to $1$ in the free group.

**Subcase 2).** The word $V_k$ has the form

$$V_k AB^{-1} \underline{q}_n a_n^{-1} a_{n-1}^{r'} a_n \underline{q}_n^{-1} BA^{-1} V_k''.$$

This has as corresponding word

$$U_k = U_k'' a_n^{-1} q_{n-1}^{-r'} a_n U_k'.$$

Considering a canonical sequence for the first part of (35), with respect to the letter $\underline{q}_n$, we find

$$V'AB^{-1}U_1q_n^{(1)}V_1'AB^{-1}U_2q_n^{(2)}V_2'\cdots$$
$$\cdots AB^{-1}U_kq_n^{(k)}V_k'AB^{-1}\underline{q}_n a_n^{-1}a_{n-1}^{r'}a_n\underline{q}_n^{-1}BA^{-1}V_k''a_n^{-1}V_{k-1}''\cdots$$
$$\cdots V_1''a_n^{-1}V'' \to \cdots \to V'AB^{-1}U_1q_n^{(1)}V_1'AB^{-1}U_2q_n^{(2)}V_2'\cdots$$
$$\cdots AB^{-1}U_kq_n^{(k)}V_k'AB^{-1}U_{k+1}\underline{q}_n V_{k+1}a_n^{-1}a_{n-1}a_n\underline{q}_n^{-1}BA^{-1}V_k''a_n^{-1}V_{k-1}''\cdots$$
$$\cdots V_1''a_n^{-1}V'' \to \cdots \to Z_1\underline{q}_n\underline{q}_n^{-1}Z_2 \to Z_1Z_2 \to \cdots \to 1. \tag{38}$$

Here the underlined letters $\underline{q}_n$ and $\underline{q}_n^{-1}$ are the first to be deleted. Consequently, the word $V_{k+1}$ does not contain letters $q_n^\sigma$, i.e. $V_{k+1} = a_{n-1}^{r_{k+1}}$. Thus $U_{k+1} = q_{n-1}^{-r_{k+1}}$.

From the sequence (38) we have the following independent transitions without insertions of letters $q_n^\sigma$:

1. $a_{n-1}^{r_{k+1}}a_n^{-1}a_{n-1}^{r'}a_n \to \cdots \to 1$.
2. $V'AB^{-1}U_1q_n^{(1)}V_1'AB^{-1}U_2q_n^{(2)}V_2'\cdots AB^{-1}U_kq_n^{(k)}V_k'AB^{-1}U_{k+1} \to \cdots \to Z_1$,
3. $BA^{-1}V_k''a_n^{-1}\cdots V_1''a_n^{-1}V'' \to \cdots \to Z_2$.

Transition 1 is a sequence of transformations in the group $\mathfrak{A}_{q_{n-1},A,B}$. By Lemma 1, §1, Chapter II, we have $a_{n-1}^{r_{k+1}}a_n^{-1}a_{n-1}^{r'}a_n = 1$ in $\mathfrak{A}_P$, whence by Lemma 2, §2, Chapter I we have $r' = 0$ and $r_{k+1} = 0$. In this case, the word $U_{k+1}$ equals 1 in the free group, and we have a sequence of transformations

$$1 \to \cdots \to AB^{-1}U_{k+1}BA^{-1} \tag{36'}$$

with index 0.

From $r' = 0$ it follows that

$$U_k = U_k''a_n^{-1}a_nU_k'.$$

By combining the beginning of the sequence (35), the beginning and ends of the sequence (38), and the transitions (36'), 2, and 3, we find the sequence

$$V \equiv V'AB^{-1}q_n^{(1)}a_n^{-1}V'' \to \cdots \to V'AB^{-1}U_1q_n^{(1)}V_1'AB^{-1}U_2q_n^{(2)}V_2'\cdots$$
$$\cdots AB^{-1}U_k''U_k'a_n^{(k)}V_k'V_k''a_n^{-1}V_{k-1}''\cdots V_1''a_n^{-1}V'' \to \cdots$$
$$\cdots \to V'AB^{-1}U_1q_n^{(1)}V_1'AB^{-1}U_2q_n^{(2)}V_2'\cdots$$
$$\cdots AB^{-1}\underbrace{U_k''a_n^{-1}a_nU_k'}_{U_k}q_n^{(k)}V_k'AB^{-1}U_{k+1}BA^{-1}V_k''a_n^{-1}V_{k-1}''\cdots$$
$$\cdots V_1''a_n^{-1}V'' \to \cdots \to Z_1BA^{-1}V_k''a_n^{-1}V_{k-1}''\cdots$$
$$\cdots V_1''a_n^{-1}V'' \to \cdots \to Z_1Z_2 \to \cdots \to 1. \tag{37'}$$

The sequence (37') has index $\lambda \leq \lambda_0 - 2$, which can be seen just as in the case of sequence (37). Consequently, the word $U$ equals 1 in the free group.

**Part III.** Assume that for $\lambda < \lambda_0$ the lemma is proved. We have proved that under this assumption, parts **I** and **II** for $\lambda = \lambda_0$ are also true.

Suppose we have a sequence of elementary transformations in the group $\mathfrak{A}_{q,A,B}$ with index $\lambda_0$:

$$Q_1 = V_1a_n^{-1}V_2a_n^{-1}\cdots V_{k-1}a_n^{-1}V_k \to \cdots \to 1, \tag{39}$$

where $k > 1$. By virtue of the fact that we have proved part **I** for $\lambda = \lambda_0$, the existence of sequence (39) implies that there exists a sequence of elementary transformations

$$V_1 a_n^{-1} V_2 a_n^{-1} \cdots V_{k-1} a_n^{-1} V_k \to \cdots \to 1 \tag{39'}$$

without insertions of letters $q_n^\sigma$, and with index $\lambda \leq \lambda_0$. If the number $m$, of letters $q_n^\sigma$ appearing in the word $Q_1$, is equal to 0, then the sequence (39′) has index $\lambda = 0$. But for $\lambda = 0$, part **III** is already proved. Suppose then, that for all words of type $Q_1$ with $m < m_0$, where $m_0 > 0$, part **III** has been proved, i.e. for all such words we have proved the impossibility of a sequence (39) with index $\lambda \leq \lambda_0$. Suppose the given word $Q_1$ contains $m_0$ letters $q_n^\sigma$, and that the sequence (39′) has index $\lambda_0$. We will show that this leads to a contradiction.

Consider the first two letters $q_n$ and $q_n^{-1}$ to be deleted in (39′). As in part **II**, which we have already proved, we have three possible cases.

    1°.  The first-deleted letters $q_n$ and $q_n^{-1}$ both appear in the word $Q_1$.

    2°.  One of them appears in $Q_1$, and the other is engendered in the sequence (39′) as a result of an application of schema (3′).

    3°.  Both letters are engendered in the sequence (39′).

**Case** 1°. If the considered letters $q_n$ and $q_n^{-1}$ both appear in the same segment $V_i$ of the word $Q_1$, then repeating the reasoning of case 1° of part **II**, we reach the conclusion that there exists a word of type $Q_1$, containing $m_0 - 2$ letters $q_n^\sigma$, and which is transformed into the empty word by a sequence of transformations with index $\lambda \leq \lambda_0$. This contradicts the assumption that for all words $Q_1$ with $m < m_0$ part **II** is proved for $\lambda \leq \lambda_0$. Consequently, we find that the letters $q_n$ and $q_n^{-1}$ are contained in different segment $V_i$ of the word $Q_1$.

We now have two possibilities:

    1)  the considered deletion is right-sided,

    2)  the considered deletion is left-sided.

**Case 1).** Suppose the letter $q_n^{-1}$ appears in the segment $V_i$. Then $q_n$ appears in $V_{i+\ell}$, where $\ell \geq 1$. The word $Q_1$ has the form

$$Q_1 \equiv Q_1' a_n \underline{q}_n^{-1} B A^{-1} V_i' a_n^{-1} V_{i+1} a_n^{-1} \cdots V_{i+\ell-1} a_n^{-1} V_{i+\ell}^{-1} A B^{-1} \underline{q}_n a_n^{-1} Q_1'',$$

where the underlined letters $q_n^\sigma$ are the first to delete in (39′), and $Q_1', Q_1''$ are some words of type $Q_1$. Replace (39′) by a canonical sequence for it with respect to $\underline{q}_n$:

$$Q_1 \equiv Q_1' a_n \underline{q}_n^{-1} B A^{-1} V_i' a_n^{-1} V_{i+1} a_n^{-1} \cdots V_{i+\ell-1} a_n^{-1} V_{i+\ell}' A B^{-1} \underline{q}_n a_n^{-1} Q_1'' \to \cdots$$

$$\cdots \to Q_1' a_n \underline{q}_n^{-1} B A^{-1} V_i' a_n^{-1} V_{i+1} a_n^{-1} \cdots V_{i+\ell-1} a_n^{-1} V_{i+\ell}' A B^{-1} U' \underline{q}_n V' a_n^{-1} Q_1'' \to \cdots$$

$$\cdots \to Z_1 \underline{q}_n^{-1} \underline{q}_n Z_2 \to Z_1 Z_2 \to \cdots \to 1. \tag{39''}$$

As the underlined letters $\underline{q}_n$ and $\underline{q}_n^{-1}$ are the first to be deleted, the words

$$V_i', V_{i+1}, \ldots, V_{i+\ell-1}, V_{i+\ell}'$$

do not contain any letters $q_n^\sigma$, i.e.

$$V_i' = a_{n-1}^{r_i}, \quad V_{i+1} = a_{n-1}^{r_{i+1}}, \ldots, V_{i+\ell-1} = a_{n-1}^{r_{i+\ell-1}},$$

$$V_{i+\ell}' = a_{n-1}^{r_{i+\ell}}.$$

In the sequence (39″) we now have the following independent transitions without insertions of letters $q_n^\sigma$:

$$B A^{-1} a_{n-1}^{r_i} a_n^{-1} a_{n-1}^{r_{i+1}} a_n^{-1} \cdots a_{n-1}^{r_{i+\ell-1}} a_n^{-1} a_{n-1}^{r_{i+\ell}} A B^{-1} U' \to \cdots \to 1. \tag{40}$$

In the sequence (40), we only have transformations in the group $\mathfrak{A}_{q_{n-1},A,B}$. Consequently, by Corollary 1 of Lemma 2, §1, Chapter II, we have

$$(a_{n-1}^{r_i} a_n^{-1} a_{n-1}^{r_{i+1}} a_n^{-1} \cdots a_{n-1}^{r_{i+\ell-1}} a_n^{-1} a_{n-1}^{r_{i+\ell}})_{\overline{S}} = a_n^{-\ell} = 1 \quad \text{in } \mathfrak{A}_P,$$

which is impossible, as $\ell \geq 1$.

Analogously, we prove the impossibility of **Case 2)**, with the only different being that in the place of Corollary 1 of Lemma 2, §1, Chapter II, we use Lemma 1, §2, Chapter I together with Lemma 1, §1, Chapter II.

**Case** $2°$. One of the first-deleted letters $q_n^\sigma$ in (39′) appears in the word $Q_1$, and the other is engendered in the sequence (39′) from some letter $q_n$, appearing in $Q_1$, as in the sequence (39′) there are no insertions of letters $q_n^\sigma$. If this letter $q_n$ appears in the same segment $V_i$ of $Q_1$ as the first of the letters $\underline{q}_n^\sigma$, then, by repeating all arguments carried out in case $2°$ of part **II**, this leads to a contradiction. Thus the only case remaining is when one of the first-deleted letters $\underline{q}_n^\sigma$ appears in $V_{i_1}$, and the other is engendered from some letter $q_n$, contained in $V_{i_2}$, where $i_1 \neq i_2$. Clearly, $i_2 < i_1$, as otherwise the engendering letter $q_n$ would lie between the two letters deleted before it is deleted itself. Now, we have two cases:

1) the considered deletion is right-sided.
2) the considered deletion is left-sided.

**Case 1).** The letter $\underline{q}_n$ appears in $V_{i_1}$, and the letter $\underline{q}_n^{-1}$ is engendered from a letter $q_n$, appearing in the segment $V_{i_2}$. The word $Q_1$ then has the form

$$Q_1' A B^{-1} q_n a_n^{-1} V_{i_2}' a_n^{-1} V_{i_2+1} a_n^{-1} \cdots a_n^{-1} V_{i_1}' A B^{-1} \underline{q}_n a_n^{-1} Q_1''.$$

Consider a canonical sequence for (39′) with respect to the individuality $\underline{q}_n$,

$$Q_1' A B^{-1} q_n a_n^{-1} V_{i_2}' a_n^{-1} V_{i_2+1} a_n^{-1} \cdots a_n^{-1} V_{i_1}' A B^{-1} \underline{q}_n a_n^{-1} Q_1'' \to \cdots$$
$$\cdots \to Q_1' A B^{-1} q_n a_n^{-1} V_{i_2}' a_n^{-1} V_{i_2+1} a_n^{-1} \cdots a_n^{-1} V_{i_1}' A B^{-1} U' \underline{q}_n V' a_n^{-1} Q'' \to \cdots \to 1.$$

The first part of this sequence, in which no transformations are performed on the underlined letter $\underline{q}_n$, save deletions, is replaced by a canonical sequence with respect to the hereditary chain engendered by the letter $\underline{q}_n^{-1}$:

$$Q_1' A B^{-1} q_n a_n^{-1} V_{i_2}' a_n^{-1} V_{i_2+1} a_n^{-1} \cdots a_n^{-1} V_{i_1}' A B^{-1} U' \underline{q}_n V' a_n Q'' \to \cdots$$
$$\cdots \to Q_1 A B^{-1} U_{1,0} q_n V_{1,0}' A B^{-1} U_{2,0} q_n V_{2,0}' \cdots$$
$$\cdots U_{k,0} q_n V_{k,0}' a_n \underline{q}_n^{-1} B A^{-1} V_{k,0}'' a_n^{-1} V_{k-1,0}'' \cdots$$
$$\cdots V_{2,0}'' a_n^{-1} V_{1,0}'' a_n^{-1} V_{i_2}' a_n^{-1} V_{i_2+1} a_n^{-1} \cdots a_n^{-1} V_{i_1}' A B^{-1} U' \underline{q}_n V' a_n^{-1} Q_1'' \to \cdots$$
$$\cdots Z_1 \underline{q}_n^{-1} \underline{q}_n Z_2 \to Z_1 Z_2 \to \cdots \to 1.$$

In this sequence, from the moment they are underlined, no transformations, except deletions, are performed on the letters $\underline{q}_n$ and $\underline{q}_n^{-1}$. We have the independent transition

$$B A^{-1} V_{k,0}'' a_n^{-1} V_{k-1,0}'' \cdots V_{2,0}'' a_n^{-1} V_{1,0}'' a_n^{-1} V_{i_2}' a_n^{-1} V_{i_2+1} a_n^{-1} \cdots$$
$$\cdots a_n^{-1} V_{i_1}' A B^{-1} U' \to \cdots \to 1, \qquad (41)$$

where there are no deletions of letters $q_n^\sigma$. Consequently,

$$V_{t,0}'' = a_{n-1}^{s_t} \ (t = 1, 2, \ldots, k), \quad V_{i_2} = a_{n-1}^{r_{i_2}}, \quad V_{i_2+1} = a_{n-1}^{r_{i_2+1}} \quad \text{etc.}$$
$$V_{i_1-1} = a_{n-1}^{r_{i_1-1}}, \quad V_{i_1}' = a_{n-1}^{r_{i_1}}.$$

The sequence (41) contains only transformations in the group $\mathfrak{A}_{q_{n-1},A,B}$. By Corollary 1 of Lemma 2, §1, Chapter II, we have

$$(a_{n-1}^{s_k} a_n^{-1} a_{n-1}^{s_{k-1}} \cdots a_{n-1}^{s_2} a_n^{-1} a_{n-1}^{s_1} a_n^{-1} a_n^{r_{i_2}} a_{n-1}^{-1} a_n^{r_{i_2+1}} a_{n-1}^{-1} \cdots a_n^{-1} a_{n-1}^{r_{i_1}})_{\overline{S}} =$$
$$= a_n^{k+i_1-i_2} = 1 \quad \text{in } \mathfrak{A}_P,$$

which is impossible, as $k \geq 0$ and $i_1 > i_2$.

**Case 2)** is proved analogously, with the only difference that in place of Corollary 1 of Lemma 2, §1, Chapter II, we use Lemma 1, §2, Chapter I.

**Case 3°.** The letters $\underline{q}_n$ and $\underline{q}_n^{-1}$ are both engendered in the sequence (39′). Then, as we have seen earlier, they must be generated from the same letter $q_n$, appearing in the word $Q_1$, as otherwise the considered deletion cannot be the first to occur in the sequence (39′). Repeating all the earlier arguments used in part **II** in case 3°, we either immediately obtain a contradiction, or else find a new sequence of transformations, which changes the same word $Q_1$ into the empty word and has index $\lambda = \lambda_0 - 2$, contradicting inductive hypothesis. This completes the argument for part **III**.

**Part IV.** Suppose we have a sequence of transformations with index $\lambda_0$:

$$Q_2 \equiv U_1 q_n G_1 A B^{-1} U_2 q_n G_2 A B^{-1} \cdots U_k q_n G_k A B^{-1} U_{k+1} \to \cdots \to 1, \qquad (42)$$

where $G_i$ is a word of the form $V_{i,1} a_n V_{i,2} a_n \cdots V_{i,t_i} a_n V_{i,t_i+1}$, and $k \geq 1$.

By virtue of the fact that part **I** has already been proved for $\lambda = \lambda_0$, we find a new sequence with index $\lambda \leq \lambda_0$:

$$Q_2 \equiv U_1 q_n G_1 A B^{-1} U_2 q_n G_2 A B^{-1} \cdots U_k q_n G_k A B^{-1} U_{k+1} \to \cdots \to 1, \qquad (42')$$

which does not contain insertions of letters $q_n^\sigma$.

By assumption, there cannot exist any sequence of type (42′) with index $\lambda < \lambda_0$. We will prove that it is also impossible for one to exist with index $\lambda_0$.

We may assume that in the sequence (42′), no transformations of type (3′) are performed on the distinguished letters $q_n$ in $Q_2$. For suppose that some such transformation was performed on the first of these letters. Then we can consider a canonical sequence for (42′) with respect to this letter $q_n$,

$$U_1 q_n G_1 A B^{-1} U_2 q_n G_2 A B^{-1} \cdots U_k q_n G_k A B^{-1} U_{k+1} \to \cdots$$
$$\cdots \to U_1 U_0 q_n U_0 G_1 A B^{-1} U_2 q_n G_2 A B^{-1} \cdots U_k q_n G_k A B^{-1} U_{k+1} \to \cdots \to 1. \qquad (42'')$$

The second part of the sequence (42″) has index $\lambda < \lambda_0$ and transforms a word of type $Q_2$ into the empty word, which is impossible by assumption.

As $k \geq 1$, the minimal number $m$ of letters $q_n^\sigma$, contained in the word $Q_2$, is one. Suppose the word in the sequence (42′) contains only a single letter $q_n$. Then (42′) has the form

$$U_1 q_n G_1 A B^{-1} U_2 \to \cdots \to 1, \qquad (43)$$

where $G_1$ does not contain any letters $q_n^\sigma$.

As, by assumption, no transformations of type (3′) are performed on the distinguished letter $q_n$, and as there are no insertions of letters $q_n^\sigma$, in the sequence (43) there will not appear any letters $q_n^{-1}$, and the letter $q_n$, contained in the initial word, cannot be deleted. Hence, for $m = 1$ the sequence (42′) cannot exist.

Suppose that for all words of type $Q_2$ with $m < m_0$ letters $q_n^\sigma$, we have proved the impossibility of the existence of a sequence of type (42′) with index $\lambda_0$. We will now prove the same for all words $Q_2$ with $m = m_0$ letters $q_n^\sigma$.

Consider the individualities of $q_n$ and $q_n^{-1}$, which in (42′) delete first. As in **III**, we have three possibilities:

1°.  The deletion of the first letters $\underline{q}_n$ and $\underline{q}_n^{-1}$ both appear in the word $Q_2$.

2°.  One of them appears in $Q_2$, and the other is engendered in the sequence (42′).

3°.  Both are engendered in the sequence (42′).

**Case 1°.** We have two subcases:

1)  the letter $\underline{q}_n$ coincides with one of the letters $q_n$ distinguished in $Q_2$;

2)  the letters $\underline{q}_n$ appears in one of the segments of type $G_i$.

**Case 1).** This case, in turn, divides into two subcases:

α) the considered deletion is right-sided,

β) the considered deletion is left-sided.

**Subcase α).** Suppose the letter $\underline{q}_n^{-1}$ appears in the segment $G_i$, which has the form

$$G_i' a_n \underline{q}_n^{-1} B A^{-1} V_{i,1} a_n V_{i,2} a_n \cdots a_n V_{i,\nu}.$$

The sequence (42′) can then be written as

$$Q_2 \equiv Q_2' q_n G_i' a_n \underline{q}_n^{-1} B A^{-1} V_{i,1} a_n V_{i,2} a_n \cdots a_n V_{i,\nu} A B^{-1} U_{i+1} \underline{q}_n G_{i+1} A B^{-1} Q_2'' \to \cdots$$
$$\cdots \to Z_1 \underline{q}_n^{-1} \underline{q}_n Z_2 \to Z_1 Z_2 \to \cdots \to 1.$$

Here no transformations are performed on the underlined letters $\underline{q}_n^{\sigma}$ until they are deleted. Hence we have the independent transitions:

$$\begin{array}{llll}
1. & Q_2' q_n G_i' a_n \to \cdots \to Z_1 & \text{with index } \lambda_1, \\
2. & G_{i+1} A B^{-1} Q_2'' \to \cdots \to Z_2 & \text{»} \quad \text{»} \quad \lambda_2, \\
3. & Z_1 Z_2 \to \cdots \to 1 & \text{»} \quad \text{»} \quad \lambda_3,
\end{array}$$

where $\lambda_1 + \lambda_2 + \lambda_3 = \lambda_0$. By combining 1, 2, and 3, we find a sequence of transformations with index $\lambda_0$:

$$Q_2' q_n G_i' a_n G_{i+1} A B^{-1} Q_2'' \to \cdots \to Z_1 G_{i+1} A B^{-1} Q_2'' \to \cdots \to Z_1 Z_2 \to \cdots \to 1,$$

which transforms into 1 a word of type $Q_2$, containing $m_0 - 2$ letters $q_n^{\sigma}$, which is impossible by assumption.

**Subcase β).** As the considered deletion is left-sided, the segment $G_i$, in which we find the considered letter $\underline{q}_n^{-1}$, is located to the right of the letter $\underline{q}_n$, and hence the word $Q_2$ has the form

$$Q_2 \equiv Q_2' \underline{q}_n V_{i,1} a_n V_{i,2} a_n \cdots a_n V_{i,\nu} a_n \underline{q}_n^{-1} B A^{-1} G_i A B^{-1} Q_2'',$$

where $\nu \geq 0$.

As, by assumption, no transformations are performed on the letters $q_n$ in $Q_2$, except deletions, it follows by (42′) that we have a transition without deletions of letters $q_n^{\sigma}$:

$$V_{i,1} a_n V_{i,2} a_n \cdots a_n V_{i,\nu} a_n \to \cdots \to 1. \tag{44}$$

Then the words $V_{i,j}$ cannot contain letters $q_n^{\sigma}$, i.e. each $V_{i,j}$ is equal to $a_{n-1}^{r_j}$. The transition (44) is a sequence of transformations in the group $\mathfrak{A}_{q_{n-1}, A, B}$. By Lemma 1, §1, Chapter II, we have

$$(a_{n-1}^{r_1} a_n a_{n-1}^{r_2} a_n \cdots a_n a_{n-1}^{r_\nu} a_n) = 1 \quad \text{in } \mathfrak{A}_P,$$

whence by Lemma 1, §2, Chapter I it follows that

$$(a_{n-1}^{r_1} a_n a_{n-1}^{r_2} a_n \cdots a_n a_{n-1}^{r_\nu} a_n)_{\overline{S}} = a_n^{\nu+1} = 1 \quad \text{in } \mathfrak{A}_P,$$

which is impossible, as $\nu + 1 \geq 1$.

**Case 2).** Suppose the letter $\underline{q}_n$ appears in the segment $G_{i_1}$, and that $\underline{q}_n^{-1}$ appears in $G_{i_2}$. Clearly, $i_1 = i_2$, as otherwise between the letters in $Q_2$ there will be some other letter $q_n$, and the considered letters $\underline{q}_n$ and $\underline{q}_n^{-1}$ cannot be the first to be deleted. Thus, the letters $\underline{q}_n$ and $\underline{q}_n^{-1}$ must lie in the same segment $G_{i_1}$. If both now appear in the same segment $V_{i_1,j}$, then repeating the reasoning carried out for part **II** (case $1°$), we find a sequence of type $(42')$ with index $\lambda \leq \lambda_0$ for a word $Q_{2,0}$ with $m = m_0 - 2$. This contradicts the inductive hypothesis. Consequently, we may assume that $\underline{q}_n$ and $\underline{q}_n^{-1}$ appear in different segments $V_{i_1,j_1}$ and $V_{i_1,j_2}$ in $G_{i_1}$.

We have two subcases:

$\alpha$) the considered deletion is right-sided,

$\beta$) the considered deletion is left-sided.

**Subcase $\alpha$).** As the deletion is right-sided, $j_1 > j_2$. The word $G_{i_2}$ then has the form

$$G'_{i_1} a_n \underline{q}_n^{-1} B A^{-1} V'_{i_1,j_2} a_n \cdots a_n V'_{i_1,j_1} A B^{-1} \underline{q}_n a_n^{-1} G''_{i_1}.$$

The sequence $(42')$ can then be written as:

$$Q_2 \equiv Q'_2 q_n G'_{i_1} a_n \underline{q}_n^{-1} B A^{-1} V'_{i_1,j_2} a_n V_{i_1,j_2+1} \cdots$$
$$\cdots a_n V_{i_1,j_1} A B^{-1} \underline{q}_n a_n^{-1} G''_{i_1} A B^{-1} Q''_2 \to \cdots$$
$$\cdots \to Z_1 \underline{q}_n^{-1} \underline{q}_n Z_2 \to \cdots \to Z_1 Z_2 \to \cdots \to 1, \qquad (45)$$

where the underlined letters $\underline{q}_n^\sigma$ are the first to be deleted. Hence we find

$$V_{i_2,j_2} = a_{n-1}^{r_{j_2}}, \ V_{i_2,t} = a_{n-1}^{r_t} \ (i = j_2+1, \ldots, j_1-1), \ V'_{i_1,j_1} = a_{n-1}^{r_{j_1}}.$$

We replace the sequence $(45)$ with a canonical sequence with respect to the letter $\underline{q}_n$:

$$Q'_2 q_n G'_{i_1} a_n \underline{q}_n^{-1} B A^{-1} a_{n-1}^{r_{j_2}} a_n a_{n-1}^{r_{j_2+1}} a_n \cdots a_n a_{n-1}^{r_{j_1}} A B^{-1} \underline{q}_n a_n^{-1} G''_{i_1} A B^{-1} Q''_2 \to \cdots$$
$$\cdots \to Q'_2 q_n G'_{i_1} a_n \underline{q}_n^{-1} B A^{-1} a_{n-1}^{r_{j_2}} a_n a_{n-1}^{r_{j_2+1}} a_n \cdots \qquad (45')$$
$$\ldots a_n a_{n-1}^{r_{j_1}} A B^{-1} U \underline{q}_n U a_n^{-1} G''_{i_1} A B^{-1} Q''_2 \to \ldots \to Z_1 \underline{q}_n^{-1} \underline{q}_n Z_2 \to Z_1 Z_2 \to \ldots \to 1.$$

From the sequence $(45')$, we find a transition without insertions of letters $q_n^\sigma$:

$$B A^{-1} a_{n-1}^{r_{j_2}} a_n a_{n-1}^{r_{j_2+1}} a_n \cdots a_n a_{n-1}^{r_{j_1}} A B^{-1} U \to \cdots \to 1,$$

whence by virtue of Corollary 1 of Lemma 2, §1, Chapter II, we have

$$(a_{n-1}^{r_{j_2}} a_n a_{n-1}^{r_{j_2+1}} a_n \cdots a_n a_{n-1}^{r_{j_1}})_{\overline{S}} = a_n^{j_1-j_2} = 1 \quad \text{in } \mathfrak{A}_P,$$

which is impossible, as $j_1 > j_2$.

**Subcase $\beta$)** is dealt with analogously, with the only difference being that in place of Corollary 1 of Lemma 2, §1, Chapter II, we use Lemma 1, §2, Chapter I.

**Case 2°.** One of the first-deleted letters $\underline{q}_n$ and $\underline{q}_n^{-1}$ in the sequence $(42')$ appears in the word $Q_2$, and the other is engendered in $(42')$.

We have two possible cases:

1) the letter $\underline{q}_n$ appears in $Q_2$, and $\underline{q}_n^{-1}$ is engendered in the sequence $(42')$;

2) the letter $\underline{q}_n^{-1}$ appears in $Q_2$, and $\underline{q}_n$ is engendered in $(42')$.

**Case 1).** In this case, we have in turn two subcases:

$\alpha$) the letter $\underline{q}_n$ is one of the letters $q_n$ distinguished in $Q_2$;

$\beta$) the letter $\underline{q}_n$ appears in some $G_i$ in the word $Q_2$.

**Subcase $\alpha$).** As, by assumption, no transformations are performed on the letters $q_n$ inside $Q_2$, except deletions, the letter $\underline{q}_n^{-1}$ must hence be engendered by some letter $q_n^{(1)}$, which appears in some segment $G_{j_1}$. The segment $G_{j_1}$ must lie to the left of the letter $\underline{q}_n$ in the word $Q_2$. The word $Q_2$ now has the form

$$Q_2' q_n G_{j_1}' AB^{-1} q_n^{(1)} a_n^{-1} V_{j_1,i_1} a_n V_{j_1,i_1+1} a_n \cdots a_n V_{j_1,i_1+\nu} AB^{-1} U_{j_1+1} \underline{q}_n G_{j_1+1} AB^{-1} Q_2''.$$

The sequence (42′) is replaced by a canonical sequence with respect to the hereditary chain engendered from the letter $\underline{q}_n^{-1}$. We find a sequence with index $\lambda_0$:

$$Q_2' q_n G_{j_1}' AB^{-1} q_n^{(1)} a_n^{-1} V_{j_1,i_1} a_n V_{j_1,i_1+1} a_n \cdots$$
$$\cdots a_n V_{j_1,i_1+\nu} AB^{-1} U_{j_1+1} \underline{q}_n G_{j_1+1} AB^{-1} Q_2'' \to \cdots$$
$$\cdots \to Q_2' q_n G_{j_1}' AB^{-1} U_1 q_n^{(1)} V_1' AB^{-1} U_2 q_n^{(2)} V_2' \cdots$$
$$\cdots U_k q_n^{(k)} V_k' a_n \underline{q}_n^{-1} BA^{-1} V_k'' a_n^{-1} V_{k-1}'' a_n^{-1} \cdots a_n^{-1} V_1'' a_n^{-1} V_{j_1,i_1}' a_n V_{j_1,i_1+1} a_n \cdots$$
$$\cdots a_n V_{j_1,i_1+\nu} AB^{-1} U_{j_1+1} \underline{q}_n G_{j_1+1} AB^{-1} Q_2'' \to \cdots$$
$$\cdots \to Z_1 \underline{q}_n^{-1} \underline{q}_n Z_2 \to Z_1 Z_2 \to \cdots \to 1. \tag{46}$$

In the second part of the sequence (46), no transformations are performed on the letters $\underline{q}_n^{\sigma}$ until their deletion. Hence, from (46) we find the independent transitions:

1. $Q_2' q_n G_{j_1} AB^{-1} U_1 q_n^{(1)} V_1' AB^{-1} U_2 q_n^{(2)} V_2' ... U_k q_n^{(k)} V_k' a_n \to ... \to Z_1$    with index $\lambda_1$,
2.                                      $G_{j_1+1} AB^{-1} Q_2'' \to \cdots \to Z_2$    »       »    $\lambda_2$,
3.                                             $Z_1 Z_2 \to \cdots \to 1$    »       »    $\lambda_3$,

We have $\lambda_1 + \lambda_2 + \lambda_3 < \lambda_0$, as the first part of the sequence (46), in which the letter $\underline{q}_n^{-1}$ is engendered, has index $\lambda' > 0$.

Combining the transitions **1**, **2**, and **3**, we find a sequence of elementary transformations:

$$Q_2' q_n G_{j_1} AB^{-1} U_1 q_n^{(1)} V_1' AB^{-1} U_2 q_n^{(2)} V_2' \cdots U_k q_n^{(k)} V_k' a_n G_{j_1+1} AB^{-1} Q_2'' \to \cdots$$
$$\cdots \to Z_1 G_{j_1+1} AB^{-1} Q_2'' \to \cdots \to Z_1 Z_2 \to \cdots \to 1, \tag{46′}$$

which transforms into 1 some word of type $Q_2$, and has index $\lambda = \lambda_1 + \lambda_2 + \lambda_3 < \lambda_0$. This contradicts the inductive hypothesis, i.e. that part **IV** of the lemma has been proved for $\lambda < \lambda_0$.

**Subcase $\beta$).** Suppose the letter $\underline{q}_n$ appears in the segment $G_i$. Then the letter $\underline{q}_n^{-1}$ is engendered by some letter $q_n^{(1)}$, which lies in the same segment $G_i$, as otherwise the considered deletion cannot be the first to occur in (42′). The word $G_i$ has the form

$$G_i' AB^{-1} q_n^{(1)} a_n^{-1} V_{i,j_1} a_n V_{i_1,j_1+1} a_n \cdots V_{i,j_1+\nu} a_n V_{i,j_1+\nu+1} AB^{-1} \underline{q}_n a_n^{-1} G_i''.$$

The sequence (42′) then becomes:

$$Q_2 \equiv Q_2' q_n G_i' AB^{-1} q_n^{(1)} a_n^{-1} V_{i,j_1} a_n V_{i_1,j_1+1} a_n \cdots$$
$$\cdots V_{i,j_1+\nu} a_n V_{i,j_1+\nu+1} AB^{-1} \underline{q}_n a_n^{-1} G_i'' AB^{-1} Q_2'' \to \cdots \to 1.$$

Construct a canonical sequence with respect to the letter $\underline{q}_n$:

$$Q_2' q_n G_i' A B^{-1} q_n^{(1)} a_n^{-1} V_{i,j_1} a_n V_{i,j_1+1} a_n \cdots$$
$$\cdots V_{i,j_1+\nu} a_n V_{i,j_1+\nu+1} A B^{-1} \underline{q}_n a_n^{-1} G_i'' A B^{-1} Q_2'' \to \cdots$$
$$\cdots \to Q_2' q_n G_i' A B^{-1} q_n^{(1)} a_n^{-1} V_{i,j_1} a_n V_{i,j_1+1} a_n \cdots$$
$$\cdots V_{i,j_1+\nu} a_n V_{i,j_1+\nu+1} A B^{-1} U \underline{q}_n V a_n^{-1} G_i'' A B^{-1} Q_2'' \to \cdots \to 1. \quad (47)$$

In the second part of this sequence, no transformations except deletions is performed on the letter $\underline{q}_n$; for this part, construct a canonical sequence with respect to the hereditary chain engendered from the letter $\underline{q}_n^{-1}$:

$$Q_2' q_n G_i' A B^{-1} q_n^{(1)} a_n^{-1} V_{i,j_1} a_n \cdots$$
$$\cdots V_{i,j_1+\nu} a_n V_{i,j_1+\nu+1} A B^{-1} U \underline{q}_n V a_n^{-1} G_i'' A B^{-1} Q_2'' \to \cdots$$
$$\cdots \to Q_2' q_n G_i' A B^{-1} U_1 q_n^{(1)} V_1' A B^{-1} U_2 q_n^{(2)} V_2' \cdots$$
$$\cdots V_k q_n^{(k)} V_k' a_n \underline{q}_n^{-1} B A^{-1} V_k'' a_n^{-1} V_k'' a_n^{-1} V_{k-1}'' a_n^{-1} \cdots a_n^{-1} V_1'' a_n^{-1} V_{i,j_1+1} a_n \cdots$$
$$\cdots V_{i,j_1+\nu} a_n V_{i,j_1+\nu+1} A B^{-1} U \underline{q}_n V a_n^{-1} G_i'' A B^{-1} Q_2'' \to \cdots$$
$$\cdots \to Z_1 \underline{q}_n^{-1} \underline{q}_n Z_2 \to Z_1 Z_2 \to \cdots \to 1. \quad (47')$$

In $(47')$, no transformations are performed on the letters $\underline{q}_n^{\sigma}$, except deletions.

From $(47)$ and $(47')$ we find the independent transitions

1. $\underline{q}_n \to \cdots \to U \underline{q}_n V$ with index $\lambda_1$.
2. $q_n^{(1)} \to \cdots \to U_1 q_n^{(1)} V_1' A B^{-1} U_2 q_n V_2' \ldots U_k q_n V_k' a_n q_n^{-1} B A^{-1} V_k'' a_n^{-1} V_{k-1}'' a_n^{-1} \cdots$
   $\cdots a_n^{-1} V_1''$ with index $\lambda_2 > 0$,
3. $Q_2' q_n G_i' A B^{-1} U_1 q_n^{(1)} V_1' A B^{-1} U_2 q_n V_2' \cdots U_k q_n V_k' a_n \to \cdots \to Z_1$ with index $\lambda_3$,
4. $V a_n^{-1} G_i'' A B^{-1} Q_2 \to \cdots \to Z_2$ with index $\lambda_4$,
5. $B A^{-1} V_k'' a_n^{-1} V_{k-1}'' a_n^{-1} \cdots a_n^{-1} V_1'' a_n^{-1} V_{i,j_1} a_n V_{i,j_1+1} a_n \cdots$
   $\cdots V_{i,j_1+\nu} a_n V_{i,j_1+\nu+1} A B^{-1} U \to \cdots \to 1$ with index $\lambda_5$,
6. $Z_1 Z_2 \to \cdots \to 1$ with index $\lambda_6$, where $\lambda_1 + \lambda_2 + \lambda_3 + \lambda_4 + \lambda_5 + \lambda_6 = \lambda_0$.

Transition 5 does not contain deletions of letters $q_n^{\sigma}$, as the considered deletion occurs first in the sequence $(47')$. Consequently, $\lambda_5 = 0$, and the words $V_t''$ and $V_{i,t}$ do not contain letters $q_n^{\sigma}$, i.e.

$$V_t'' = a_{n-1}^{s_t} \quad (t = 1, 2, \ldots, k), \qquad V_{i,t} = a_{n-1}^{r_t} \quad (t = j_1, j_1 + 1, \ldots, j_1 + \nu + 1).$$

Transition 5 is a sequence of transformations in the group $\mathfrak{A}_{q_{n-1},A,B}$. By Lemma 2, §1, Chapter II, it follows that the sequence of exponents corresponding to the word $U$ is degenerate. Hence, $V = 1$ in the free group, i.e. we have a sequence of transformations consisting only of insertions:

$$1 \to \cdots \to a_n V a_n^{-1}. \quad (48)$$

Clearly, the index $\lambda'$ of the sequence $(48)$ is no greater than the index $\lambda_1$ of transition 1. Combining the sequences $(48)$, 3, 4, and 6, we find a sequence of transformations

$$Q_{2,0} \equiv Q_2' q_n G_i' A B^{-1} U_1 q_n V_1' A B^{-1} U_2 q_n V_2' \cdots U_k q_n V_k' G_i'' A B^{-1} Q_2'' \to \cdots$$
$$\cdots \to Q_2' q_n G_i' A B^{-1} U_1 q_n V_1' A B^{-1} U_2 q_n V_2' \cdots U_k q_n V_k' a_n V a_n^{-1} G_i'' A B^{-1} Q_2'' \to \cdots$$
$$\cdots \to Z_1 V a_n^{-1} G_i'' A B^{-1} Q_2'' \to \cdots \to Z_1 Z_2 \to \cdots \to 1, \quad (47'')$$

which has index

$$\lambda = \lambda' + \lambda_3 + \lambda_4 + \lambda_6 \le \lambda_1 + \lambda_3 + \lambda_4 + \lambda_6 < \lambda_1 + \lambda_3 + \lambda_4 + \lambda_6 + \lambda_2 = \lambda_0.$$

The sequence $(47'')$ transforms a word $Q_{2,0}$ of type $Q_2$ into the empty word and has index $\lambda < \lambda_0$, which contradicts the inductive hypothesis (of **IV**).

**Case 2).** The letter $\underline{q}_n^{-1}$ appears in $Q_2$, and $\underline{q}_n$ is engendered in the sequence $(42')$. Suppose $\underline{q}_n^{-1}$ appears in the segment $G_i$. Then the letter $q_n^{(1)}$ in the word $Q_2$, which engenders the letter $\underline{q}_n$, must appear in the same segment $G_i$, as otherwise the considered deletion cannot be the first one. In short, the word $G_i$ must have the form

$$G_i' AB^{-1} q_n^{(1)} a_n^{-1} a_{n-1}^{r_1} a_n a_{n-1}^{r_2} a_n \cdots a_{n-1}^{r_\nu} a_n a_{n-1}^{r_{\nu+1}} a_n \underline{q}_n^{-1} BA^{-1} G_i'',$$

where $q_n^{(1)}$ is the individuality which engenders in the sequence $(42')$ the letter $\underline{q}_n$ which deletes with $\underline{q}_n^{-1}$. Construct a canonical sequence with respect to the hereditary chain engendered by the letter $\underline{q}_n$:

$$Q_2 \equiv Q_2' q_n G_i AB^{-1} q_n^{(1)} a_n^{-1} a_{n-1}^{r_1} a_n a_{n-1}^{r_2} a_n \cdots$$
$$\cdots a_{n-1}^{r_\nu} a_n a_{n-1}^{r_{\nu+1}} a_n \underline{q}_n^{-1} BA^{-1} G_i'' AB^{-1} Q_2'' \to \cdots$$
$$\cdots \to Q_2' q_n G_i' AB^{-1} U_1 q_n^{(1)} V_1' AB^{-1} U_2 q_n V_2' \cdots V_{k-1}' AB^{-1} U_k \underline{q}_n V_k a_n^{-1} V_{k-1}'' a_n^{-1} \cdots$$
$$\cdots a_n^{-1} V_1'' a_n^{-1} a_{n-1}^{r_1} a_n a_{n-1}^{r_2} a_n \cdots a_{n-1}^{r_\nu} a_n a_{n-1}^{r_{\nu+1}} a_n \underline{q}_n^{-1} BA^{-1} G_i'' AB^{-1} Q_2'' \to \cdots$$
$$\cdots \to Z_1 \underline{q}_n \underline{q}_n^{-1} Z_2 \to Z_1 Z_2 \to \cdots \to 1. \qquad (49)$$

As the considered deletion is the first to occur, we have

$$V_j'' = a_{n-1}^{s_j} \quad (j = 1, 2, \dots, k-1) \text{ and } V_k = a_{n-1}^{s_k}.$$

Hence $U_k = q_n^{-s_k}$. From the sequence $(49)$, we find the independent transitions:

1. $q_n^{(1)} \to \cdots \to U_1 q_n^{(1)} V_1' AB^{-1} U_2 q_n V_2' \cdots$
   $\cdots V_{k-1} AB^{-1} U_k q_n V_k a_n^{-1} V_{k-1}'' a_n^{-1} \cdots a_n^{-1} V_1''$ with index $\lambda_1$,
2. $a_{n-1}^{s_k} a_n^{-1} a_{n-1}^{s_{k-1}} a_n^{-1} \cdots a_n^{-1} a_{n-1}^{s_1} a_n^{-1} a_{n-1}^{r_1} a_n a_{n-1}^{r_2} a_n \cdots$
   $\cdots a_{n-1}^{r_\nu} a_n a_{n-1}^{r_{\nu+1}} a_n \to \cdots \to 1$ with index $\lambda_2$,
3. $Q_2' q_n G_i' AB^{-1} U_1 q_n^{(1)} V_1' AB^{-1} U_2 q_n V_2' \cdots V_{k-1}' AB^{-1} U_k \to \cdots \to Z_1$
   with index $\lambda_3$,
4. $BA^{-1} G_i'' AB^{-1} Q_2'' \to \cdots \to Z_2$ with index $\lambda_4$,
5. $Z_1 Z_2 \to \cdots \to 1$ with index $\lambda_5$.

We have $\lambda_1 > 0$, as in the case under consideration the letter $\underline{q}_n$ does not appear in $Q_2$, and $\lambda_1 + \lambda_2 + \lambda_3 + \lambda_4 + \lambda_5 = \lambda_0$. As there are no deletions of letters $q_n^\sigma$ in the transition 2, this is a sequence of transformations in the group $\mathfrak{A}_{q_{n-1}, A, B}$. By Lemma 1, §1, Chapter II we find that in the group $\mathfrak{A}_P$ the equality

$$a_{n-1}^{s_k} a_n^{-1} a_{n-1}^{s_{k-1}} a_n^{-1} \cdots a_n^{-1} a_{n-1}^{s_1} a_n^{-1} a_{n-1}^{r_1} a_n a_{n-1}^{r_2} a_n \cdots a_{n-1}^{r_\nu} a_n a_{n-1}^{r_{\nu+1}} a_n = 1.$$

By Lemma 2, §2, Chapter I we have $s_k = 0$, i.e. $U_k = 1$ in the free group. This implies that there is a sequence consisting only of insertions:

$$1 \to \cdots \to AB^{-1} U_k BA^{-1}, \qquad (50)$$

the index of which is equal to zero, as the word $U_k$ does not contain any letters $q_n^\sigma$.

By combining the sequences (50), 3, 4, and 5, we find a sequence of transformations

$$Q_{2,0} \equiv Q_2' q_n G_i' AB^{-1} U_1 q_n V_1' AB^{-1} U_2 q_n V_2' \cdots V_{k-1}' G_i'' AB^{-1} Q_2'' \rightarrow \cdots$$

$$\cdots \rightarrow Q_2' q_n G_i' AB^{-1} U_1 q_n V_1' AB^{-1} U_2 q_n V_2' \cdots V_{k-1}' AB^{-1} U_k BA^{-1} G_i'' AB^{-1} Q_2'' \rightarrow \cdots$$

$$\cdots \rightarrow Z_1 BA^{-1} G_i'' AB^{-1} Q_2'' \rightarrow \cdots \rightarrow Z_1 Z_2 \rightarrow \cdots \rightarrow 1 \quad (49')$$

with index $\lambda = \lambda_3 + \lambda_4 \lambda_5 < \lambda_1 + \lambda_2 + \lambda_3 + \lambda_4 + \lambda_5 = \lambda_0$.

The sequence (49') transforms the word $Q_{2,0}$ of type $Q_2$ into the empty word, and has index $\lambda < \lambda_0$, which contradicts the inductive hypothesis.

**Case 3°** is dealt with analogously to case 3° of the proof of part **III** for $\lambda = \lambda_0$. This completes the proof of the Main Lemma. $\square$

**Lemma 4.** *Suppose $A = pXpX^+$ and $B = pYpY^+$, where $X$ and $Y$ are words written using letters from $S$. If $A = B$ in the group $\mathfrak{A}_P$, then the group $\mathfrak{A}_{q,A,B}$ is the trivial group. If $A \neq B$, then the group $\mathfrak{A}_P$ is a subgroup of the group $\mathfrak{A}_{q,A,B}$.*

*Proof.* Suppose $A = B$ in $\mathfrak{A}_P$. Clearly, it suffices to prove that in this case all letters of the positive alphabet

$$a_1, a_2, \ldots, a_n, q_1, q_2, \ldots, q_n$$

of the group $\mathfrak{A}_{q,A,B}$ are equal to the identity element in $\mathfrak{A}_{q,A,B}$. The system of defining relations of the group $\mathfrak{A}_{q,A,B}$ is obtained by combining the systems (3), (4), (5), and (6) of relations (see §1, Chapter II).

From $A = B$ in $\mathfrak{A}_P$, it clearly follows that $A = B$ in $\mathfrak{A}_{q,A,B}$. Then $AB^{-1} = 1$, and the relations (4) can be written as:

$$a_i q_i = q_i \qquad (i = 1, 3, 5, \ldots, n-2). \quad (51)$$

From (51) it follows that in the group $\mathfrak{A}_{q,A,B}$ all generators $a_i$ are equal to 1 for all odd $i < n$. As the system of defining relations of the group $\mathfrak{A}_{q,A,B}$ contains all defining relations of the group $\mathfrak{A}_P$, it follows that in $\mathfrak{A}_{q,A,B}$ we also have that $a_i$ is equal to 1 for even indices $i$ (see Remark 1', p. 25). Substituting in the relations (5) the identity for every $a_i$, we find that all $q_i$ are equal to the identity for $i = 1, 2, \ldots, n-1$. From the relations (6) we have $q_n = q_n q_n^{-1} = 1$. Consequently, the group $\mathfrak{A}_{q,A,B}$, when $A = B$ in $\mathfrak{A}_P$, will be the trivial group.

Suppose $A \neq B$ in $\mathfrak{A}_P$. We will prove that then $\mathfrak{A}_P$ is a subgroup of the group $\mathfrak{A}_{q,A,B}$. As is known, for this it suffices to show that for an arbitrary word $E$, consisting only of letters from the alphabet of the group $\mathfrak{A}_P$, we have $E = 1$ in $\mathfrak{A}_{q,A,B}$ if and only if $E = 1$ in $\mathfrak{A}_P$. Suppose we have a sequence of transformations with index $\lambda$:

$$E \rightarrow \cdots \rightarrow 1. \quad (52)$$

By the Main Lemma (part **I**) we then have a sequence of elementary transformations without insertions of letters $q_n^\sigma$:

$$E \rightarrow \cdots \rightarrow 1. \quad (52')$$

In the sequence (52') the letters $q_n^\sigma$ do not appear, as they did not appear in the word $E$. This implies that in (52'), there are no transformations from the schema (3'), i.e. (52') is a sequence of transformations in the group $\mathfrak{A}_{q_{n-1},A,B}$. As the word $E$ contains only letters from the alphabet of the group $\mathfrak{A}_P$, so by Lemma 1, §1, Chapter II we have $E = 1$ in $\mathfrak{A}_P$. Lemma 4 is proved. $\square$

### §3. Theorem on unsolvability of algorithmic problems in the theory of groups

Theorem 2′, §2, Chapter II claims that there does not exist an algorithm which decides for an arbitrary pair of words of the form $pXpX^{+}$ − where $p$ is the pivot letter of the group $\mathfrak{A}_P$, and $X$ is any word consisting of letters from the class $S$ − whether or not they are equal in the group $\mathfrak{A}_P$. In §1 of the present chapter we constructed a group $\mathfrak{A}_{q,A,B}$, for every pair of words $(A, B)$, given by a finite number of defining relations.

By Lemma 4, §2, Chapter II, it follows that there is no algorithm which decides for every group $\mathfrak{A}$, defined by a finite number of generators and defining relations, whether or not it is the trivial group, as such an algorithm cannot even exist for the system of groups $\mathfrak{A}_{q,A,B}$.

**Theorem 1.** *Let $F_0$ be an arbitrary group, given by a finite number of generators and defining relations. There does not exist an algorithm which decides for an arbitrary group $F$, given by a finite number of generators and defining relations, whether or not it is isomorphic to $F_0$.*

*Proof.* Denote by $\mathfrak{A}_{q,A,B}^{(0)}$ the free product of the group $F_0$ by a group $\mathfrak{A}_{q,A,B}$ from the system of groups $\mathfrak{A}_{q,A,B}$ constructed in §1, Chapter II. Every group $\mathfrak{A}_{q,A,B}^{(0)}$ is given by a finite number of generators and defining relations. We will prove that the group $\mathfrak{A}_{q,A,B}^{(0)}$ is isomorphic to the group $F_0$ if and only if $A = B$ in $\mathfrak{A}_P$.

If $A = B$ in $\mathfrak{A}_P$, then by Lemma 4, §2, Chapter II the group $\mathfrak{A}_{q,A,B}$ is the trivial group, and hence the group $\mathfrak{A}_{q,A,B}^{(0)}$ is isomorphic to $F_0$. If this group $\mathfrak{A}_{q,A,B}^{(0)}$ is isomorphic to the group $F_0$, then $\mathfrak{A}_{q,A,B}^{(0)}$ is the free product of the two groups $\mathfrak{A}_{q,A,B}$ and $F_0$, of which neither is trivial. The group $F_0$, being isomorphic to the group $\mathfrak{A}_{q,A,B}^{(0)}$, itself decomposes as a free product with two factors $F_1'$ and $F_1$, isomorphic to the groups $\mathfrak{A}_{q,A,B}$ and $F_0$, respectively. Then the group $\mathfrak{A}_{q,A,B}^{(0)}$ decomposes as a free product with three factors: $\mathfrak{A}_{q,A,B}, F_1'$, and $F_1$, etc.

Continuing this reasoning, we find that if $A \neq B$ in $\mathfrak{A}_P$ and the groups $\mathfrak{A}_{q,A,B}^{(0)}$ and $F_0$ are isomorphic, then the group $\mathfrak{A}_{q,A,B}^{(0)}$ decomposes as a free product with arbitrarily many factors. But this is impossible, as it is known that the number of factors in a free product decomposition of a group with a finite number of generators is bounded above by the number of generators (see [3], p. 253). Consequently, when $A \neq B$ in $\mathfrak{A}_P$, the group $\mathfrak{A}_{q,A,B}^{(0)}$ cannot be isomorphic to the group $F_0$.

By Theorem 2, §2, Chapter I it now follows from what we just proved that there cannot be an algorithm which decides for every group $\mathfrak{A}_{q,A,B}^{(0)}$ whether or not it is isomorphic to $F_0$. Clearly, then, neither can such an algorithm exist for all groups given by a finite number of generators and defining relations. □

**Theorem 2.** *Let $\alpha_0$ be an arbitrary invariant property, and $\beta_0$ an arbitrary non-trivial hereditary property. If there exists a group $F_1$ given by a finite number of generators and defining relations and possessing the property $(\alpha_0 \,\&\, \beta_0)$, then the recognition problem for the property $(\alpha_0 \,\&\, \beta_0)$ is unsolvable.*

*Proof.* For the group $F_1$, which has the property $(\alpha_0 \,\&\, \beta_0)$, we will construct the system of groups $\mathfrak{A}_{q,A,B}^{(1)}$, just as we did above for the group $F_0$.

Theorem 2 will be proved if we prove that the property $(\alpha_0 \,\&\, \beta_0)$ is possessed by the group $\mathfrak{A}_{q,A,B}^{(1)}$ if and only if $A = B$ in $\mathfrak{A}_P$. If $A = B$ in $\mathfrak{A}_P$, then the group $\mathfrak{A}_{q,A,B}^{(1)}$ is isomorphic to the group $F_1$. Therefore, the invariant properties $\alpha_0$ and $\beta_0$, possessed

by the group $F_1$, will also be possessed by the group $\mathfrak{A}^{(1)}_{q,A,B}$, i.e. in $\mathfrak{A}^{(1)}_{q,A,B}$ the property $(\alpha_0 \ \& \ \beta_0)$ will also be possessed.

Suppose $A \neq B$ in $\mathfrak{A}_P$. Then by Lemma 4, §2, Chapter II the group $\mathfrak{A}_P$ will be a subgroup of the group $\mathfrak{A}_{q,A,B}$, and hence also of the group $\mathfrak{A}^{(1)}_{q,A,B}$. The group $\mathfrak{A}_P$ has a free subgroup on two generators (see the Corollary of Lemma 2, §2, Chapter I), which will hence also be a subgroup of the group $\mathfrak{A}^{(1)}_{q,A,B}$. In this case the non-trivial hereditary property $\beta_0$ cannot be possessed by the group $\mathfrak{A}^{(1)}_{q,A,B}$, as otherwise $\beta_0$, being hereditary, would be possessed by the free group on two generators, which is impossible, as by assumption $\beta_0$ is a non-trivial hereditary property. Consequently, if $A \neq B$ in $\mathfrak{A}_P$, the property $(\alpha_0 \ \& \ \beta_0)$ is not possessed by the group $\mathfrak{A}^{(1)}_{q,A,B}$. Theorem 2 is proved.   □

Note that the following properties are all non-trivial hereditary properties of groups:

(1) being a finite group,
(2) being a periodic group,
(3) being an abelian group,
(4) being a nilpotent group,
(5) being a solvable group, etc.

From Theorem 2 one can extract, for example, the following consequences:

(a) There does not exist an algorithm which decides for every group $F$, given by finitely many generators and defining relations, whether or not this group is a finite group of order $k$.

(b) There does not exist an algorithm which decides for every group $F$, given by finitely many generators and defining relations, whether or not this group is a nilpotent group of some fixed class $k$ or not, etc.

**Theorem 3.** *Let $\alpha_0$ be an arbitrary special (see p. 13) invariant property. If there exists a group $F_2$, given by a finite number of generators and defining relations, which has property $\alpha_0$, then the recognition problem for $\alpha_0$ is unsolvable.*

*Proof.* We construct, analogously to earlier, the system of groups $\mathfrak{A}^{(2)}_{q,A,B}$. The unsolvability of the recognition problem for $\alpha_0$ will be proved, if we prove that the property $\alpha_0$ is possessed by the group $\mathfrak{A}^{(2)}_{q,A,B}$ if and only if $A = B$ in $\mathfrak{A}_P$.

If $A = B$ in $\mathfrak{A}_P$, then the group $\mathfrak{A}^{(2)}_{q,A,B}$ is isomorphic to $F_2$, and the invariant property $\alpha_0$ will by assumption be possessed by this group.

Suppose $A \neq B$ in $\mathfrak{A}_P$. Then the group $\mathfrak{A}_P$ will, by Lemma 4, §2, Chapter II, be a subgroup of the group $\mathfrak{A}_{q,A,B}$, and hence also of the group $\mathfrak{A}^{(2)}_{q,A,B}$. Consequently, the group $\mathfrak{A}^{(2)}_{q,A,B}$ will have unsolvable identity problem, and hence cannot possess the special invariant property $\alpha_0$.   □

The *meta-identity problem* for groups is formulated in the following way: find an algorithm, which decides for every group $F$, with finitely many generators and defining relations, whether or not it has solvable identity problem. As the property of having solvable identity problem is a special invariant property, from Theorem 3 we find:

**Corollary.** *The meta-identity problem for groups is unsolvable.*

**Theorem 4.** *The recognition problem is unsolvable for the following properties of groups:*

1) *being a simple group,*
2) *having a free subgroup of rank $k$.*

*Proof.*

1) Construct the system of groups $\mathfrak{A}_{q,A,B}^{(0)}$, taking $F_0$ to be the simple group of order two. If $A = B$ in $\mathfrak{A}_P$, the group $\mathfrak{A}_{q,A,B}^{(0)}$ is isomorphic to $F_0$, i.e. is a simple group. If $A \neq B$ in $\mathfrak{A}_P$, then the group $\mathfrak{A}_{q,A,B}^{(0)}$ cannot be a simple group, as the normal subgroup in $\mathfrak{A}_{q,A,B}^{(0)}$ generated by $\mathfrak{A}_{q,A,B}$ cannot be trivial, and cannot coincide with the entire group $\mathfrak{A}_{q,A,B}^{(0)}$.

2) Consider the system of groups $\mathfrak{A}_{q,A,B}$. If $A = B$ in $\mathfrak{A}_P$, the group $\mathfrak{A}_{q,A,B}$ will be the trivial group and cannot contain a free subgroup. When $A \neq B$ in $\mathfrak{A}_P$, the group $\mathfrak{A}_{q,A,B}$ contains a free subgroup on two generators (Lemma 4, §2, Chapter II and the Corollary of Lemma 2, §2, Chapter I). But the free group on two generators contains a free subgroup on any finite number of generators. $\square$

S. I. ADIAN

# FINITELY PRESENTED GROUPS AND ALGORITHMS

*(Communicated by academician I. M. Vinogradov, 7 May 1957)*

A group $F$ will be said to be *finitely presented*, if it can be defined using a finite number of generators and defining relations. A property $\alpha$ of groups will be called *invariant* if it, when possessed by a group $F$, is also possessed by any group $F_1$ isomorphic to $F$.

In [3], the impossibility of an algorithm for the recognition of certain classes of invariant properties of finitely presented groups. A. A. Markov proved the impossibility of an algorithm for the recognition of a very broad class of properties of associative systems [2]. The first part of the present article is a continuation of the work in the article [3]. We will prove a theorem analogous to the theorem of A. A. Markov on associative systems.

**Theorem 1.** *Suppose $\alpha$ is an invariant property of groups. If there exists some finitely presented group $F_1$ which has property $\alpha$, and some finitely presented group $F_2$ which does not embed in any finitely presented group with this property, then there does not exist an algorithm which decides for every finitely presented group $F$ whether or not it has the property $\alpha$.*

In the article [1], P. S. Novikov constructed a finitely presented group with an unsolvable identity problem. Let this group, which we will denote by $F_0$, be given by the generators $a_1, a_2, \ldots, a_n$ and defining relations*

$$A_i = A_i'. \tag{1}$$

In the proof of Theorem 1, we will use the fact that the group $F_0$ is torsion-free (this is not proved in [1], but can without much difficulty be proved with [1] as a basis).

The following lemma is nearly obvious.

**Lemma 1.** *Any finitely presented group $F_2$ embeds isomorphically in a finitely presented group $F_2'$ defined by a system of generators in which every generator has infinite order.*

Let the group $F_2'$, in which the group $F_2$ from the statement of Theorem 1 is isomorphically embedded, be given by the generators $b_1, b_2, \ldots, b_m$ and the defining relations

$$B_j = B_j'. \tag{2}$$

Let $F$ denote the free product

$$F \equiv F_0 * F_2' * F_3,$$

where $F_3$ is the free group on a single generator $p$.

---

*In defining a finitely presented group we will only ever write out its positive alphabet and its non-trivial relations, i.e. those relations which do not hold in the free group.

The group $F$ will have the generators

$$a_1, a_2, \ldots, a_n, b_1, b_2, \ldots, b_m, p \qquad (3)$$

and defining relations (1) together with (2).

Two letters of the alphabet (3) will be said to be *homotypic*, if either both are generators of the group $F_0$, or if both are generators of the group $F_2'$. Otherwise, they will be said to be *heterotypic*.

**Lemma 2.** *Any two heterotypic letters in the alphabet of the group $F$ freely generate the subgroup of $F$ generated by them.*

**Lemma 3.** *The alphabet of the group $F$ can be arranged in a sequence with repetitions*

$$e_1, e_2, e_3, \ldots, e_k, \qquad (4)$$

*which satisfies the following properties: 1) the set of all letters which appear with odd indices before $e_{k-1}$ in (4), will be the entire alphabet (3); 2) the letters $e_i$ and $e_{i+2}$ are heterotypic, and the letters $e_k$ and $e_{k-1}$ are heterotypic and different from the letter $p$.*

We can write down the defining relations of the group $F$ over the alphabet (4). To do this, it is necessary to add a number of relations of the form $e_i = e_j$, reflecting the repetitions of letters present in (4).

Let the defining relations of the group $F$ over the alphabet (4) be

$$C_\nu = D_\nu. \qquad (5)$$

Consider an arbitrary word $A$ of the group $F$, not containing any letters $p$ or $p^{-1}$. For every word of this form, we will construct a finitely presented group, which we will denote by $F_{qA}$.

The positive alphabet of the group $F_{qA}$ will contain the letters from the alphabet (4), together with $k$ new letters $q_1, q_2, \ldots, q_k$. As defining relations of the group $F_{qA}$ we will take all relations (5), together with the following defining relations:

$$q_i q_{i+1} = q_{i+1} e_i \qquad (i = 1, 2, \ldots, k-1); \qquad (6)$$

$$e_{2j-1} q_{2j-1} = q_{2j-1} E \qquad (j = 1, 2, \ldots, \big[\tfrac{k-1}{2}\big]); \qquad (7)$$

$$e_k q_k E = q_k e_k q_k^{-1}, \qquad (8)$$

where $E$ denotes the word $p^2 A p^{-1} A p^{-1}$.

**Main Lemma.** *If $A = 1$ in the group $F$, then the corresponding group $F_{qA}$ is trivial. If the word $A$ has infinite order in the group $F$, then the group $F$ will be a subgroup of $F_{qA}$.*

The first part of this lemma is easily proved by using the relations (6), (7), (8), and Lemma 3. The proof of the second part is difficult, and represents the main difficulty of the present article.

Using the Main Lemma, the proof of Theorem 1 is not difficult.

Every word $A$ of the group $F_0$ can be put into correspondence with a finitely presented group

$$F_{qA}' = F_{qA} * F_1,$$

where $F_{qA}$ is the group construct from the word $A$, and $F_1$ is a finitely presented group with property $\alpha$, the existence of which is assumed in the statement of Theorem 1.

We have constructed a class of finitely presented groups. If $A = 1$ in $F_0$, then $A = 1$ in $F$, and – by the Main Lemma – the group $F_{qA}$ is trivial. Then $F_{qA}'$ is isomorphic to $F$, and the invariant property $\alpha$ will be possessed by $F_{qA}'$.

Suppose $A \neq 1$ in $F_0$. Then the word $A$ has infinite order in $F_0$, as $F_0$ is a torsion-free group. Consequently, the word $A$ has infinite order in $F$. By the Main Lemma, the group $F$ is a subgroup of $F'_{qA}$, and hence also of $F'_{qA}$. The group $F$ has $F'_2$ as a subgroup, and this in turn has $F_2$ as a subgroup. As the group $F_2$ by the assumptions of Theorem 1 cannot be embedded in any finitely presented group with property $\alpha$, it follows that $F'_{qA}$ does not have property $\alpha$. As the identity problem is unsolvable in the group $F_0$, there cannot exist an algorithm which decides whether or not a given group of the form $F'_{qA}$ has property $\alpha$ or not. Theorem 1 is proved.

We cannot weaken the assumption in Theorem 1, requiring the existence of a finitely presented group which does not embed in any finitely presented group with property $\alpha$, by replacing it with the existence of a finitely presented group without property $\alpha$. Indeed, there are finitely presented groups which coincide with their commutator subgroup, and there are finitely presented groups which do not coincide with their commutator subgroup. At the same time, there exists an algorithm which decides for any given finitely presented group, whether or not it coincides with its commutator subgroup. This is the algorithm solving the identity problem in commutative groups.

Notice that the same algorithm allows one to decide for every arbitrary finitely presented solvable (or arbitrary nilpotent) group whether or not it is trivial. This is even more interesting when considering the fact that the identity problem in solvable groups has not been solved.

Every invariant property $\alpha$ gives rise to the class of finitely presented groups which possess $\alpha$. This class will be called the *class of $\alpha$-groups*. A class of finitely presented groups $\alpha$ will be called *complete*, if every finitely presented group is isomorphic to a subgroup of some group from the class $\alpha$. Theorem 1 can then be reformulated in the following way.

**Theorem 1'.** *If $\alpha$ is a non-empty and incomplete class of finitely presented groups, then there does not exist an algorithm which decides for every finitely presented group whether or not it lies in the class $\alpha$.*

Requiring that $\alpha$ be incomplete is not a necessary condition, as there are infinitely many complete classes for which deciding membership is impossible. For example, the class of finitely presented groups which can be decomposed into a direct product of $k$ groups, or the class of groups which can be decomposed into a free product of $k$ groups, etc. The completeness of these classes is obvious.

The second part of this article is devoted to the proof of completeness of some classes of finitely presented groups. From what was mentioned above regarding the class of groups coinciding with their commutator subgroup, it follows that this class is complete. We note right away that every complete class contains a group with unsolvable identity problem.

As was proved in [3], there does not exist an algorithm which decides for every finitely presented group whether or not it is simple. If it is known in advance that a given group is simple, then its identity problem is very easy to solve.* This implies that the class of simple groups is incomplete in the sense defined above.

For a finitely presented group $F$ to be simple, it is necessary and sufficient that if any relation $A = B$ is added to the defining relations to $F$, then this transforms the group into the trivial group.

---

*This algorithm was communicated to the author by A. V. Kuznetsov.

Consider a finite system of words

$$L_1, L_2, \ldots, L_r, \tag{9}$$

written over the letters $x_1, x_2, \ldots, x_s, x_1^{-1}, x_2^{-1}, \ldots, x_s^{-1}$. The finitely presented group $F$ will be called *conditionally trivial with respect to the system of identical relations*

$$L_i \equiv 1 \qquad (i = 1, 2, \ldots, r), \tag{10}$$

if adding to the defining relations of $F$ any of the identical relations (10) to $F$ transforms the group $F$ into the trivial group. Analogously, we define *conditionally finite* groups, *conditionally abelian* groups, etc. If the class of $\alpha$-groups contains the trivial group, then the class of conditionally $\alpha$-groups contains the entire class of finitely presented groups.

An identical relation will be called *non-trivial*, if it does not hold in free groups.

**Theorem 2.** *For every finite system of non-trivial identities (10), the class of conditionally trivial groups with respect to this system of identities will be complete.*

We established earlier that any finitely presented group can be isomorphically embedded in a group of type $F_{qA}$. For this, it suffices that the word $A$ does not contain any letters $p$ or $p^{-1}$, and that $A$ has infinite order in the group $F$.

To prove Theorem 2, it suffices to find a word $A_0$ which has infinite order in the group $F$, and becomes trivial when any of the identical relations (10) are added.

The elements $a_1$ and $b_1$ are free generators of the subgroup $F'$ of $F$ generated by them. Consequently, the group $F'$ contains a free subgroup of every finite rank. Let $F''$ be a free subgroup with rank $s + 1$ in $F'$. Let a set of free generators for it be $E_1, E_2, \ldots, E_s, E_{s+1}$. Replace the letters $x_1, x_2, \ldots, x_s, x_1^{-1}, x_2^{-1}, \ldots, x_s^{-1}$ by the corresponding words $E_1, E_2, \ldots, E_s, E_1^{-1}, E_2^{-1}, \ldots, E_s^{-1}$ in the words (9). We denote the resulting words by $L_1', L_2', \ldots, L_r'$.

As every identical relation in (10) is by assumption non-trivial, none of the words $L_i'$ $(i = 1, 2, \ldots, r)$ can be equal to 1 in the group $F''$. The word

$$A_0 = [\cdots[[E_{s+1}L_1', L_2'], L_3'], \ldots, L_r'] \tag{11}$$

is therefore not equal to 1 in $F''$, as in a free group two elements do not commute if they are both not equal to the identity element and not equal to one another. As the group $F''$ is free, the word $A_0$ has infinite order in $F''$, and consequently, also in the group $F$. This word $A_0$ becomes equal to the identity element if any of the identical relations (11) are added. Theorem 2 is proved.

**Theorem 3.** *The class of finitely presented groups, defined by finite systems of conjugate generators, is complete.*

Since every generator of the group $F_{qA}$ is conjugate to the word $A$, to prove Theorem 3 it suffices to take $A = a_1$.

S. I. ADIAN

# ON ALGORITHMIC PROBLEMS IN EFFECTIVELY COMPLETE CLASSES OF GROUPS

*(Communicated by academician I. M. Vinogradov, 16 May 1958)*

The purpose of the present article is to further generalise some theorems on the un-solvability of algorithmic problems in the theory of groups [1, 2, 3].* We will use the notation from the articles [2, 3]. In the first part we will give some remarks on the article [3], allowing us to completely restore the proof of the main lemma of that article.** The following remarks will yield the proof of a somewhat stronger assertion, from which the indicated lemma follows as a corollary. Specifically, in the first part we will prove the following lemma about the group $F_{qA}$.

**Lemma B.** *Suppose the word $A$ of the group $F$ does not contain letters $p^\sigma$. If $A = 1$ in the group $F$, then $F_{qA}$ is the trivial group. If $A \neq 1$, then $F$ is a subgroup of $F_{qA}$.*

In the second part, using Lemma B, we will prove some theorems about the un-solvability of algorithmic problems in effectively complete classes of finitely presented groups (abbreviated as f.p. groups).

**I. Some remarks on the article [3].** The analogy between the defining relations of the groups $\mathfrak{A}_{q,A,B}$, constructed in [2], and the defining relations of the groups $F_{qA}$, constructed in [3], is obvious. In the defining relations of the group $F_{qA}$ the role of the word $AB^{-1} \equiv pXpX^+Y^{+^{-1}}p^{-1}Y^{-1}p^{-1}$ is played by $E \equiv p^2Ap^{-1}Ap^{-1}$. This analogy runs deep.

The first two sections of Chapter II in [2] are devoted to proving Lemma 4, §2, Chapter II. As part of the proof, a series of auxiliary assertions are used, which appear in Chapter I. The statements contained in §1, Chapter I are all of a general nature, and apply equally well to any group. The statements contained in §2, Chapter I of [2] related, instead, to concrete properties of the group $\mathfrak{A}_p$, and are connected with the peculiarities of the construction of this particular group. Nevertheless, analogous statements remain true for the group $F = F_0 * F_2' * F_3$ constructed in [3].

For every statement, contained in §2, Chapter I, of [3] and used in first to sections of Chapter II, we will indicate the corresponding statement about the group $F$, which will be used in the proof of Lemma B.

1) Lemma 1, §2, Chapter I (p. 22). This corresponds to the following clear lemma.

**Lemma 1\*.** *Suppose $D$ is a word in the group $F$, and that it does not contain any letters $p^\sigma$. If $D = 1$ in $F$, then the projections of $D$ onto the alphabets of the group $F_0$ and $F_2'$, respectively, will be equal to 1 in the group $F$.*

2) Lemma 3, §2, Chapter I (p. 24). This lemma corresponds to Lemma 3 of [3] (see p. 77). This lemma is also obvious.

3) Lemma 4, §2, Chapter I (p. 25). This lemma will be replaced by Lemma 2 of the article [3] (p. 77). The proof of Lemma 2 of [3] is based on Lemma 1 of the same article, which we will prove first.

---

*In [3], on p. 77 after formula (8) the following should be added: "where $E$ is the word $p^2Ap^{-1}Ap^{-1}$".

**Theorem 1 of [3] has also been proved by Rabin in [5].

*Proof of Lemma 1.* Suppose the f.p. group $F_2$ has generators $a_1, a_2, \ldots, a_n$ and defining relations $A_i = 1$. As our group $F_2'$, we may take the group with generators

$$b_1, b_2, \ldots, b_n, c \tag{1}$$

and defining relations $B_i = 1$, where the word $B_i$ is obtained from the word $A_i$ by replacing each letter $a_j^\sigma$, $\sigma = \pm 1$, by the word $(cb_j)^\sigma$. The subgroup of $F_2'$ generated by the elements $cb_j$ will be isomorphic to the group $F_2$. It is easy to see that the generators (1) have infinite order in $F_2'$.        □

*Proof of Lemma 2.* Every generator of the group $F$ has infinite order in $F$ (for the generators of the group $F_0$ this follows from Theorem 2, p. 69, of [4]; for the generators of $F_2'$ it follows from Lemma 1 which we just proved; and the letter $p$ does not appear in any defining relation of $F$). Any two heterotypic letters of the alphabet of $F$ will be generators from different free factors of the group $F$. Consequently, they will be free generators of the subgroup of $F$ generated by them. Lemma 2 is proved.        □

4) Lemma 2′, §2, Chapter I (p. 25). This lemma is a consequence of the fact that $a_n$ and $a_{n-1}$ are free generators of the subgroup of $\mathfrak{A}_P$ generated by them. By Lemma 2 and Lemma 3 in [3], $e_k$ and $e_{k-1}$ are free generators of the subgroup of $F$ generated by them. Hence, it follows that the statement of Lemma 2′ holds true for $e_k$ and $e_{k-1}$.

5) Lemma 7, §2, Chapter I (p. 33) is replaced by the following statement.

**Lemma 7\*.** *Suppose $A, C, D$ are arbitrary words of the group $F$, not containing the letter $p^\sigma$, and with $A \neq 1$ in $F$. If in this case we have*

$$p^2 A p^{-1} A p^{-1} D p A^{-1} p A^{-1} p^{-2} C = 1 \quad in \ F,$$

*then in $F$ the projections of the words $C$ and $D$ onto the alphabets of $F_0$ and $F_2'$ are equal to 1.*

This lemma follows from the fact that $F_0$, $F_2'$, and $F_3 = \{p\}$ are free factors of $F$.

6) Lemma 6, §2, Chapter I and its corollary (pp. 26 and 33) will be replaced by the following lemma.

**Lemma 6\*.** *Suppose the word $A$ of the group $F$ does not contain the letter $p$ and is not equal to 1 in $F$. Then a necessary and sufficient condition for the word*

$$Q \equiv e^{k_1}(p^2 A p^{-1} A p^{-1})^{s_1} e^{k_2} \cdots e^{k_\mu}(p^2 A p^{-1} A p^{-1})^{s_\mu} e^{k_{\mu+1}},$$

*where $e$ is an arbitrary letter in the alphabet of $F$, to be equal to 1 in $F$, is for the corresponding sequence of numbers $k_1, s_1, k_2, \ldots, k_\mu, s_\mu, k_{\mu+1}$ to be degenerate.*

Lemma 6\* and its corollaries are proved in the same way that, in Chapter I of [2], Lemma 6 and its corollaries are proved.

Thus, for each of the auxiliary statements about the group $\mathfrak{A}_P$ ((1)–(6)), we have produced a corresponding statement about the group $F$. Now the proof of Lemma B is a verbatim repetition of the arguments of the first two sections of Chapter II in [2].

**II. Algorithmic problems in effectively complete classes of groups.** By the problem of recognising the group property $\alpha$ for groups from the class $K$, we will mean the problem of finding an algorithm which for every group of the class $K$ decides the question of whether the group has the property $\alpha$ or not. In [1, 2, 3] we studied a particular type of such problems, when the class $K$ consisted of all f.p. groups.

A class of f.p. groups will be said to be *effectively complete* if there exists an algorithm $\mathfrak{N}$, which converts the presentation of an arbitrary f.p. group $F$ into a presentation

of some group from the class $K$ with a subgroup isomorphic to $F$. We will call such an algorithm $\mathfrak{N}$ an *algorithm which effectivizes the completeness of $K$*. We will say that the algorithm $\mathfrak{N}$ *fixes the group $F_0$*, if it converts every presentation for $F_0$ into a presentation for the same group $F_0$.

**Condition T.** We will say that the pair $(K, \alpha)$, where $K$ is some class of f.p. groups and $\alpha$ is some invariant property of groups, satisfies condition T, if one can find some f.p. group $F_1$ in the class $K$ which has property $\alpha$, and some f.p. group $F_2$ which cannot be embedded into any f.p. group with property $\alpha$.

**Theorem 1.** *Let $K$ be an effectively complete class of f.p. groups. If the invariant property $\alpha$ is such that the pair $(K, \alpha)$ satisfies condition T, and the algorithm $\mathfrak{N}$, which effectivizes the completeness of the class $K$, fixes some group $F$ with property $\alpha$, then the recognition problem for the property $\alpha$ in the class $K$ is unsolvable.*

*Proof.* Consider the group $F_2$, which by virtue of condition T does not embed in any f.p. group with property $\alpha$. Using Lemma 1 of [3], we can embed $F_2$ in a f.p. group $F_2'$ in which every generator has infinite order. Next, using the group $F = F_0 * F_2' * \{p\}$, where $F_0$ is some f.p. group with unsolvable identity problem, construct, just as in [3]*, the system of groups $F_{qA}$. Above, we proved that Lemma B is true for the groups $F_{qA}$. Let $F_A = F_{qA} * F_1$, where $F_1$ is a group in the class $K$ possessing the property $\alpha$. Consider the system of groups $G_A = \mathfrak{N}(F_A)$, where $A$ is an arbitrary word in the group $F_0$. As the algorithm $\mathfrak{N}$ effectivizes the completeness of $K$, every group $G_A$ belongs to the class $K$. Theorem 1 will be proved, if we can prove that there does not exist an algorithm which decides for every group $G_A$ whether or not it has property $\alpha$. For this, it suffices to prove that $G_A$ possesses property $\alpha$ if and only if $A = 1$ in $F_0$.

Suppose $A = 1$ in $F_0$. Then $A = 1$ in $F$. By Lemma B, in this case the group $F_{qA}$ is the trivial group. Consequently, $F_A$ is isomorphic to the group $F_1$, and hence $G_A$ is also isomorphic to $F_1$. Hence, the invariant property $\alpha$ is possessed by $G_A$.

Suppose $A \neq 1$ in $F_0$. Then $A \neq 1$ in $F$. In this case, by Lemma B the group $F$ will be a subgroup of the group $F_{qA}$. Hence we find $F_2 \subset F_2' \subset F \subset F_{qA} \subset F_A \subset \mathfrak{N}(F_A) = G_A$, i.e. the group $F_2$ embeds in $G_A$. Hence, by virtue of condition T, the group $G_A$ cannot possess the property $\alpha$. The theorem is proved.  □

The following is an interesting particular case of Theorem 1.

**Theorem 2.** *Let $K$ be an effectively complete class of f.p. groups, with an algorithm effectivizing its completeness fixing every group in this class. If the invariant property $\alpha$ is such that the pair $(K, \alpha)$ satisfies condition T, then the recognition problem for $\alpha$ for groups in $K$ is unsolvable.*

It is easy to see that Theorem 1 of [3] is a special case of this theorem.

A class of f.p. groups $K$ will be said to be *recursive*, if there exists an algorithm which decides for every f.p. group whether or not this group lies in $K$.

**Theorem 3.** *Suppose $K$ is a recursive and effectively closed class of f.p. groups, and $\alpha$ is some invariant property of groups. If the pair $(K, \alpha)$ satisfies condition T, then the recognition problem for the property $\alpha$ for groups in the class $K$ is unsolvable.*

*Proof.* Let $\mathfrak{N}$ be an algorithm effectivizing the completeness of $K$, and $\mathfrak{M}$ an algorithm which decides for every f.p. group whether or not it lies in $K$ or not. Using these two algorithms, it is easy to construct an algorithm $\mathfrak{N}^*$ which effectivizes the completeness of $K$ and which leaves every group in $K$ fixed, hence proving Theorem 3.  □

---

*See the footnote on p. 80.

As a non-trivial example of a recursive and effectively-closed class of f.p. groups, we have the class of f.p. groups which coincide with their commutator subgroup.

Every system of non-trivial identical relations

$$L_i \equiv 1 \qquad (i = 1, 2, \ldots, k), \tag{2}$$

written using the symbols $x_1, x_2, \ldots, x_r$, defines the class consisting of all groups which are conditionally trivial with respect to this system (see [3], p. 79).

**Theorem 4.** *Let $K_1$ be a class of conditionally trivial groups with respect to the system (2) of non-trivial identical relations. If the invariant group property $\alpha$ is such that the pair $(K_1, \alpha)$ satisfies condition* T, *then the recognition problem for property $\alpha$ for groups in $K_1$ is unsolvable.*

*Proof.* We first embed the group $F_2$, by using Lemma 1 from [3], into a group $F_2'$, in which every generator has infinite order. We will take $F_2'$ to be such that it also contains a free subgroup on two generators $b_1$ and $b_2$. Let $B_1, B_2, \ldots, B_r$ be free generators of the subgroup they generate of the free group on $b_1$ and $b_2$. Replace everywhere in the words $L_i$ the generators $x_1, x_2, \ldots, x_r$ by the words $B_1, B_2, \ldots, B_r$. The resulting words of the group $\{b_1, b_2\}$ will be denoted by $Y_i$ $(i = 1, 2, \ldots, k)$. As every identical relation (2) is non-trivial by assumption, none of the words $Y_i$ are equal to 1 in the group $\{b_1, b_2\}$. Thus $Y_i \neq 1$ also in $F_2'$. Let $C$ be an arbitrary word from the group $F_0$. We set

$$A = [\cdots [[C, Y_1], Y_2] \cdots Y_k], \tag{3}$$

where the square brackets denote taking the commutator. Then $A$ is a word in the group $F = F_0 * F_2' * \{p\}$.

It is easy to see that $A = 1$ in $F$ if and only if $C = 1$ in $F_0$. Hence there does not exist an algorithm which decides for every word $A$ of type (3), whether or not it is equal to 1 in $F$. We note also that adding any of the identical relations (2) to the defining relations of $F$ will transform the corresponding word $Y_i$, and hence also the word $A$, into 1.

Consider now the system of groups $F_{qA}$, just as it was done in [3], which correspond to the words (3) for all possible $C \in F_0$. By virtue of Lemma B, each of these groups will be conditionally trivial with respect to the system (2), i.e. all these groups belong to the class $K_1$. Clearly, every group $F_A = F_{qA} * F_1$ belongs to $K_1$, where $F_1$ is a group which by condition T belongs to $K_1$ and possesses property $\alpha$. Theorem 4 will be proved, if we can prove that there does not exist an algorithm, which for every group $F_A$ decides whether or not it has property $\alpha$. To prove this, it suffices to show that the group $F_A$ will possess property $\alpha$ if and only if $A = 1$ in $F$. This is proved in the same manner as in Theorem 1.  □

A. A. MARKOV

# IMPOSSIBILITY OF SOME ALGORITHMS IN THE THEORY OF ASSOCIATIVE SYSTEMS

*(Communicated by academician I. M. Vinogradov, 4 January 1951)*

**1.** When considering associative systems, defined by a finite system of defining relations on some given alphabet $A$, the following problems arise naturally.

*The isomorphism problem.* Give an algorithm, by which one can decide for any two finite systems of relations on $A$ whether or not the associative systems defined by them are isomorphic.

*The triviality problem.* Give an algorithm, by which one can decide for any finite system of relations on $A$ whether or not the associative system defined by it is the system with a single element.

*The group recognition problem.* Give an algorithm, by which one can decide for any finite system of relations on $A$ whether or not the associative system defined by it is a group.

*The group embeddability problem.* Give an algorithm, by which one can decide for any finite system of relations on $A$ whether or not the associative system defined by it is embeddable in a group.

*The semigroup recognition problem.* This is formulated in the same way as the group recognition problem, substituting "semigroup" for "group".

*The meta-identity problem.* Give an algorithm, by which one can decide, for any finite system of relations on $A$, whether or not the identity problem for the associative system defined by it is solvable.

A solution to the isomorphism problem would, of course, also give a solution to the triviality problem. The concepts of groups, embeddings into groups, and semigroups, are all well-known [3]. Concerning the meta-identity problem, it is known [1, 2, 4] that there exist associative systems of the type considered for which the identity problem is unsolvable. Therefore, it is natural to look for a general effective method for recognising such "bad" associative systems from their defining systems of relations.

**2.** We managed to establish the following result.

*For every alphabet $A$ containing more than three letters, each of the problems specified above are unsolvable: the desired algorithm is in each case impossible.*

The concept of algorithm has the same sense here as in our previous notes [1, 2]. We now briefly outline how to prove the stated result.

**3.** We will start with the system of relations $\mathfrak{C}_5$ on the alphabet $A_5$, constructed in [2], which defines an associative system with unsolvable identity problem. This associative

system is not a semigroup, as the relation $h \longleftrightarrow f$ clearly does not hold in $\mathfrak{C}_5$, even though $dh \longleftrightarrow df$ is one of the relations of $\mathfrak{C}_5$.

We match the letters $a, b, \ldots, m$ of the alphabet $A_5$ to the words

$$aba, abba, \ldots, abbbbbbbbbbbba$$

over the alphabet $A_0 = \{a, b\}$, and replace in $\mathfrak{C}_5$ each letter by its corresponding word, obtaining a system of relations $\mathfrak{D}$ over the alphabet $A_0$. In view of the fact that every left-hand side and every right-hand side of every defining relation in $\mathfrak{C}_5$ is non-empty, it easily follows that the associative system defined by the system $\mathfrak{D}$ has unsolvable identity problem, and that this associative system is not a semigroup.[1]

**4.** Let now $G$ and $H$ be two arbitrary words over $A_0$. We construct from these words a new system of relations $\mathfrak{D}_{G,H}$ over the alphabet $B = \{a, b, c, d\}$, by adding to $\mathfrak{D}$ the new relations $cGd \longleftrightarrow 0$ and $\xi cHd \longleftrightarrow cHd$ ($\xi = a, b, c, d$).

We can establish the following lemmas.

**Lemma 1.** *If the relation*

$$(1) \qquad\qquad\qquad\qquad G \longleftrightarrow H$$

*follows from $\mathfrak{D}$, then the associative system defined by the system of relations $\mathfrak{D}_{G,H}$ consists of a single element.*

**Lemma 2.** *If the relation (1) does not follow from $\mathfrak{D}$, then for any two words $P$ and $Q$ in the alphabet $A_0$ the relation $P \longleftrightarrow Q$ follows from $\mathfrak{D}_{G,H}$ if and only if it follows from $\mathfrak{D}$.*

From Lemma 2, it follows that in the case that (1) does not follow from $\mathfrak{D}$, then the identity problem is unsolvable for the associative system defined by $\mathfrak{D}_{G,H}$, and that in this case this associative system is not a semigroup.

By putting it all together, we see that a solution for the alphabet $B$ of any of the five problems: triviality, group recognition, group embeddability, semigroup recognition, or meta-identity problem; would entail a solution to the identity problem for the associative system defined by $\mathfrak{D}$ over the alphabet $A_0$. Such a solution is, however, impossible. Consequently, these problems are unsolvable for the four-letter alphabet $B$. By virtue of the unsolvability of the triviality problem for $B$, the isomorphism problem is also unsolvable for this alphabet. Changing to an arbitrary alphabet with at least four letters is trivial.

**Leningrad Department of the**
**Steklov Institute of Mathematics**
**of the USSR Academy of Sciences**

---

[1]This method of reduction to an alphabet with two letters was invented by Post [5].

A. A. MARKOV

# IMPOSSIBILITY OF ALGORITHMS FOR RECOGNISING SOME PROPERTIES OF ASSOCIATIVE SYSTEMS

*(Communicated by academician I. M. Vinogradov, 23 February 1951)*

**1.** When our previous article [3] had already appeared in print, we were able to obtain a very general theorem on unsolvability, from which the published results can be directly obtained as a special case. The present note is mainly devoted to establishing this theorem.

**2.** We will call a *K-system* any associative system which is defined by a finite system of relations over some alphabet [1].

**3.** We will say that an associative system is *trivial* if it consists only of its identity element. A trivial associative system is a K-system.

**4.** Let $n$ be a natural number. We will say that a K-system $S$ is *definable* by $n$ letters in the following two cases: first, when $n > 0$ and $S$ is defined by a finite system of relations over an alphabet with $n$ letters; and second, when $n = 0$ and $S$ is a trivial associative system.

**5.** We will say that an associative system $S$ is *embedded* in the associative system $T$, if $T$ has a subsystem[1] isomorphic to $S$.

**6.** We will be interested in various properties which can be possessed by K-systems (in short – "properties of K-systems") such as, for example, being trivial, being finite, being embeddable in a group, being a group, or being a semigroup.

We will say that a property $\mathfrak{P}$ of K-systems is *invariant*, if every K-system isomorphic to some K-system $\mathfrak{R}$ with property $\mathfrak{P}$, itself has property $\mathfrak{P}$. All five aforementioned properties of K-systems are clearly invariant.

**7.** We will say that a property $\mathfrak{P}$ of K-systems is *hereditary*, if every K-system which is embeddable in some K-system with property $\mathfrak{P}$ also has the property. Every hereditary property of K-systems is obviously invariant. Being trivial, finite, embeddable in a group, and the property of being a semigroup are all hereditary, but the property of being a group is not hereditary.

**8.** Every finite system of relations over any alphabet defines a K-system uniquely up to isomorphism. If $\mathfrak{P}$ is an invariant property of K-systems, then given any finite system of relations $\mathfrak{R}$ over an alphabet $A$ we may ask the question, whether or not this property is satisfied by the K-system defined by the system $\mathfrak{R}$. Naturally, for any fixed alphabet $A$ the following question is natural to ask.

Give an algorithm which, for every given finite system of relations over $A$, decides whether or not the K-system it defines has property $\mathfrak{P}$.

This problem will be called the *recognition problem for property $\mathfrak{P}$ for the alphabet A*. The term "algorithm" is used in the same sense as in our previous work [1, 2, 3].

---

[1] Although all of the associative systems considered by us have an identity element [1], it may well happen that the identity element of a subsystem of $T$ does not coincide with the identity element of $T$. This must be kept in mind in what follows.

**9.** The following general impossibility theorem holds.

*Let $\mathfrak{P}$ be an invariant property of $K$-systems. If there exists some $K$-system with this property, and some $K$-system which does not embed in any $K$-system with this property, then there exists an alphabet for which the recognition problem for property $\mathfrak{P}$ is unsolvable (in the sense that it is impossible for the desired algorithm to exist). If at the same time there is a $K$-system with property $\mathfrak{P}$, definable by $n$ letters, then for every alphabet with $n + 4$ or more letters, the recognition problem for property $\mathfrak{P}$ is unsolvable.*

We will now give a sketch of the proof of this theorem.

**10.** We first introduce the following terminology.

We will say that a relation is *regular*, if both its sides are non-empty. We will say that a system of relations is *regular*, if all its relations are regular. We will say that a $K$-system is *regularly definable*, if there exists a regular system of relations which defines it.

**11.** It is easy to prove the following lemmas.

**Lemma 1.** *Every $K$-system can be embedded in some regularly definable $K$-system.*

**Lemma 2.** *Every regularly definable $K$-system can be embedded in some $K$-system definable by two letters.*

Hence, every $K$-system is embeddable in some $K$-system definable by two letters.

**12.** Let now $S_0$ be a $K$-system which does not embed in any $K$-system with $\mathfrak{P}$. Such a $K$-system exists, by the assumptions of the theorem. Let this system be defined by a system of relations $\mathfrak{R}_0$ with alphabet $B_0$.

As is known, there exists a $K$-system with unsolvable identity problem [1, 2, 4]. Let $S_1$ be such a $K$-system and let $S_1$ be defined by a system of relations $\mathfrak{R}_1$ with alphabet $B_1$. This alphabet will be chosen to have no letters in common with $B_0$, which can obviously be done.

The union of the alphabets $B_0$ and $B_1$ will be denoted by $B_2$; the union of the systems $\mathfrak{R}_0$ and $\mathfrak{R}_1$ will be denoted by $\mathfrak{R}_2$. The system $\mathfrak{R}_2$ is a system of relations over the alphabet $B_2$ and, as such, defines some $K$-system $S_2$ (being nothing else than the free product of the systems $S_0$ and $S_1$).

It is not difficult to see that $S_0$ embeds in $S_2$. Hence $S_2$ also does not embed in any $K$-system with property $\mathfrak{P}$. It is not difficult to see that the identity problem is unsolvable for $S_2$, just as it is in $S_1$.

According to what we have proved, $S_2$ embeds in some $K$-system $S_3$ which can be defined by some system of relations $\mathfrak{R}_3$ over the alphabet $\{a, b\}$. The $K$-system $S_3$, just as $S_2$, does not embed in any $K$-system with property $\mathfrak{P}$. For $S_3$, just as for $S_2$, the identity problem is unsolvable.

**13.** Suppose further that $S_4$ is a $K$-system with property $\mathfrak{P}$, definable by $n$ letters. Such a $K$-system exists, by the assumptions of the theorem. Suppose that if $n > 0$, the $K$-system $S_4$ is defined by a system of relations $\mathfrak{R}_4$ with alphabet $B_4$. We will choose this alphabet so that it does not include the letters $a, b, c, d$, which is obviously possible.

For $n > 0$ we denote by $B_5$ the alphabet obtained from joining to $B_4$ the letters $a, b, c, d$. For $n = 0$, we denote by $B_4$ the alphabet $\{a, b, c, d\}$. The alphabet $B_5$ always consists of $(n + 4)$ letters.

Let $G$ and $H$ be arbitrary words over the alphabet $\{a, b\}$. We will construct a system of relations $\mathfrak{R}_{G,H}$ over the alphabet $B_5$ in the following manner. When $n > 0$, we define it as the union of the systems $\mathfrak{R}_3$, $\mathfrak{R}_4$, and the system of five relations

(1)     $cGd \longleftrightarrow 0; \quad \xi cHd \longleftrightarrow cHd \qquad (\xi = a, b, c, d);$

and when $n > 0$ we take the union of $\mathfrak{R}_3$ with the system of relations (1).

Denote by $S_{G,H}$ the $K$-system defined by the system $\mathfrak{R}_{G,H}$ with alphabet $B_5$.
We can now establish the following lemmas.

**Lemma 3.** *If the relation*

$$(2) \qquad\qquad\qquad\qquad G \longleftrightarrow H$$

*follows from $\mathfrak{R}_3$, then $S_{G,H}$ is isomorphic to $S_4$.*

**Lemma 4.** *If the relation (2) does not follow from $\mathfrak{R}_3$, then for any two words $P$ and $Q$ over the alphabet $\{a, b\}$, the relation $P \longleftrightarrow Q$ holds in $\mathfrak{R}_{G,H}$ if and only if it holds in $\mathfrak{R}_3$.*

From Lemma 4 it follows that in the case that (2) does not follow from $\mathfrak{R}_3$, we have that $S_3$ embeds in $S_{G,H}$. Taking into account the properties of the $K$-systems $S_3$ and $S_4$, we conclude from Lemma 3 that $S_{G,H}$ has property $\mathfrak{P}$, if (2) follows from $\mathfrak{R}_3$, and that it does not have the property if (2) does not follow from $\mathfrak{R}_3$. Hence a solution to the recognition problem for property $\mathfrak{P}$ with alphabet $B_5$ would entail a solution to the identity problem in $S_3$, which, however, is impossible. Consequently, the recognition problem for property $\mathfrak{P}$ for the $(n + 4)$-lettered alphabet $B_5$ is impossible. Switching to arbitrary alphabets with at least $(n + 4)$ letters is obvious.

**14.** When applied to hereditary properties, our theorem can be simplified. Indeed, if the property $\mathfrak{P}$ is hereditary, then any $K$-system without this property cannot be embedded in any $K$-system with this property. On the other hand, if there is in this case a $K$-system with property $\mathfrak{P}$, then this implies that the trivial $K$-system also has this property, as this system embeds in every $K$-system. Taking this all into account, we find, as a corollary of the proved theorem, the following result.

*Let $\mathfrak{P}$ be a hereditary property of $K$-systems. If there exists some $K$-system with this property, and some $K$-system without this property, then for every alphabet with 4 or more letters, the recognition problem for property $\mathfrak{P}$ is unsolvable.*

**15.** By applying our result to the hereditary properties of being trivial, being embeddable in a group, and being a semigroup, we directly recover the results of our previous article regarding these properties. Its application to the hereditary property of being finite gives, *that for every alphabet with $\geq 4$ letters, the recognition problem is unsolvable for finite $K$-systems.*

Furthermore, taking into account the fact that the property of having solvable word problem is a hereditary property of $K$-systems, we find the result from our previous article regarding the meta-identity problem.

**16.** In order to recover the result on the recognition problem for groups for $K$-systems from our previous article, we must exhibit a $K$-system which does not embed in any group. Such a $K$-system can be, for example, the $K$-system defined by the relation $aa \longleftrightarrow a$ over the single-letter alphabet $\{a\}$. As, on the other hand, the trivial $K$-system is a group, we can apply our general theorem for the property of being a group with $n = 0$, which yields the corresponding result from our previous article.

**17.** We now move onto considering a wider class of associative system, for which we are able to find a strong theorem on unsolvability.

Let $\mathfrak{R}$ be a finite system of relations over the alphabet $A$, and let $B_1, \ldots, B_k$ be words over this alphabet. We will consider the $K$-system $S$ defined by the system of relations $\mathfrak{R}$. The elements of $S$ which can be expressed as a product $B_{i_1} \cdots B_{i_m}$ of words $B_i$ ($m \geq 0, 1 \leq i_\lambda \leq k$) will form an associative subsystem of the system $S$. This subsystem will be called the *associative system defined by the alphabet $A$, the system of relations $\mathfrak{R}$, and the words $B_1, \ldots, B_k$.* Any associative system, defined in this manner,

will be called a $\Pi$-*system*. We will say that a $\Pi$-system, defined in this way over the alphabet $A$, will be said to be *definable over $A$*. It can easily proved that every $\Pi$-system is definable over an alphabet with two letters.

**18.** Every $K$-system is clearly a $\Pi$-system. The concept of invariant property can be extended in the obvious way to $\Pi$-systems. For $\Pi$-systems, definable over the alphabet $A$, we can define the *recognition problem for property* $\mathfrak{P}$ in the following way.

Let $\mathfrak{P}$ be an invariant property of $\Pi$-systems. Give an algorithm which, for every finite system of relations over $A$ and any finite system of words in $A$, decides whether or not the $\Pi$-system defined by them has property $\mathfrak{P}$.

**19.** Using a method analogous to the above, we may prove the following theorem on unsolvability.

*Let $\mathfrak{P}$ be an invariant property of $\Pi$-systems. If there exists some $\Pi$-system with this property, and some $\Pi$-system without it, then for every alphabet containing more than one letter, the recognition problem for property $\mathfrak{P}$ is unsolvable: the desired algorithm is impossible.*

**Leningrad Department of the**
**Steklov Institute of Mathematics**
**of the USSR Academy of Sciences**




\end{document}